\documentclass[12pt]{amsart} 

\usepackage{amsthm,amssymb}
\usepackage{fullpage}
\usepackage[utf8]{inputenc}
\usepackage{graphicx} %%% Do not specify dvips or pdftex as an option %% Do not use psfrag.
\usepackage{subcaption}
\usepackage{amssymb}
\usepackage{xcolor}
\usepackage{float}

\newcommand{\amph}{\operatorname{Amph}}
\newcommand{\embed}{L}

\newcommand{\stw}{\mathbb{S}^2}
\newcommand{\sth}{\mathbb{S}^3}

\newcommand{\rth}{\mathbb{R}^3}

\newtheorem{proposition}{Proposition}
\newtheorem{theorem}{Theorem}
\newtheorem{lemma}{Lemma}
\newtheorem{corollary}{Corollary}
\newtheorem{question}{Question}
\newtheorem{problem}{Problem}
\newtheorem{remark}{Remark}

\begin{document}
\title[Self-dual maps II: links and symmetry]{Self-dual maps II: links and symmetry} 

\thanks{$^1$ Partially supported by CONACyT 166306 and PAPIIT-UNAM IN112614}
\thanks{$^2$ Partially supported by grant PICS07848 and INSMI-CNRS}
\author[Luis Montejano]{Luis Montejano$^1$}
\address{Instituto de Matem\'aticas, Universidad Nacional A. de M\'exico at Quer\'etaro
Quer\'etaro, M\'exico, CP. 07360}
\email{luis@im.unam.mx}
\author[Jorge L. Ram\'irez Alfons\'in]{Jorge L. Ram\'irez Alfons\'in$^2$}
\address{
%IMAG, Univ.\ Montpellier, CNRS, Montpellier, France and 
UMI2924 - J.-C. Yoccoz, CNRS-IMPA, Brazil and IMAG, Univ.\ Montpellier, France }
\email{jorge.ramirez-alfonsin@umontpellier.fr}
\author[Ivan Rasskin]{Ivan Rasskin}
\address{IMAG, Univ.\ Montpellier, CNRS, Montpellier, France}
\email{ivan.rasskin@umontpellier.fr}



\begin{abstract}  In this paper, we investigate representations of links that are either {\em centrally symmetric} in $\mathbb{R}^3$ or {\em antipodally symmetric} in $\mathbb{S}^3$.  By using the notions of {\em antipodally self-dual} and {\em antipodally symmetric} maps, introduced and studied by the authors in \cite{MRAR1}, we are able to present sufficient combinatorial conditions for a link $L$ to admit such representations. The latter naturally arises sufficient conditions for $L$ to be {\em amphichiral}.  

We also introduce another (closely related) method yielding again to sufficient conditions for $L$ to be amphichiral. We finally prove that a link $L$, associated to a map $G$, is amphichiral if the {\em self-dual pairing} of $G$ is not one of 6 specific ones among the classification of the 24 {\em self-dual pairing} $Cor(G) \rhd Aut(G)$.
\end{abstract}

\maketitle
\tableofcontents
\newpage

\section{Introduction}
Finding symmetrical diagrams for a link $L$ can be a challenging task. In this paper, we are interesting in
questions concerning the {\em symmetry} and the {\em amphichirality} of $L$. 
\smallskip

We focus our attention to study the existence of both {\em centrally symmetric} embeddings in $\mathbb{R}^3$ and  {\em antipodally symmetric} embeddings in $\mathbb{S}^3$ of $L$. It turns out that the notions of {\em antipodally self-dual map} and {\em antipodally symmetric map}, introduced and  studied by the authors in \cite{MRAR1}, are helpful combinatorial tools for the above mentioned embeddings. These able us to come up with a number of contributions on the amphichirality of links.
\smallskip

The paper is organized as follows. In the next two sections we give a brief overview of some basic notions and definitions on knots and maps needed for the rest of the paper. 
\smallskip
 
In Section \ref{sec;R3sym}, we present combinatorial conditions for $L$ to admit a {\em centrally symmetric} embedding in $\mathbb{R}^3$ when the {\em Tait graph}, associated  to $L$, is antipodally self-dual (Theorem \ref{thm:antipodalR3}). 
As a straightforward consequence, we obtain sufficient combinatorial conditions for a link to be {\em amphichiral} (Corollary \ref{cor:neg-amphi}). The latter allows us to present infinite families of amphichiral {\em alternating} links (Corollaries \ref{cor:antipodal}, \ref{cor;w}, \ref{cor;Turk'shead}, \ref{cor;e}, \ref{cor;pan} and \ref{cor;sum}). For instance, we show the amphichirality of the {\em $(3,n)$-Turk's head} links for any integer $n\ge 1$. We also present conditions on $L$ to admit a {\em antipodally symmetric} embedding in $\sth$ when the {\em Tait graph}, associated  to $L$, is antipodally self-dual (Theorem \ref{thm:antipodalS3}). This allows us to present infinite families of links antipodally symmetric in $\sth$ (Corollaries \ref{cor:antipodal3sph}, \ref{cor;wantipol3}, \ref{cor;Torus} and \ref{cor;sumS3}).  For instance, we show that the torus knot $T(2,n)$ admits such embedding for any $n\ge 2$.
\smallskip

The approach given in Section \ref{sec;R3sym} does not detect amphichirality when the Tait graphs is not antipodally self-dual. For instance, it is not able to determine the amphichirality of the {\em Figure-eight} knot (whose associated Tait graph is not antipodally self-dual), see Remark \ref{rem:approach}. 
In Section \ref{sec;gamma}, we propose a different method, by introduce and investigate the notion of 
{\em $\Gamma$-curve} in an incidence graph $I(G)$ where $G$ is a Tait's graph associated to a link $L$. This tool allows us to give sufficient combinatorial conditions on $I(G)$ for $L$ to be amphichiral. The latter is done via sphere isometries (Theorem \ref{theo:amphi}). This approach shows the amphichirality of new families of links (in particular the Figure-eight knot) and lead to an easy way to construct infinite many amphichiral links.  We also discuss a closely related method in connection to {\em rigid} and {\em invertible} knots. We finally introduce the {\em amphichiral number} of a link $L$ (the minimal number of crossing switches to become $L$ amphichiral) and present a first rough upper bound for some particular families (Proposition \ref{pro:amphinumber}).
\smallskip

In Section \ref{sec;isomap}, we focus our attention to self-dual pairings. After recalling basic notions and definitions, we study some properties of {\em equivalent} maps and we give necessary conditions for the associated link to be amphichiral (Theorem \ref{thm:self-dual-amphi}). We then consider  the groups $Cor(G)$ (generated by the set of all self-dualities of $G$) and $Aut(G)$ (generated by the set of all automorphisms of $G$). It is known that any map of $G$  belong to one of the 24 {\em self-dual pairings} $Cor(G) \rhd Aut(G)$ classification. We show that the link $L$, associated to a self-dual map $G$ which its self-pairing is other than 6 specific ones, is amphichiral (Theorem \ref{theo;amphi1}). 
%We also present a result (Theorem \ref{thm:symgroup}) on the same line of the {\em symmetry group} of a knot introduced by Gr\"unbaum and Sherphard \cite{GS}.
\smallskip

We finally end with some concluding remarks.

%%%%%%%%%%%%%%%%%%%%%%%%%%%%%%%%%%%%%%%%%%%
\section{Knot : preliminaries}\label{sec:knot}

We review some knot notions. We refer the reader to \cite{Adam, Liv} for standard background on knot theory. 

\begin{figure}[H]
\centering
\includegraphics[width=.51\linewidth]{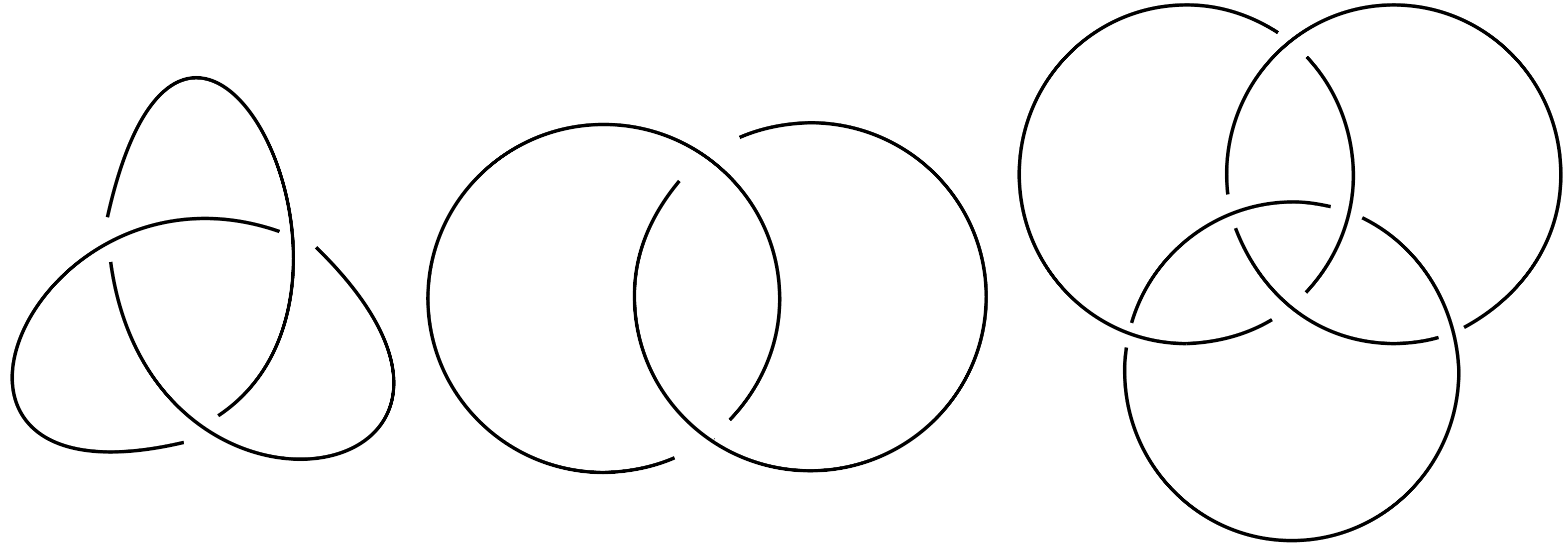}
\caption{(From left to right) Trefoil, Hopf link, denoted by $2_1$, and Borromean rings (3 components).}
\label{fig1}
\end{figure}

A {\em link} with $k$ components consists of $k$ disjoint simple closed curves ($\mathbb{S}^1$) in $\mathbb{R}^3$. A {\em knot} $K$ is a link with one component.  A {\em link diagram} $D(L)$ of a link $L$ is a regular projection of $L$ into $\mathbb{R}^2$ in such a way that the projection of each component is smooth and at most two curves intersect at any point. At each crossing point of the link diagram the curve which goes over the other is specified, see Figure \ref{fig1}.

\begin{figure}[H]
\centering
\includegraphics[width=.7\linewidth]{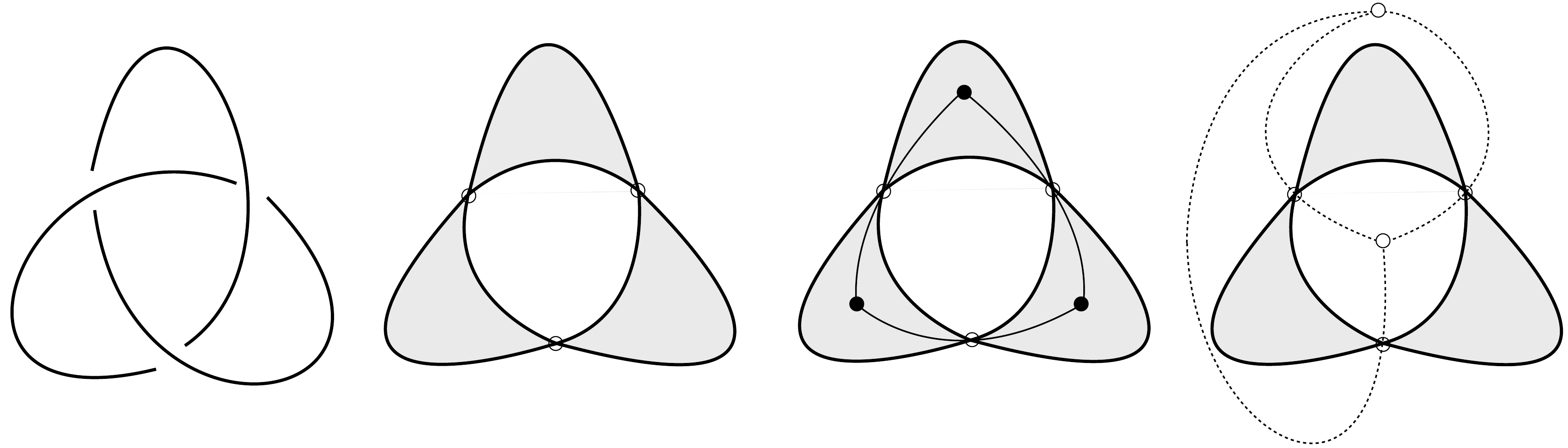}
\caption{(From left to right) A diagram of the Trefoil, its shadow with a 2-colored faces (vertices on white crossed circles), corresponding Black graph (bold edges and black circles) and White graph (dotted edges and white circles).}
\label{fig2}
\end{figure}

A {\em shadow} of a link diagram $D$ is a 4-regular graph if the over/under passes of $D$ are ignored. Since the shadow is Eulerian (4-regular) then its faces can be 2-colored, say with colors black and white. We thus have that each vertex is incident to 4 faces alternatively colored around the vertex, see Figure \ref{fig2}.

Given such a coloring, we can define two graphs, one on the faces of each color. Let $B_D$ denote the graph with black faces as its vertices and two vertices are joined if the corresponding faces share a vertex ($B_D$ is called the {\em checkerboard graph} of $D$). We define the graph $W_D$ on the white faces of the shadow analogously. The two faces are called {\em faces graphs} of the shadow. Since the shadow of a knot is connected then $B_D$ and $W_D$ are also connected, and it is not hard to see that $W_D=B_D^*$ and $B_D=W_D^*$, that is, the two faces graphs are duals of each other. 
\smallskip

It is not hard to convince oneself that a connected 4-regular planar graph $G$ is determined by either of its face graphs. For, recall that the {\em medial graph} of $H$, denoted by $med(H)$ is the graph obtained by placing one vertex on each edge of $H$ and joining two vertices if the corresponding edges are consecutive on a face of $H$. We notice that $med(H)$ is 4-regular since each edge is shared by exactly two faces. It is thus suffices to notice that $med(B_D)$ and $med(W_D)$ are the same (since they $B_D$ and $W_D$ are duals) and that $med(W_D)$ is exactly the shadow of $D$.
\smallskip

If $D$ is a link diagram with more than one component, then the shadow of $D$ may not be connected. In this case, only one of the corresponding faces graphs would be connected and $B_D^*\neq W_D$. We thus would not be able to determine uniquely the shadow of $L$ from the faces graphs. The latter can be overcome by considering each component of the shadow separately, and define the corresponding black face graph to be the union of the black face graphs of its components, and analogously for the white face graph.
\smallskip

%A \textit{signature} of a planar graph $G=(V,E)$ is a mapping $S:E\rightarrow \{+,-\}$. 
An \textit{edge-signed} planar graph, denoted by $(G,S_E)$, is a planar graph $G$ equipped with a signature on its edges $S_E:E\rightarrow \{+,-\}$. We will denote by $-S_E$ the signature of $G$ satisfying $-S_E(e)=-(S_E(e))$ for every $e\in E$. We write $S_E^+$ (resp. $S_E^-$) when all the signs of $S_E$ are $+$ (resp. $-$). 
\smallskip

\begin{figure}[H]
\centering
\includegraphics[width=.6\linewidth]{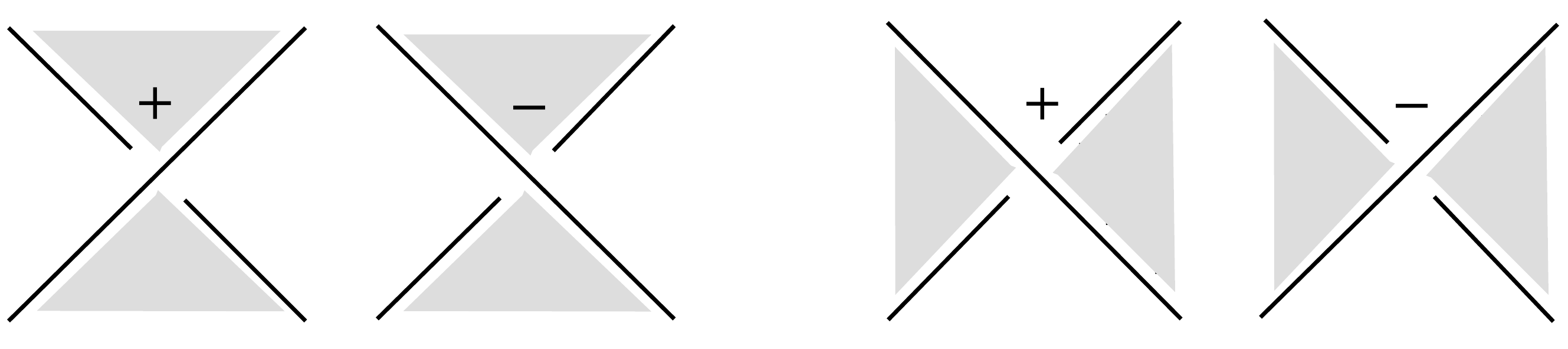}
\caption{(Left) Left-over-right rule from black point of view. (Right) Right-over-left rule from white point of view.}
\label{fig3}
\end{figure}

Given a crossing of the link diagram we sign {\em positive} or {\em negative} according to the {\em left-over-right} and {\em right-over-left} rules from point of view of black around the crossing, see Figure \ref{fig3} (a). The latter induce an opposite signing on each crossing by the same rules but now from point of view of white around the crossing, see Figure \ref{fig3} (b).

%Let $G_D$ be the graph where the set of vertices correspond to black faces and two vertices are joined if the corresponding faces share a vertex. We sign each edge of $G_D$ with the same sign of the corresponding shared vertex. We thus have that $G_D$ is a planar graph with signs on its edges. 

If the crossing is positive, relative to the black faces, then the corresponding edge is declared to be positive in $B_D$ and negative in $W_D$. Therefore, in this fashion, a link diagram $D$ determine a dual pair of signed planar graphs $(B_D,S_E)$ and $(W_D,-S_E)$ where the signs on edges are swapped on moving to the dual.  

\begin{remark}\label{ram;blackwhite} A link diagram can be uniquely recovered from either $(B_D,S_E)$ or $(W_D,-S_E)$.
\end{remark}

\begin{figure}[H]
\centering
\includegraphics[width=.7\linewidth]{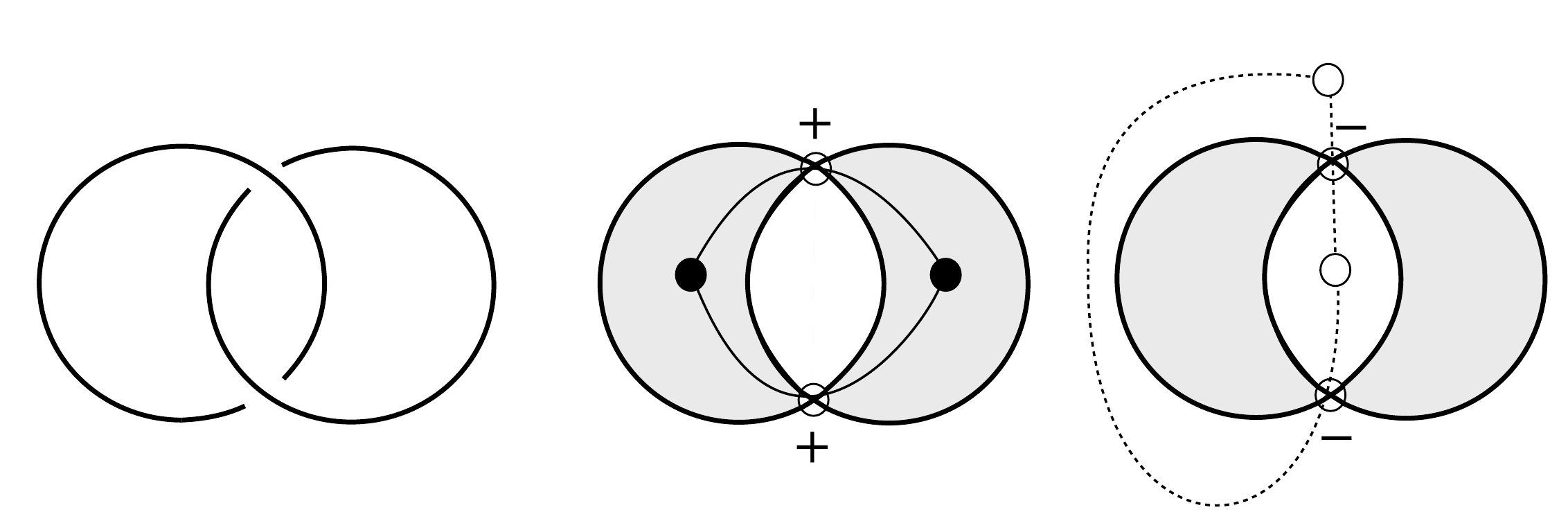}
\caption{(From left to right) diagram $D$ of link $2_1$, signed Black graph $(B_D,S_E)$ and signed White graph $(W_D,-S_E)$.}
\label{fig4}
\end{figure}

We thus have that given an edge-signed planar graph $(G,S_E)$, we can associate to it (in a canonical way) a link diagram $D(L)$ such that $(B_D,S_E)$ (and $(W_D,-S_E)$) gives $(G,S_E)$. The unsigned graph $G$ is called the {\em Tait graph} of the link $L$ with diagram $D$. The construction is easy, we just consider the $med(H)$ with signatures on the vertices (induced by the edge-signature $S_E$ of $G$). The desired diagram, denoted by $D(G,S_E)$, is obtained by determining the under/over pass at each crossing according to Left-over-right (or Right-over-left) rule associated to the sign of the corresponding edge of $G$, see Figure \ref{fig4}.

\begin{figure}[H]
\centering
\includegraphics[width=.5\linewidth]{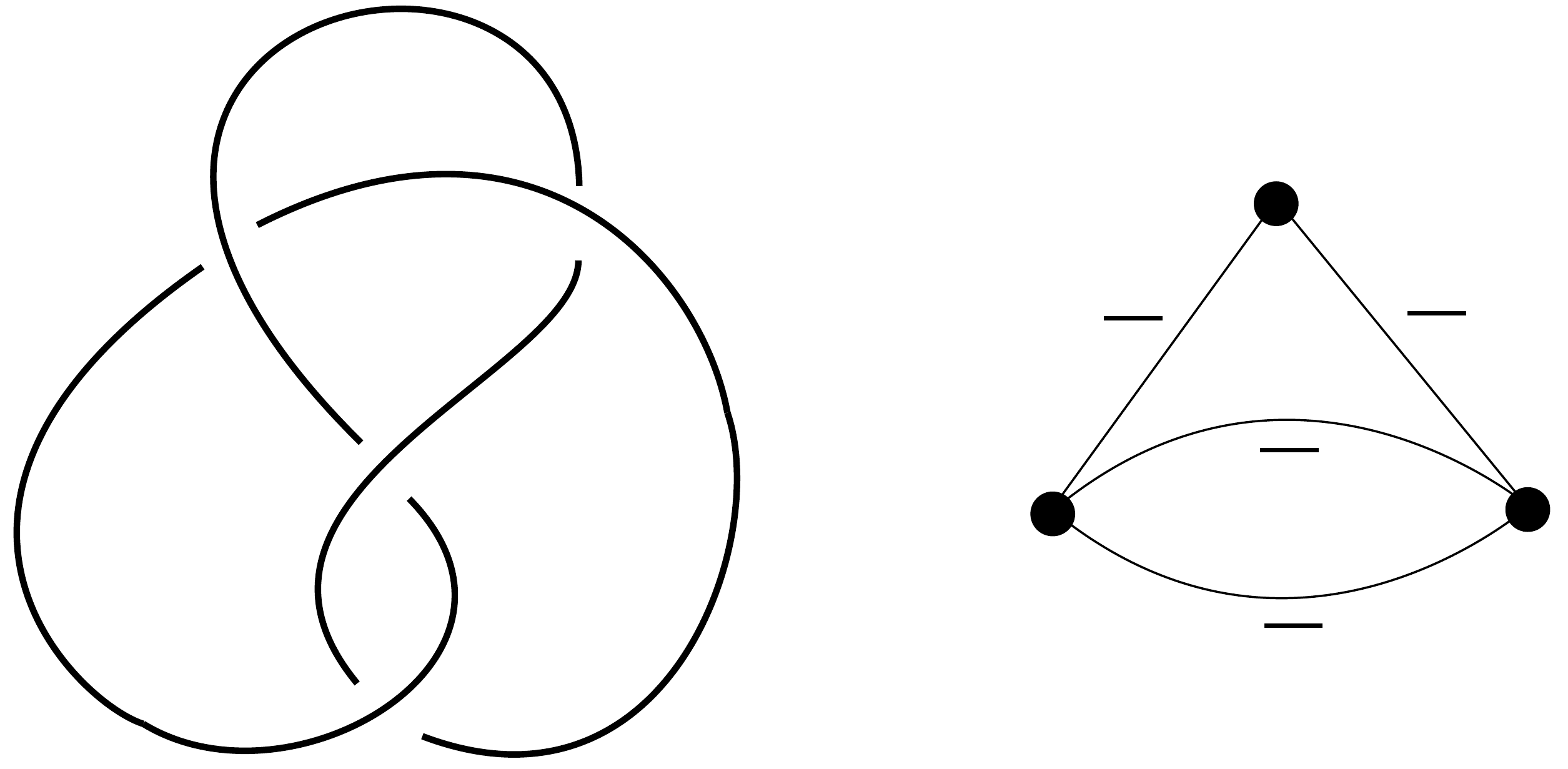}
\caption{A diagram $D$ of the {\em Figure-eight} knot (alternating knot, denoted by $4_1$) and its signed Black graph $(B_D,S_E^-)$.}
\label{fig5}
\end{figure}

A link $L$ is {\em alternating} if it admits a diagram $D$ such that the crossings alternate under/over while we go through the link. We notice that a link $L$ is alternating if and only if $L=D(G,S_E^+)=D(G,S_E^-)$ where $G$ is the Tait graph of link $L$, see Figure \ref{fig5}.

\subsection{Amphichirality}\label{sec:amphi}

Two links $L_1$ and $L_2$ are {\em equivalent} if there is an orientation-preserving homeomorphism 
$\varphi :\mathbb{R}^3\rightarrow  \mathbb{R}^3$ with $\varphi(L_1)=L_2$. 
%If the map $\phi$ is orientation-preserving then it is necessarily an {\em isotopy} this means that there is a family $\phi_t$ of piecewise linear homeomorphisms with $t\in [0,1]$ such that $\phi_0$ is the identity, $\phi_1=\phi$ and  the map $$(x,t)\mapsto (\phi_t(x),t)$$ is a piecewise-linear homeomorphism from $\mathbb{S}^3\times [1,0]$ to itself.   
\smallskip

We notice that, in the same way as infinitely many signed plane graphs represent equivalent links  for any link $L$,  there are infinitely many link diagrams that can represent equivalent links.
We shall denote by $[L]$ (called {\em link-type}) the class of links equivalent to $L$.
\smallskip

The mirror-image $L^*$ of a link $L$ is obtained by reflecting it in a plane through the origin in $\mathbb{R}^3$. It may also be defined given a diagram $D$ of $L$ by simply exchanging all the crossing of $D$ (this is clear if one consider reflecting the knot in the plane where the diagram is drawn). Although such a reflection is a bijective map, it is not orientation preserving (its determinant is negative). Therefore, it may happen that $L$ and $L^*$ are not equivalent (that is, $L\not\in [L^*]$ and $L^*\not\in [L]$).  
\smallskip

Link $L$ is called {\em amphichiral} (also known as {\em achiral}) if there is $L\in [L^*]$. Equivalently, $L$ is amphichiral if there is an automorphism of $\rth$ (or $\sth$) preserving $L$ and reversing the orientation. To see this equivalence, we only need to compose the preserving automorphism with a reflection in a plane to get $L^*$ (such composition reverses twice the orientation).    

%%%%%%%%%%%%%%%%%%%%%%%%%%%%%%%%%%%%%%%
\section{Maps: background}\label{sec;maps}

We review some notions and properties on maps needed for the rest of the paper.
\smallskip

A {\em map} of $G=(V,E,F)$ is the image of an embedding of $G$ into $\mathbb{S}^2$ where the set of vertices are a collection of distinct points in $\mathbb{S}^2$ and the set of edges are a collection of Jordan curves joining two points in $V$ satisfying that $\alpha\cap\alpha'$ is either empty or a point in the endpoints for any pair of Jordan curves $\alpha$ and $\alpha'$.  Any embedding of the topological realization of $G$ into $\mathbb{S}^2$ partitions the 2-sphere into simply connected regions of $\mathbb{S}^2\setminus G$ called the {\em faces} $F$ of the embedding.
\smallskip

Let us define the {\em antipodal} function as

$$\begin{array}{llll}
\alpha_n: & \mathbb{S}^n& \rightarrow & \mathbb{S}^n\\ 
& x & \mapsto & -x
\end{array}$$

We say that $Y\subseteq\mathbb{S}^n$ is $n$-{\em antipodally symmetric} if $\alpha_n(Y)=Y$.
\smallskip

A self-dual map $G$ is called {\em antipodally self-dual} if the dual map $G^*$ is antipodally embedded in $\stw$ with respect to $G$, that is, $\alpha_2(G)=G^*$. We say that $G$ is 2-{\em antipodally symmetric} map if it admits an embedding in $\stw$ such that $\alpha_2(G)=G$.
%for any $x$ in $G$ if $x\in V(G)$ then $-x\in V(G)$ and if $x\in E(G)$ then $-x\in E(G)$.
\smallskip
  
\begin{lemma} \label{lem:key} \cite[Lemma 1]{MRAR1} If G is an antipodally self-dual map then $med(G)$ is 2-antipodally symmetric.
\end{lemma}

Let $G$ be a 2-antipodally symmetric map. If $v\in V(G)$ then its {\em antipodal} vertex is given by $\alpha_2(v)=-v$. We call them {\em antipodal pair} of vertices. 

\begin{remark}\label{rem:2-aut} If $G$ is 2-antipodally symmetric map then its number of faces must be even. Moreover, the function $\alpha_2$ naturally matches the pairs of {\em antipodal faces}, say $f$ and $\alpha_2(f)$ (we may refer $\alpha_2(f)$ as the {\em $f$-antipodal} face of $f$). The latter naturally induces a permutation of the faces that turns out to be an automorphism of $G^*$ (that is, $\alpha_2\in Aut(G)$, and thus $G^*$ is also 2-antipodally symmetric). 
\end{remark}

\begin{proposition} Let $G$ be 2-antipodally symmetric map where its faces are 2-colored properly (that is, two faces sharing an edges have different colors). Then, if one pair of antipodal faces have the same (resp. different color) then all pairs of antipodal faces have the same (resp. different color). 
\end{proposition}

\begin{proof} Suppose that $f$ and $f$-antipodal are both colored with the same color, say black, then any face $f_1$ adjacent to $f$ is colored white and thus $f_1$-antipodal must also be colored white since it is adjacent to $f$-antipodal, by carry on this argument, we end with $f$ and $f$-antipodal having the same color for all faces $f$. Similar argument can be applied when $f$ and $f$-antipodal have different colors.
\end{proof}

%We say that a 2-coloring of faces of a 2-antipodally symmetric map $G$ is {\em preserving} (resp.{\em reversing}) if $f$ and $f$-antipodal have the same color (resp. different colors) for each pair of antipodal faces, see Figure \ref{fig14}.

A {\em bicolored map} is a map  $G=(V,E,F)$ together with a coloring $C_X:X\rightarrow \{black, white\}$ where $X$ is either $V(G), E(G)$ or $F(G)$. A {\em signed map} is map $G=(V,E,F)$ together with a signature $S_Y:Y\rightarrow \{+,-\}$ where $Y$ is either $V(G), E(G)$ or $F(G)$. 
\smallskip

Throughout the paper, we will consider bicolored signed maps $(G,C_X,S_Y)$, that is, maps together with both a {\em vertex-, edge-} or {\em face-coloring} $C_X$ and a {\em vertex-, edge-} or {\em face-signature} $S_Y$.
\smallskip

Let $(G,C_F,S_V)$ be a colored-face vertex-signed map. We say that an automorphism $\sigma(G)\in Aut(G)$ is \textit{color-preserving} (resp. \textit{color-reversing}) if each pair of faces $f$ and $\sigma(f)$ have the same (resp. different) color. Similarly, $\sigma$ is said to be \textit{sign-preserving} (resp. \textit{sign-reversing}) if each pair of vertices $v$ and $\sigma(v)$ have the same (resp. different) sign.

\begin{figure}[H]
\centering
\includegraphics[width=.7\linewidth]{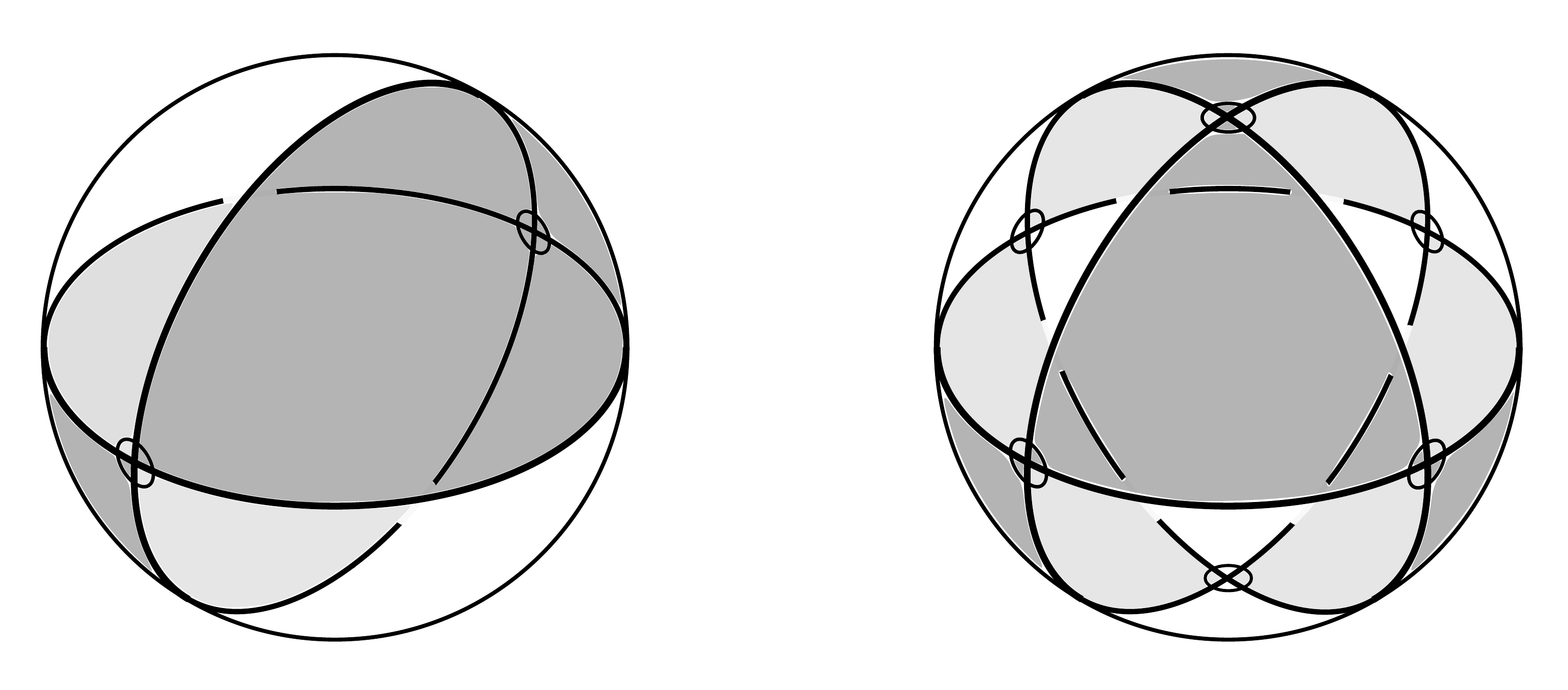}
\caption{(Left) A $2$-antipodally symmetric map where the antipodal mapping is color-preserving. (Right) A $2$-antipodally symmetric map where the antipodal mapping is color-reversing.}\label{fig14}
\end{figure}

Notice that in the case when $(G,C_F,S_V)$ is a 2-antipodally symmetric map, the automorphism $\alpha_2$ can be either color-preserving or color-reversing, see Figure \ref{fig14}.

%Let $(G,C_F,S_V)$ be a 2-antipodally symmetric vertex-signed map.  We say that the antipodal mapping $\alpha_2\in Aut(M)$ is \textit{sign-preserving} (resp. \textit{sign-reversing}) if the pair of vertices $v$ and $\alpha_2(v)$ have same (resp. different) signs. Similarly, we say that $\alpha_2$ is \textit{color-preserving} (resp. \textit{color-reversing}) if for every face $f$ the pair $f$ and $\alpha_2(f)$ have the same (resp. different) color.

%Let $G$ be a connected $2$-antipodally symmetric map where its faces are $2$-colored. 
%Then, the faces of $M$ can be $2$-colored such that every two faces sharing en edge have different colors. 

%This coloring of the faces is usually called the \textit{checkerboard coloring}. We say that the antipodal mapping $\alpha_2\in Aut(M)$ is \textit{color-preserving} (resp. \textit{color-reversing}) if for every face $f$ the pair $f$ and $\alpha_2(f)$ have same (resp. different) colors.

\begin{remark}\label{rem;color} If $G$ is antipodally self-dual then $med(G)$ admits a 2-antipodally symmetric embedding $\alpha_2$ in $\stw$ with $\alpha_2$ coloring-reversing. Indeed, since $G$ is antipodally self-dual then
antipodal faces of $med(G)$ correspond to a pair of antipodal vertices, say $v\in V(G)$ and $\alpha_2(v)\in V(G^*)$. Therefore, if we color all faces corresponding to vertices in $G$ (resp. in $G^*$) in black (resp. in white) we obtain that  $\alpha_2$ is coloring-reversing.
%In fact, in \cite[Theorem 1]{MRAR1} was proved that if $G$ is an antipodally self-dual map then its {\em vertex-face incident} graph $I(G)$ always admits at least one {\em symmetric cycle} of length $2n$ with $n\ge1$ odd. The latter implies that $med(G)$ admits a reversing 2-coloring.
\end{remark}

%We notice that if $M=med(G)$ then its dual is always $2$-colorable. 
%with say black the vertices of $G$ and with white the vertices of $G^*$.

%
%\begin{figure}[H]
%\centering
%\includegraphics[width=.5\linewidth]{2-colorings}
%\caption{(a) The antipodally symmetric map $med(K_4)$. The antipodal faces $(a,A), (b, B), (c, C)$ and $(d,D)$ with $A$ the exterior face with a reversing 2-coloring. (b) The antipodally symmetric  graph $med(K_2^2)$ where $K_2^2$ consists of 2 vertices and 2 parallel edges. The antipodal faces $(a,A), (b,B)$ with $A$ the exterior face with  a preserving 2-coloring.}
%\label{fig14}
%\end{figure}

%\subsection{Antipodally symmetric links in $\rth$}

\subsection{Incidence graph} 
\begin{figure}[H]
    \centering
    \includegraphics[width=0.4\textwidth]{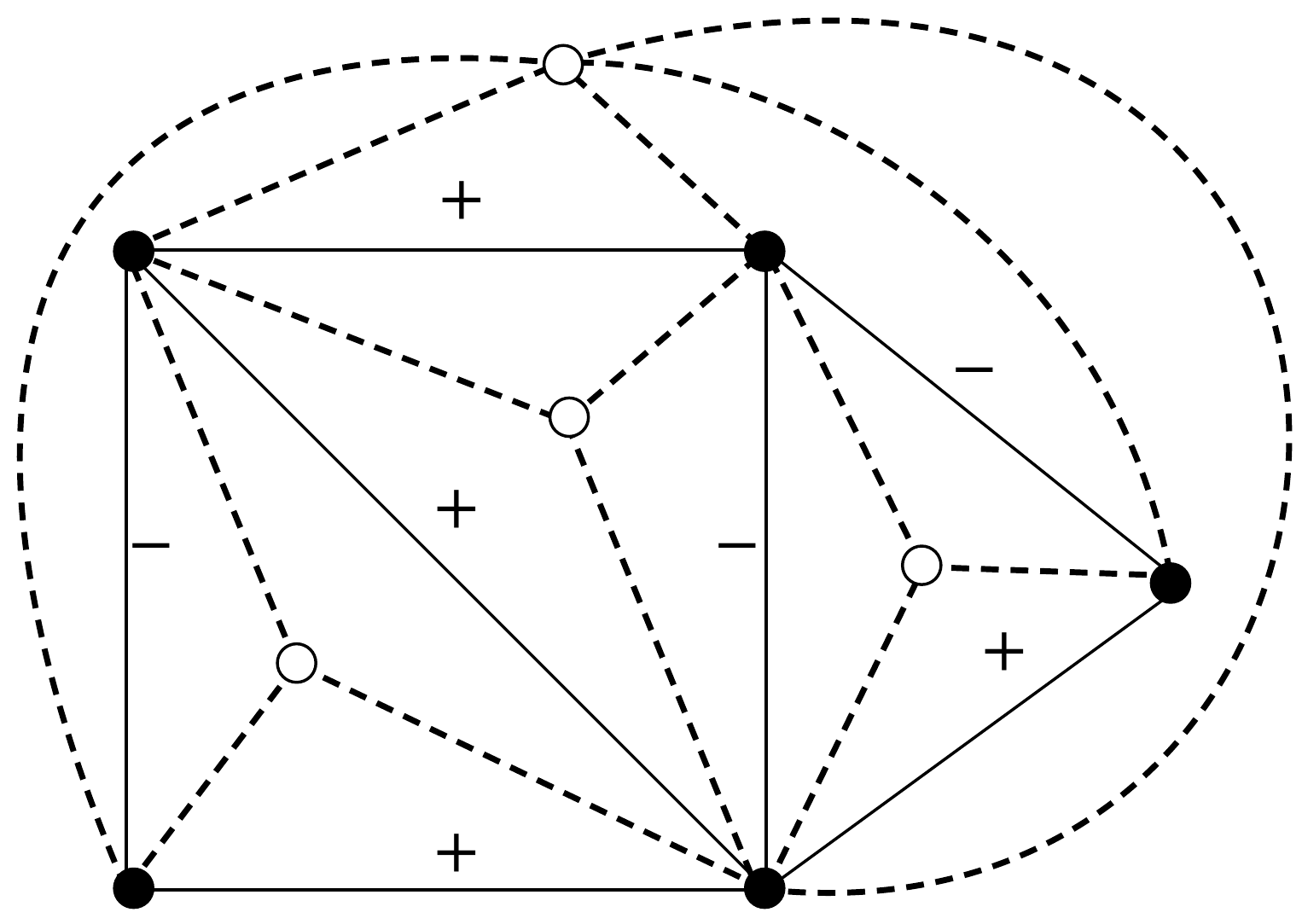}
  \caption{Edge-signed graph $(G,S_E)$ (light edges) and its face-signed, vertex-colored incidence graph $(I(G),S_F)$ (dashed edges).} \label{fig255}
\end{figure}

Let $G$ be a plane graph and $G^*$ its geometric dual. We recall that  the {\em (vertex-face) incidence} graph $I(G)=(V^{I},E^{I})$ 
%is the graph obtained from a planar graph $G$ having as set of vertices $V^{I}=V(G)\cup V(G^*)$ and edges all incidence diagonals of faces in $G^\square$. 
 has as vertices $V^{I}=V(G)\cup V(G^*)$ and two vertices $v\in V(G)$ and $v^*\in V(G^*)$ are adjacent if $v$ is a vertex of the face corresponding to $v^*$,  see Figure \ref{fig255}
 %The latter naturally induce a 2-coloring on the vertices of $I(G)$. 

%Also  can be constructed through the White and Black graphs by placing one white (resp. black) vertex in each white (resp. black) face and joining a white vertex with a black vertex if and only if the corresponding faces share an edge, see Figure \ref{fig25}.

\begin{remark}\label{rem;I(G)} Let $G$ be a plane graph.

(a) $I(G)=med(G)^*$.

(b) $I(G)$ is bipartite, the color of a vertex $x$ is black (resp. white) if  $x\in V(G)$ (resp. $x\in V(G^*)$.

(c) Each face of $I(G)$ is of length four (each face contains exactly one edge of $G$ as a diagonal). 
\end{remark}

%Let $G$ be the Tait signed graph of a link $L$. As noticed in Remark \ref{rem;I(G)}, each face of $I(G)$ contains a vertex of $med(G)$ corresponding to a crossing of the associated diagram $D(L)$ of $L$. The latter induce two different signatures to the faces of $I(G)$ (one opposite to the other) according to the edge signature of either $(B_D,S)$ or $(W_D,-S)$. 

%Let $(G,S_E)$ be an edge-signed map. We denote by $(I(G), S_F)$ the graph $I(G)$ equipped with the face-signature $S_F$ induced by $S_E$.
%\smallskip

\begin{figure}[H]
    \centering
    \includegraphics[width=0.7\textwidth]{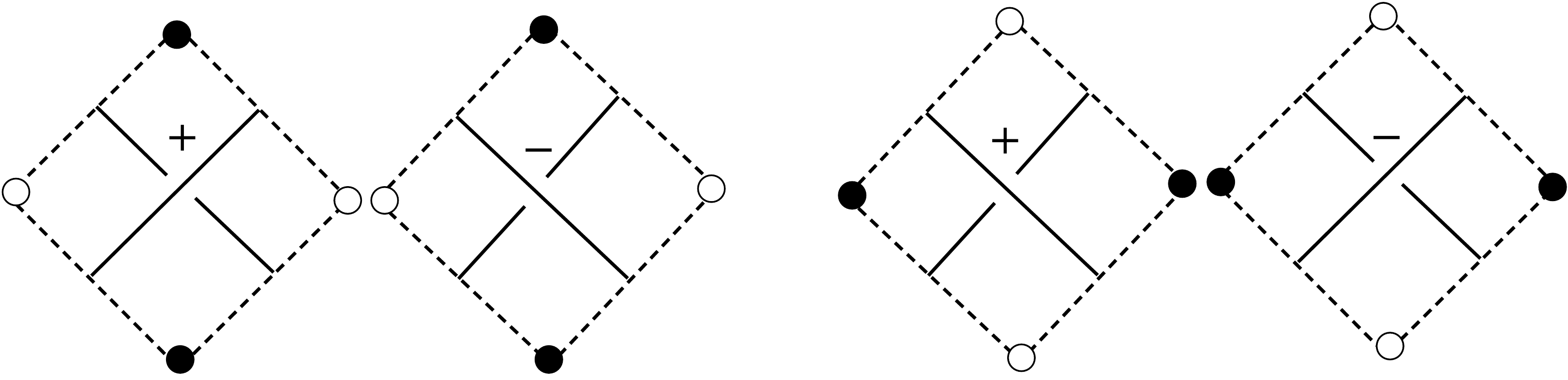}
    \caption{(From left to right) Left-over-right rule from black point of view and  Right-over-left rule from white point of view.}
    \label{fig27}
\end{figure}

We have that $(I(G), C_V, S_F)$ also determines (in a canonical way) the link diagram $D(G,S_E)$. The construction is easy, we just determine the under/over pass at each each crossing (in each face) according to either the Left-over right rule or Right-over left rule, see Figure \ref{fig27}. 

We denote by $C_V^o$ the {\em opposite} vertex-coloring of $C_V$, that is, white vertices become black and black ones become white.

\begin{remark}\label{rem:I(G)} Let $(G,S_E)$ be an edge-signed map and let $(I(G), C_V,S_F)$ be its vertex-face incidence graph equipped with a face-signature $S_F$ (arising from $S_E$) and vertex-coloring $C_V$. We have that if $(I(G), C_V, S_F)$ determines link $D(G,S_E)$  then both $(I(G), C_V^o,S_F)$ and $(I(G), C_V,-S_F)$ determine $D(G,S_E)^*$.
\end{remark}  

\subsection{Special embedding of a link}
\begin{figure}[H]
    \centering
    \includegraphics[width=1\textwidth]{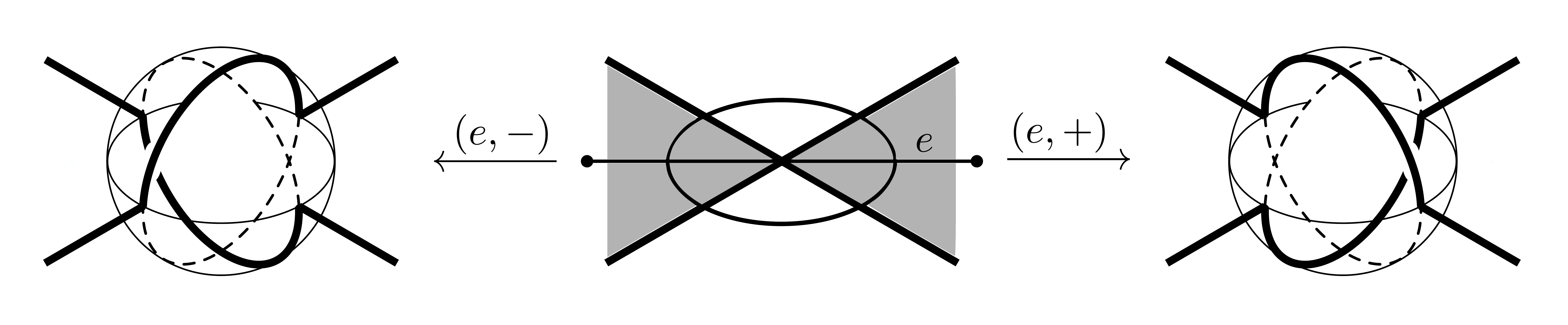}
    \caption{The crossing spheres rules.}
    \label{crossing-spheres}
\end{figure}

Let $D(G,S_E)$ be the link diagram in $\stw$ obtained from a map $(med(G),C_F,S_V)$.
We shall construct a specific embedding of $D(G,S_E)$ in $\rth$
by modifying (locally) the diagram around each crossing. We proceed as follows. Take a small sphere $\mathbb{S}_1$ around each crossing (say, with center the crossing itself). and move (locally) the piece of arc of the diagram passing over (resp. passing under) around $\mathbb{S}_1$ outside (resp. inside) of $\mathbb{S}^2$ according with the {\em crossing sphere rules}, see Figure \ref{crossing-spheres}. 

\begin{figure}[H]
\centering
\includegraphics[width=.7\linewidth]{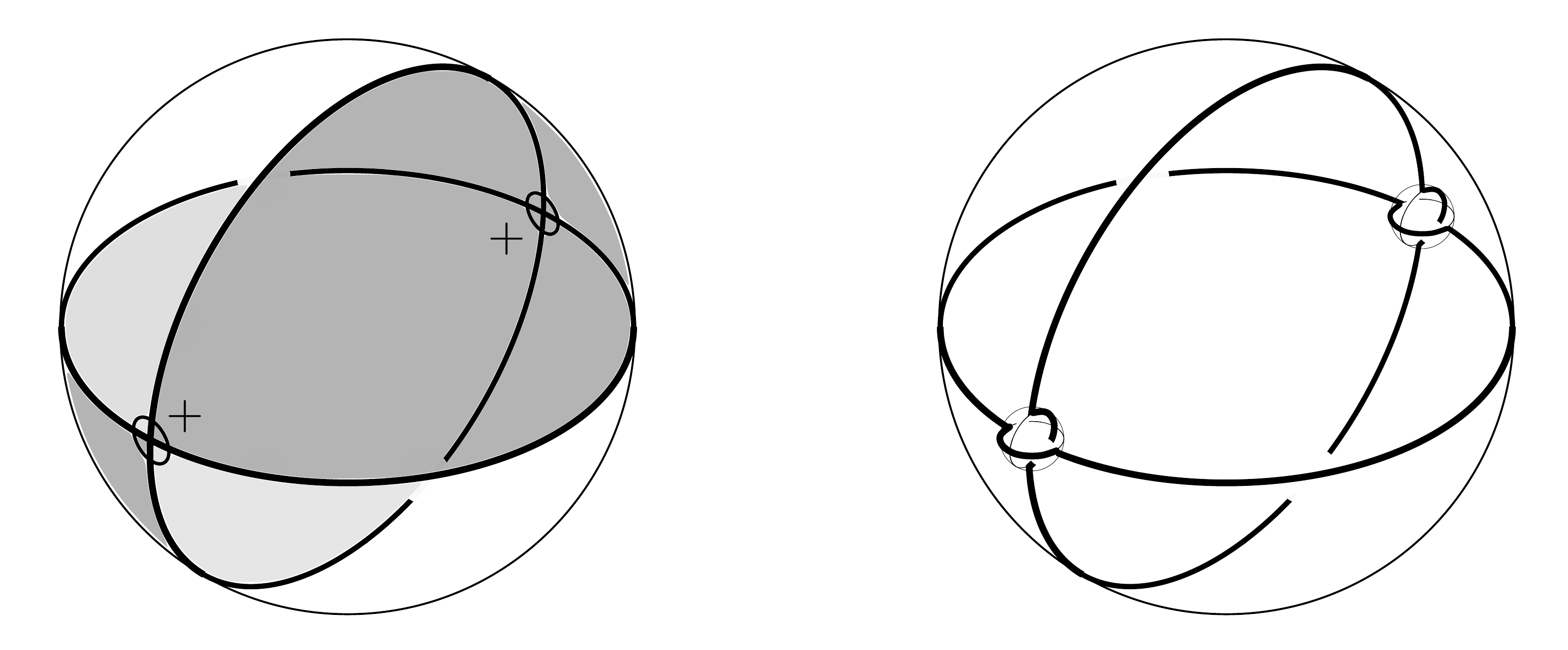}
\caption{(Left) A 2-antipodally symmetric map of the graph on 2 vertices and four parallel edges. This graph corresponds to a shadow of the Hopf link, see Figure \ref{fig4} (second row). In this case, $\alpha_2$ is color-preserving and sign-preserving. (Right) $\embed(2_1)$.}
\label{fig16}
\end{figure}

The rest of the diagram $D(G,S_E)$ remains the same in $\mathbb{S}^2$. We denote by $\embed(G,S_E)$ such embedding,
see Figure \ref{fig16}.

From the above discussion we have the following diagram.

$$\begin{array}{rcl}
(med(G), C_F, S_V) & \longleftrightarrow & (I(G),C_V,S_F)\\
\searrow \hspace{-.41cm}\nwarrow& &\swarrow \hspace{-.41cm}\nearrow \\
& (G,S_E) &\\
 &\downarrow &\\
 & D(G,S_E)&\\
  &\updownarrow &\\
 & L(G,S_E)&\\
\end{array}$$

where $A\rightarrow B$ means {\em $B$ can be constructed from $A$}.

\section{Central symmetry in $\mathbb{R}^3$}\label{sec;R3sym}

%From now on, the link $D(G,S_E)$ will be understood as the link corresponding to the diagram in $\stw$ arising from a edge-signed map $(G,S_E)$. More precisely, $D(G,S_E)$ is the diagram in $\mathbb{S}^2$ induced by $(med(G),S_V)$ ($S_V$ correspond to the edge-signature $S_E$ of $G$). 

Let us define the {\em centrally symmetric} function as

$$\begin{array}{llll}
c_n: & \mathbb{R}^n& \rightarrow & \mathbb{R}^n\\ 
& x & \mapsto & -x
\end{array}$$

We say that $Y\subseteq\mathbb{R}^n$ is $n$-{\em centrally symmetric} if $c_n(Y)=Y$. We say that $L$ is 3-centrally symmetric if there is $\hat L\in [L], \hat L\in \mathbb{R}^3$ such that $c_3(\hat L)=\hat L$.

\begin{theorem} \label{thm:antipodalR3} Let $(G,S_E)$ be an face-colored edge-signed self-dual map and suppose that $med(G)$ is a $2$-antipodally symmetric map (realized by $\alpha_2$). If either
%Then there is $M$ isomorphic to $med(G)$ such that the antipodal mapping of $\stw$ $a_{\stw}:x\mapsto-x$ is in $Aut(M)$. Then we have
\smallskip

a) $\alpha_2$ is color-preserving and sign-reversing  or

b) $\alpha_2$ is color-reversing and sign-preserving 
\smallskip

\noindent then $L(G,S_E)$ is $3$-centrally symmetric. 
\end{theorem}

\begin{proof} Let $D(G,S_E)$ be the link diagram induced by $(med(G),C_F,S_V)$.  Since $med(G)$ is 2-antipodally symmetric then it is also 3-centrally symmetric. Consider the embedding  $\embed(G,S_E)$. It can be checked that if $\alpha_2$ is either (a) color-preserving and sign-reversing or (b) color-reversing and sign-preserving then the piece of arc of the diagram passing over (resp. passing under) in crossing $v$ correspond to the piece of arc of the diagram passing over (resp. passing under) in crossing $\alpha_2(v)$ in $\embed(G,S_E)$, see Figure \ref{fig:3centsym}.  

\begin{figure}[H]
\centering
\includegraphics[width=\linewidth]{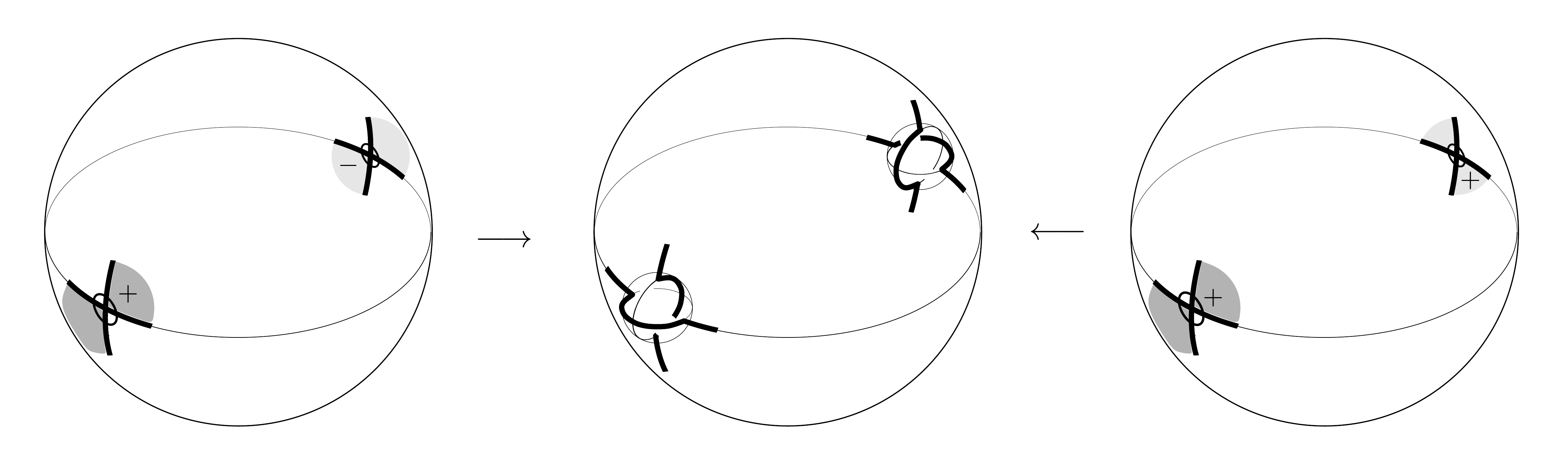}
\caption{(Left) Antipodal pair of vertices of $med(G)$ in case (a)  (Right)  Antipodal pair of vertices of $med(G)$ in case (b) (Center) The local modifications around an antipodal pair of vertices following the crossing sphere rules.}
\label{fig:3centsym}
\end{figure}

Therefore, the modifications around $v$ are centrally symmetric with respect to those done around $\alpha_2(v)$. We thus obtain that $\embed(G,S_E)$ is a $3$-centrally symmetric.
\end{proof}

%\begin{figure}[H]
%\centering
%\includegraphics[width=.38\linewidth]{obspoint}
%\caption{Viewing crossing from above.}
%\label{fig15}
%\end{figure}

%\begin{figure}[H]
%\centering
%\includegraphics[width=.42\linewidth]{localperturbation}
%\caption{Local perturbation around a crossing.}
%\label{fig15a}
%\end{figure}

%\begin{figure}[H]
%\centering
%\includegraphics[width=.8\linewidth]{mouvementsym}
%\caption{Two possible perturbations around crossings $v$ and $-v$ (a) A preserving 2-coloring and antipodal pair of vertices with opposite signs. (b) A reversing 2-coloring and antipodal pair of vertices with the same sign.}\label{fig17}
%\end{figure}
 
Figure \ref{fig16} illustrates an example in which neither condition (a) nor (b) of Theorem \ref{thm:antipodalR3} are full filled. In this case, $\embed(2_1)$ does not give a 3-centrally symmetric embedding of $2_1$. In fact, it can be showed that the Hopf link does not admit a 3-centrally symmetric embedding.

\begin{corollary}\label{cor:neg-amphi} Let $(G,S_E)$ be an edge-signed map and suppose that $med(G)$ is a 2-antipodally symmetric map (realized by $\alpha_2$). If either
%Then there is $M$ isomorphic to $med(G)$ such that the antipodal mapping of $\stw$ $a_{\stw}:x\mapsto-x$ is in $Aut(M)$. Then we have
\smallskip

a) $\alpha_2$ is color-preserving and sign-reversing  or

b) $\alpha_2$ is color-reversing and sign-preserving 
\smallskip

\noindent then $L(G,S_E)$ is amphichiral.
\end{corollary}

\begin{proof} By Theorem \ref{thm:antipodalR3}, 
%$[D(G,S_E)]$ is 3-centrally symmetric and in particular 
$c_3(\embed(G,S_E))=\embed(G,S_E)$. We shall perform 3 reflections to this embedding with respect to the 3 orthogonal planes $x=0, y=0$ and $z=0$. Let $(x,y,z)\in L'$. We have that the first reflection maps $(x,y,z)$ into $(-x,y,z)$ the second maps $(-x,y,z)$ into $(-x,-y,z)$ and the third one maps $(-x,-y,z)$ into $(-x,-y,-z)$. Since $c_3(\embed(G,S_E))=\embed(G,S_E)$ then these three reflections give a nonpreserving orientation homeomorphism mapping $\embed(G,S_E)$ into itself, implying that $L(G,S_E)$ is amphichiral. 

%Finally, we notice that each reflexion change the orientation of knot $D(G,S)$ and thus arriving with reversed orientation after the three reflexions.
\end{proof}

\subsection{Families of amphichiral links}
\begin{corollary}\label{cor:antipodal} Let $(G,S_E^+)$ be an edge-signed map with $G$ antipodally self-dual. Then, $L(G,S_E^+)$ is amphichiral. 
\end{corollary}

\begin{proof} Since $G$ is antipodally self-dual then, by Lemma \ref{lem:key}, $med(G)$ is a 2-antipodally symmetric map realized by, say $\alpha_2$. Since all the signs in $S_E$ are the same then, in particular,  each pair of antipodal vertices of $med(G)$ have the same sign, that is, $\alpha_2$ is sign-preserving. Moreover, by  Remark \ref{rem;color}, $\alpha_2$ is color-reversing with respect to $med(G)$. The result follows by Corollary \ref{cor:neg-amphi} (b). 
\end{proof}

\begin{figure}[H]
\centering
\includegraphics[width=0.6\linewidth]{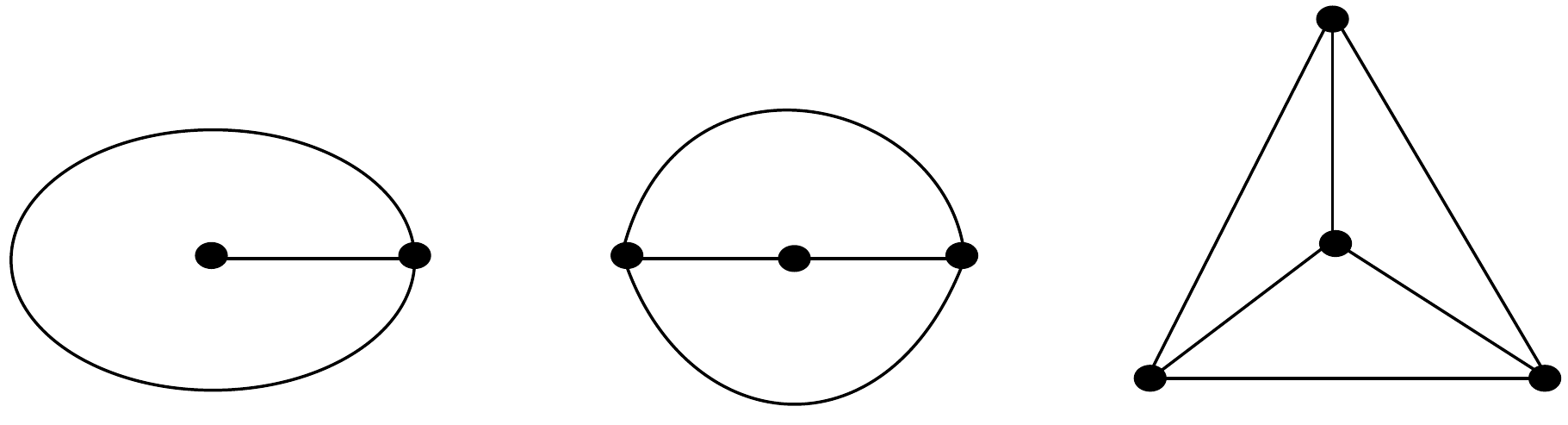}
\caption{$W_1$, $W_2$ and $W_3$.}
\label{fig201}
\end{figure}

Let $n\ge 1$ be an integer. An {\em $n$-wheel}, denoted by $W_n$, is the graph consisting of an $n$-cycle with a center joined to each vertex of the cycle, see Figure \ref{fig201}

%Let $n\ge 1$ be an odd integer. Let $L(W_n,S_E)$ be the family of link diagrams, called {\em $n$-wheel links}, having $W_n$ as Tait's graph and an arbitrary edge-signature $S_E$. 
 
\begin{corollary}\label{cor;w} Let $n\ge 1$ be an odd integer and let $(W_n,S_E)$ be a map such that $S_E(v)=S_E(-v)$ for each pair of antipodal vertices of $med({W_n})$. Then, $L(W_n,S_E)$ is amphichiral. 
\end{corollary}

\begin{proof} By \cite[Proposition 1]{MRAR1}, $W_n$ is an antipodally self-dual map and thus, by Lemma \ref{lem:key}, $med(W_n)$ is a 2-antipodally symmetric map realized by, say $\alpha_2$. Each pair of antipodal vertices of $med({W_n})$ have the same sign, that is, $\alpha_2$ is sign-preserving. Moreover, by  Remark \ref{rem;color}, $\alpha_2$ is color-reversing with respect to $med(W_n))$. The result follows by Corollary \ref{cor:neg-amphi} (b). 
\end{proof}

%By Remark \ref{rem;color}, $med({W_n})$ admits a reversing 2-coloring for any odd integer $n\ge 3$, see Figure \ref{fig19a} for the case $n=5$.

%\begin{figure}[H]
%\centering
%\includegraphics[width=0.32\linewidth]{2-Col-medW5}
%\caption{Reversing 2-colorings of $med({W_5})$ where antipodal faces are given by $i$ and $i^*$.}
%\label{fig19a}
%\end{figure}

%\jorge{- Maybe it would be helpful to  illustrate the corresponding embedding in the 2-sphere ??}
\smallskip

\begin{figure}[H]
\centering
\includegraphics[width=0.7\linewidth]{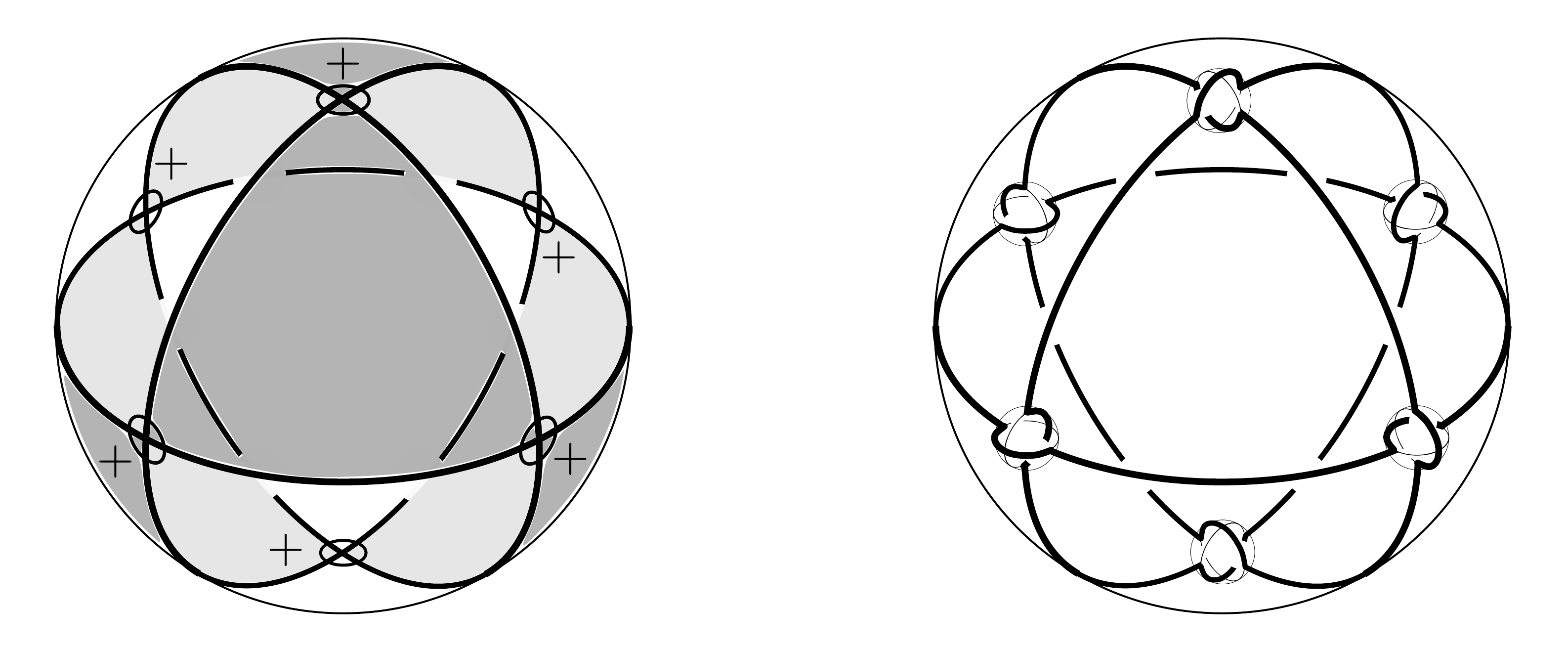}
\caption{A 2-colored map $(W_3, S_E^+)$ and  $\embed(W_3, S_E^+)$.}
\label{fig20a}
\end{figure}

In the case $n=3$, $L(W_3, S_E^+)$ is the well-known {\em Borromean rings}, a 3-centrally symmetric embedding is illustrated in Figure \ref{fig20a}.

\begin{figure}[H]
\centering
\includegraphics[width=0.3\linewidth]{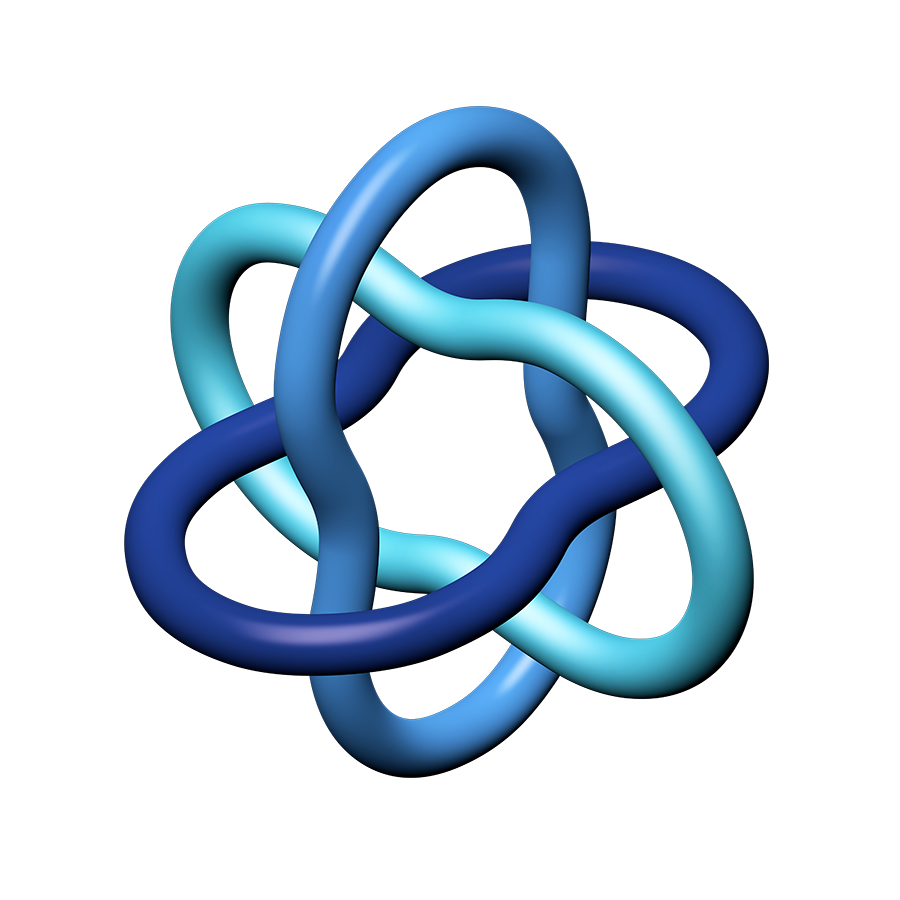}
\caption{International Mathematical Union logo.}
\label{fig20a1}
\end{figure}

An appropriate slightly straightening around the crossings arise the IMU's logo\footnote{See https://www.mathunion.org/outreach/imu-logo}, see Figure \ref{fig20a1} .

\begin{figure}[H]
\centering
\includegraphics[width=0.5\linewidth]{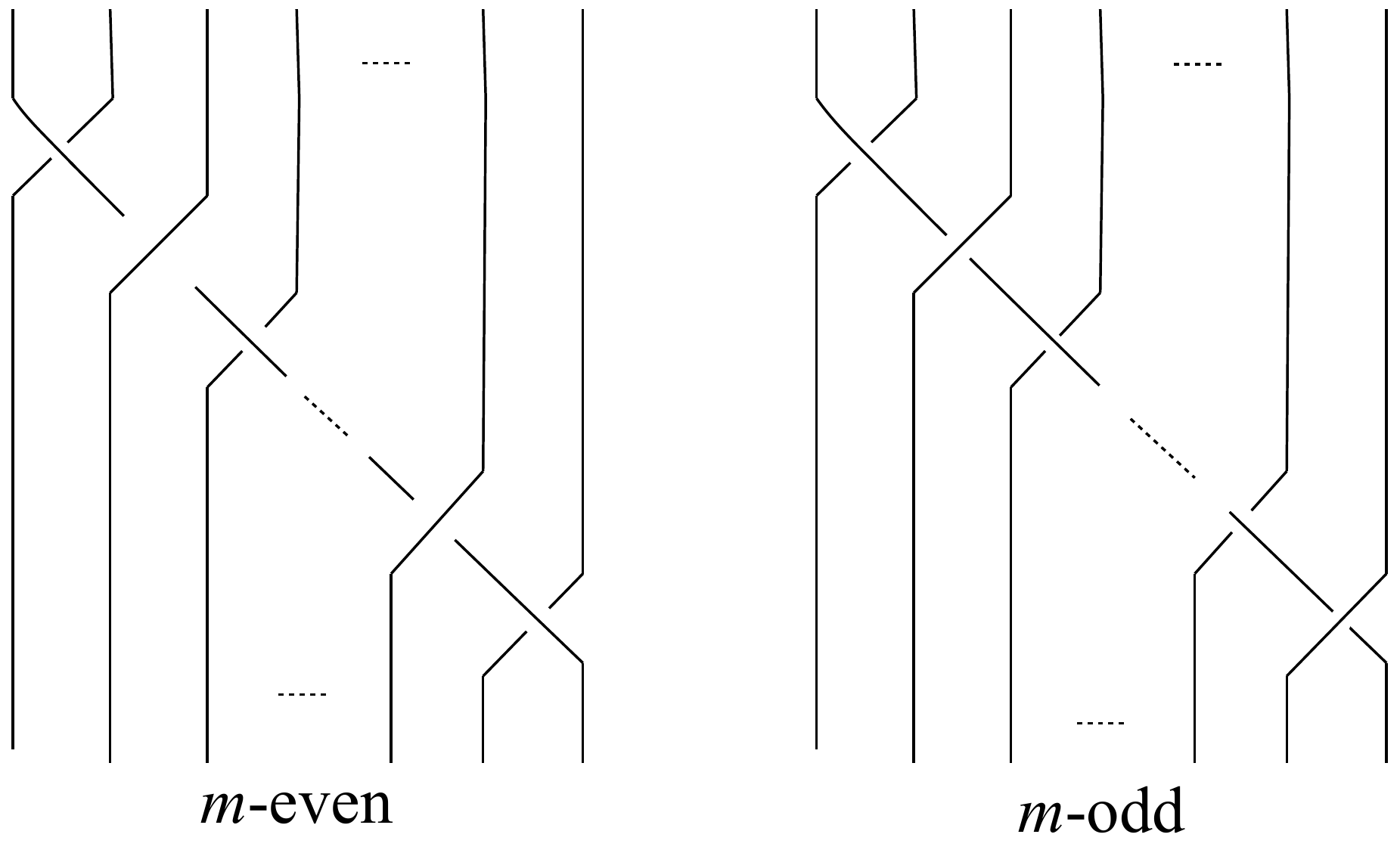}
\caption{Braid $A_m$ according with the parity of $m$.}
\label{fig19b}
\end{figure}

Let $A_m$ be the $m$-strand braid as shown in Figure \ref{fig19b}

\begin{figure}[H]
\centering
\includegraphics[width=0.6\linewidth]{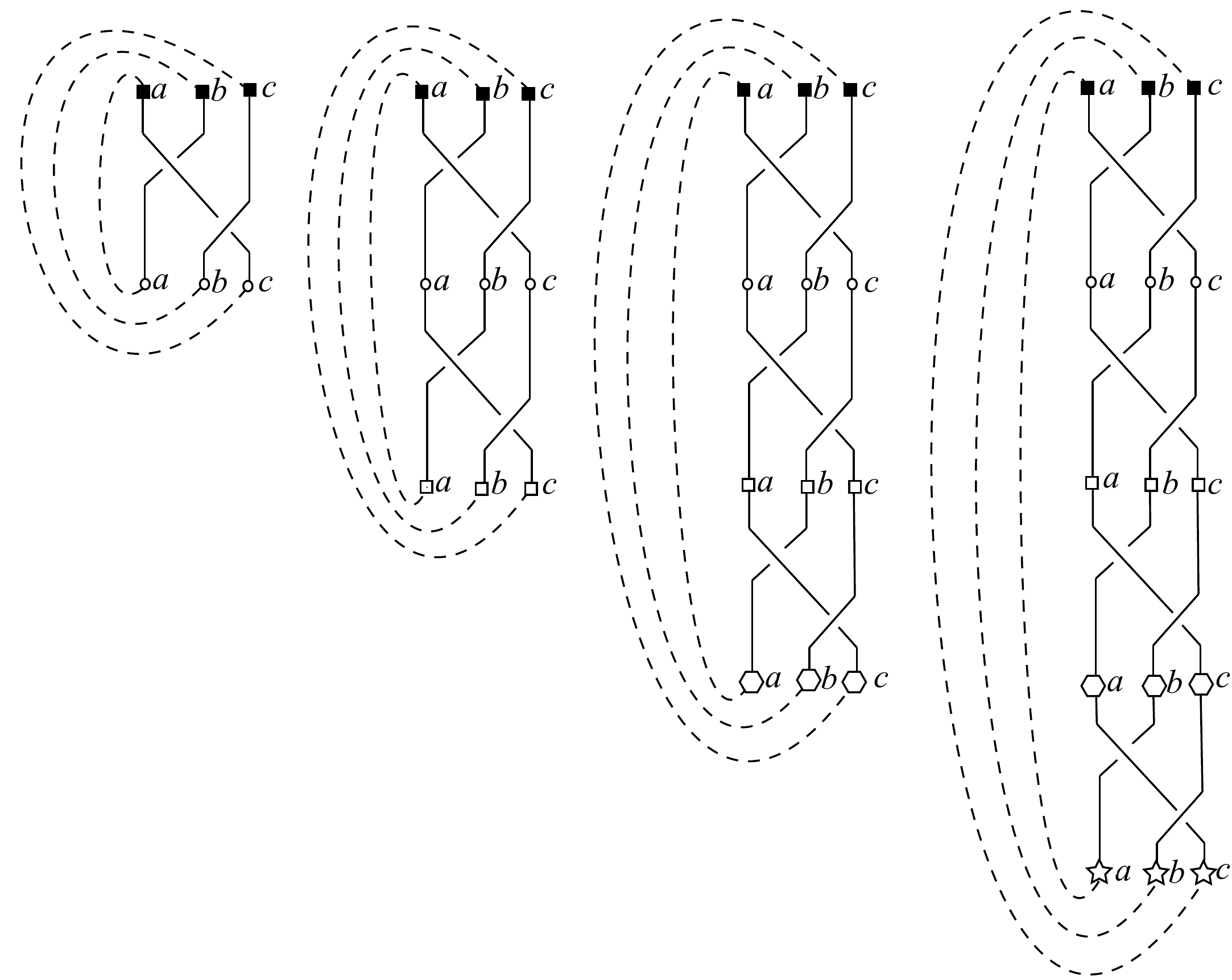}
\caption{$TH(3,1), TH(3,2), TH(3,3)$ and $TH(3,4)$. We label the extremes of each $A_3$ brick.}
\label{fig19c}
\end{figure}

For $m,n\ge1$, the $(m,n)$-{\em Turk's head} link, denoted by $TH(m,n)$, is the closure of the $m$-strand braid $(A_m)^n$, see Figure \ref{fig19c}.

\begin{corollary}\label{cor;Turk'shead} Let $n\ge 1$ be an integer. Then, $TH(3,n)$ is amphichiral. 
\end{corollary}

\begin{figure}[H]
\centering
\includegraphics[width=0.65\linewidth]{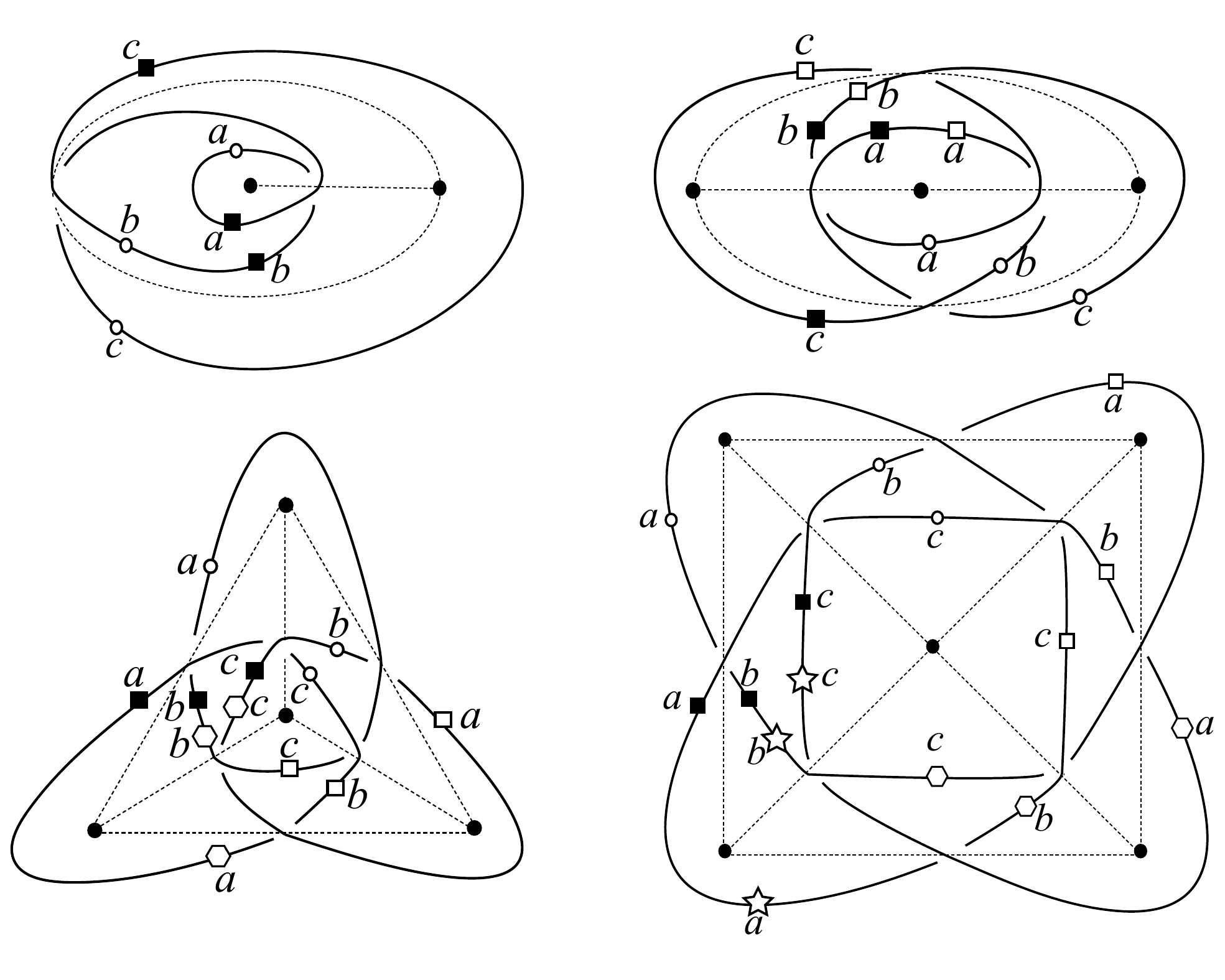}
\caption{Wheel links with labelings showing the equivalence to the corresponding Turk's head links.}
\label{fig19d}
\end{figure}

\begin{proof} It can be checked that $D(W_n,S_E^+)=TH(3,n)$, see Figure \ref{fig19d} where  the first four cases are illustrated. The result then follows by Corollary \ref{cor;w}.
\end{proof}

%\jorge{For which signature $S$ the link $L(W_n,S)$ is a torus knot?}

\begin{figure}[H]
\centering
\includegraphics[width=0.55\linewidth]{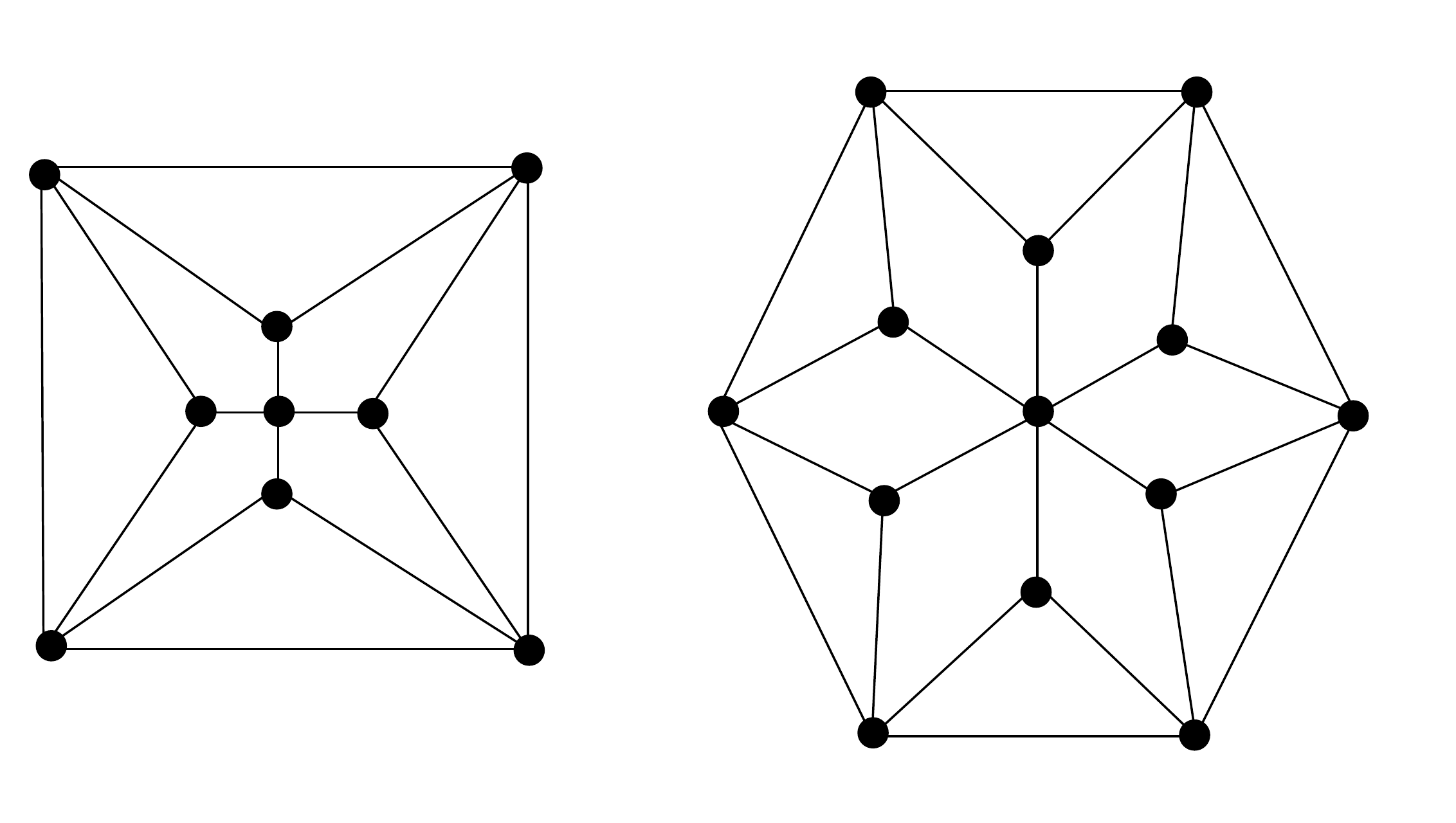}
\caption{4-ear and 6-ear.}
\label{fig20aa}
\end{figure}

Let $n\ge 3$ be an integer. The {\em $n$-ear}, denoted by $E_n$, is the graph consisting
of a $n$-cycle with an ear added on each edge and a center joined to each ear, see Figure \ref{fig20aa}.

%Let $n\ge 4$ be an even integer. Let $L(E_n,S_E)$ be the family of link diagrams, called {\em $n$-ear links}, having $E_n$ as Tait's graph and an arbitrary edge-signature $S_E$.  

\begin{corollary}\label{cor;e} Let $n\ge 4$ be an even integer. Let $L(E_n,S_E)$ be a map such that $S_E(v)=S_E(-v)$ for each pair of antipodal vertices of $med({E_n})$. Then, $L(E_n,S_E)$ is amphichiral. 
\end{corollary}

\begin{proof} By \cite[Proposition 2]{MRAR1}, $E_n$ is an antipodally self-dual map and thus, by Lemma \ref{lem:key}, $med(E_n)$ is a 2-antipodally symmetric map realized by, say $\alpha_2$. Each pair of antipodal vertices of $med({E_n})$ have the same sign and,  that is, $\alpha_2$ is sign-preserving. Moreover, by  Remark \ref{rem;color}, $\alpha_2$ is color-reversing with respect to $med(E_n)$. The result follows by Corollary \ref{cor:neg-amphi} (b). 
\end{proof}

%By Remark \ref{rem;color}, $med({E_n})$ admits a reversing 2-coloring for any even integer $n\ge 4$, see Figure \ref{fig19bb} for the case $n=4$.

%\begin{figure}[H]
%\centering
%\includegraphics[width=0.42\linewidth]{Medial4-ear}
%\caption{Reversing 2-colorings of $med({E_4})$ where antipodal faces are given by $i$ and $i^*$.}
%\label{fig19bb}
%\end{figure}

\begin{figure}[H]
\centering
\includegraphics[width=0.55\linewidth]{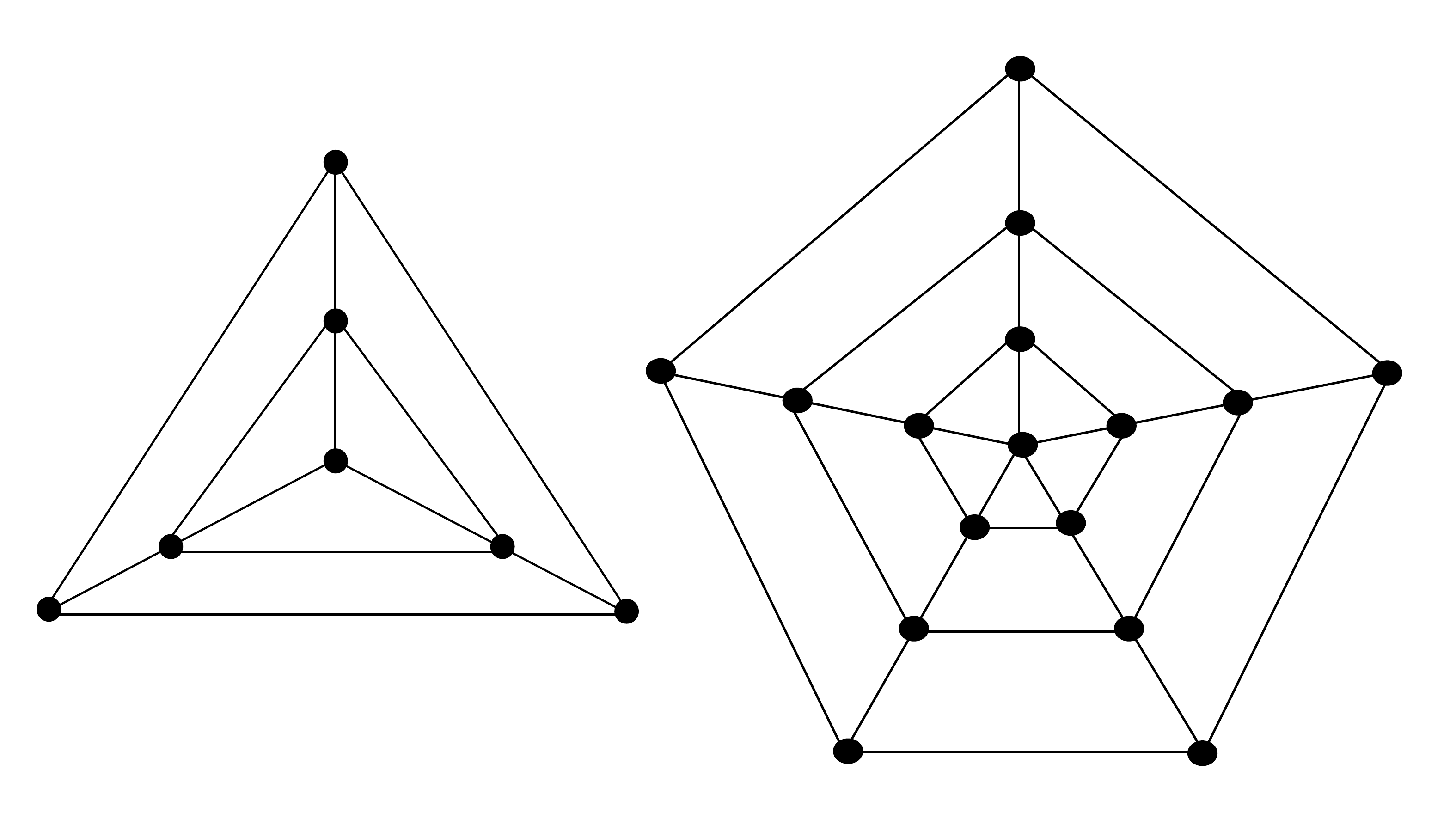}
\caption{$P_3^2$ and $P_5^3$.}
\label{fig20aaa}
\end{figure}

Let $n\ge 3$ and $\ell\ge 1$ be integers. The $(n,\ell)$-pancake, denoted by
$P^\ell_n$, is the graph consisting of $\ell$ cycles $\{v_1^1,\dots ,v_n^1\},\dots , \{v_1^\ell,\dots ,v_n^\ell\}$, a vertex $v_0$ and edges $\{v_i^{j-1},v_i^{j}\}$ for each $j=1,\dots , n$ and all $i$, see Figure \ref{fig20aaa}.

%Let $n\ge 3$ and $\ell\ge 1$ be integers with $n$ odd. Let $L(P^\ell_n, S_E)$ be the family of link diagrams, called {\em $(n,\ell)$-pancake links}, having $P_n^\ell$ as Tait's graph and an arbitrary edge-signature $S_E$.  

\begin{corollary}\label{cor;pan} Let $n\ge 3$ and $\ell\ge 1$ be integers with $n$ odd. Let $L(P_n^\ell,S_E)$ be a map such that $S_E(v)=S_E(-v)$ for each pair of antipodal vertices of $med(P_n^\ell)$. Then,$L(P_n^\ell,S_E)$ is amphichiral. 
\end{corollary}

\begin{proof} By \cite[Proposition 3]{MRAR1}, $P_n^\ell$ is an antipodally self-dual map and thus, by Lemma \ref{lem:key}, $med(P_n^\ell)$ is a 2-antipodally symmetric map realized by, say $\alpha_2$. Each pair of antipodal vertices of $med(P_n^\ell)$ have the same sign and,  that is, $\alpha_2$ is sign-preserving. Moreover, by  Remark \ref{rem;color}, $\alpha_2$ is color-reversing with respect to $med(P_n^\ell)$. The result follows by Corollary \ref{cor:neg-amphi} (b). 
\end{proof}

%By Remark \ref{rem;color}, $med({P_n^\ell})$ admits a reversing 2-coloring for any odd integer $n\ge 3$ and any $\ell\ge 1$, see Figure \ref{fig19c} for the case $n=4$.

%\begin{figure}[H]
%\centering
%\includegraphics[width=0.44\linewidth]{medialP52}
%\caption{Reversing 2-colorings of $med({P_5^2})$ where antipodal faces are given by $i$ and $i^*$.}
%\label{fig19c}
%\end{figure}

Let $K_1\# K_2$ denote the {\em sum} of knots $K_1$ and $K_2$.

\begin{corollary}\label{cor;sum} Let $K$ be a knot. Then, $K\# K^*$ is amphichiral.
\end{corollary}

\begin{proof} Let $K=D(G,S_E)$ be the knot diagram obtained from some edge-signed map $(G,S_E)$. We have that $K^*=D(G,-S_E)=D(G^*,S_E)$.  We notice that $K\# K^*=D(G \diamond G^*, S_E\cup S_E)$. Now, by \cite[Theorem 3]{MRAR1}, $G \diamond G^*$ is antipodally self-dual thus, by Lemma \ref{lem:key}, $med(G \diamond G^*)$ is a 2-antipodally symmetric map realized by, say $\alpha_2$. Each pair of antipodal vertices of $med(G \diamond G^*)$ have the same sign and,  that is, $\alpha_2$ is sign-preserving. Moreover, by  Remark \ref{rem;color},  $\alpha_2$ is color-reversing with respect to $med(G \diamond G^*)$. The result follows by Corollary \ref{cor:neg-amphi} (b). 
\end{proof}

\begin{remark}\label{rem:approach} There are amphichiral knots that are not detected by Theorem \ref{thm:antipodalR3} (when the Tait graph of the diagram is not antipodally self-dual). For instance, $L(W_2, S_E^+)$ turns out to be the Figure-eight knot which is well-known to be amphichiral however the 2-wheel $W_2$ is not antipodally self-dual map (and thusTheorem \ref{thm:antipodalR3} cannot be applied). 
\end{remark}

In Section \ref{sec;gamma}, we propose a different method to detect amphichirality. This new approach shows, in particular, the amphichirality of the Figure-eight knot (see first example in Figure \ref{fig31}).

\section{Antipodal symmetry in $\sth$}\label{sec;S3sym}

Let us recall that $\mathbb{S}^3$ can be  thought as the {\em 1-point compactification} of $\mathbb{R}^3$, that is, we take $\mathbb{R}^3$ and an additional point denoted by $\infty$. By using the {\em stereographic projection}, $\phi_n: \mathbb{S}^n \rightarrow \mathbb{R}^n$, it can be showed that $\mathbb{R}^3\cup\{\infty\}$ is equivalent to $\mathbb{S}^3$.  
Indeed, the stereographic projection sends the {\em South pole} $(0,0,0,-1)$ of $\sth$ to $(0, 0, 0)$, the {\em Equator} of $\sth$ to the unit 2-sphere, the Southern hemisphere to the region inside the unit 2-sphere, and the Northern hemisphere to the region outside of it.  The stereographic projection is not defined at the projection point $(0, 0, 0, 1)$ (the {\em North pole} of $\sth$) and the closer a point in $\sth$ is to $(0,0,0,1)$, the more distant its image is from $(0,0,0)$ in the space. For this reason, we may speak of $(0,0,0,1)$ as mapping to {\em infinity} in the 3-dimensional space, and $\sth$ as completing $\mathbb{R}^3$ by adding a point at infinity.

%\begin{observation}\label{obs;lift} 
%Let $G$ be a antipodally self-dual map and let $med(G)$ be the corresponding antipodally symmetric map of  $med(G)$ in $\mathbb{S}^2$ (obtained by \eqref{lem:key}). By construction, the map $med(G)$ is centrally symmetric in $\mathbb{R}^3$. We suppose that the map $med(G)$ is embedded in the unit 2-sphere $\mathbb{S}$ and thus inducing a central symmetric embedding of map $med(G)$ in $\mathbb{S}^3$ by taking the inverse of the stereographic projection lifting the unit 2-sphere $\mathbb{S}$ together with the map $med(G)$ into the equator of $\mathbb{S}^3$.  
%\end{observation}
\smallskip

 We say that $L$ is 3-antipodally symmetric if there is $\hat L\in [L], \hat L\in \mathbb{S}^3$ such that $\alpha_3(\hat L)=\hat L$.
 
\begin{theorem} \label{thm:antipodalS3} Let $(G,S_E)$ be an edge-signed self-dual map and suppose that $med(G)$ is a $2$-antipodally symmetric map (realized by $\alpha_2$). If either
\smallskip

a) $\alpha_2$ is color-preserving and sign-preserving  or 

b) $\alpha_2$ is color-reversing and sign-reversing
\smallskip

\noindent then $L(G,S_E)$ is $3$-antipodally symmetric.
\end{theorem}

\begin{proof}  
%Let $D(G,S_E)$ be the link diagram induced by $(med(G),C_F,S_V)$.  
Since $med(G)$ is 2-antipodally symmetric then it is also 3-centrally symmetric. Consider the embedding $\embed(G,S_E)$. It can be checked that if $\alpha_2$ is either (a) color-preserving and sign-preserving or (b) color-reversing and sign-reversing then the piece of arc of the diagram passing over (resp. passing under) in crossing $v$ correspond to the piece of arc of the diagram passing under (resp. passing overer) in crossing $\alpha_2(v)$ in $\embed(G,S_E)$ see Figure \ref{fig:2antisym} (center). 

\begin{figure}[H]
\centering
\includegraphics[width=\linewidth]{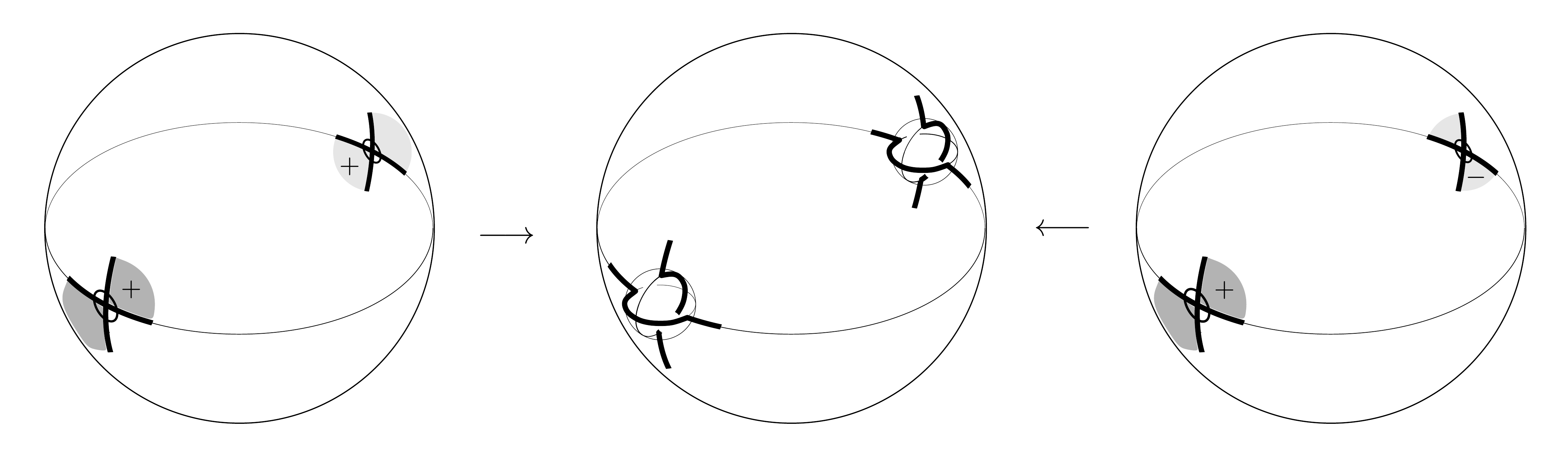}
\caption{(Left) Antipodal pair of vertices of $med(G)$ in case (a). (Right)  Antipodal pair of vertices of $med(G)$ in case (b) (Center) The local modifications around the antipodal pair of vertices following the color-sphere rules.}
\label{fig:2antisym}
\end{figure}

Let $H$ be the hyperplane in $\mathbb{R}^4$ containing the Equator of $\sth$.  Let $\pi$ be the stereographic projection of $\sth$ into $H$. We suppose that the Equator is the unit 2-sphere. Notice that 

$$\|\pi(x)\|\left\{\begin{array}{ll}
>1 & \text{ if } x \text{ is in the Northern hemisphere of } \sth,\\ 
=1 & \text{ if } x \text{ is in the Equator of } \sth,\\ 
<1 & \text{ if } x \text{ is in the Southern hemisphere of } \sth.\\ 
\end{array}\right.$$

We shall lift $\embed(G,S_E)$ into $\mathbb{S}^3$ by taking the inverse $\pi^{-1}$. We claim that $\pi^{-1}(\embed(G,S_E)$ gives the desired 3-antipodally symmetric embedding of $L(G,S_E)$. To see this, we define the {\em inversion} function (with respect to $\mathbb{S}^{n-1}$) as

$$\begin{array}{lllc}
i_n: & \mathbb{R}^n& \longrightarrow & \mathbb{R}^n\\ 
& x & \mapsto & \frac{x}{\|x\|^2}
\end{array}$$

Notice that 

$$\|i(x)\|\left\{\begin{array}{ll}
>1 & \text{ if } \|x\|<1,\\ 
=1 & \text{ if } \|x\|=1,\\ 
<1 & \text{ if } \|x\|>1.\\ 
\end{array}\right.$$

Let $x\in\sth$. We have that $$\pi\circ\alpha_3(x)=\pi(-x)=c_3\circ i_3\circ \pi(x).$$

We thus have that the perturbed pieces in each pair of antipodal vertices are mapped into diametrically antipodal pieces in $\mathbb{S}^3$. Indeed, for instance, let us take the piece $p$, say inside of $\stw$ around crossing $v$. We have then that $i_3\circ c_3(p)$ gives the piece $p'$ outside of  $\stw$ around crossing $\alpha_2(v)$. Therefore, $\pi^{-1}$ maps $p$ and $p'$ to the Northern and Southern hemisphere of $\sth$ respectively such that $\pi^{-1}(p)=\alpha_3\circ \pi^{-1}(p)=\pi^{-1}(p')$.
\smallskip

We thus obtain that $\pi^{-1}(\embed(G,S_E)$ is an antipodally symmetric embedding in $\sth$.

%\begin{figure}[H]
%\centering
%\includegraphics[width=0.85\linewidth]{mouvementsym3sphere}
%\caption{Two possible cases around crossings $v$ and $-v$ (a) A preserving 2-coloring and antipodal pair of vertices  with the same sign (b) Reversing 2-coloring and antipodal pair of vertices with opposite signs (we denote by $-i$ the antipodal of point $i$).}
%\label{fig19aa}
%\end{figure}

%For each crossing of $L(S,G)$, viewed from above of it (in $\mathbb{R}^3$), we push (locally) the piece of arc of the diagram passing over into $\mathbb{S}_{ext}$ and push (again locally) the piece of arc of the diagram passing under into $\mathbb{S}_{int}$. Both pushes are again performed by radial projection from the origin as done in the proof of Theorem \ref{thm:antipodalR3}, see Figure \ref{fig19b}. 

%\begin{figure}[H]
%\centering
%\includegraphics[width=.45\linewidth]{localperturbationS3}
%\caption{Local perturbations around a crossing and its antipodal.}
%\label{fig19b}
%\end{figure}

%\begin{observation}
%Notice that the corresponding pushed pieces of arc are not longer symmetric in $\mathbb{R}^3$. Indeed, in this case the arcs $[b,d]$ and $[-c,-a]$ are pushed to $S_{ext}$ while the arcs $[-b,-d]$ and $[c,a]$ are pushed to $S_{int}$.
%\end{observation}
 
\end{proof}

The following result is easily derived from Theorem \ref{thm:antipodalS3} (a).

\begin{corollary}\label{cor:antipodal3sph} Let $(G,S_E^+)$ be an edge-signed map and suppose that $med(G)$ is a 2-antipodally symmetric map (realized by $\alpha_2$). If $\alpha_2$ is color-preserving then $L(G,S_E^+)$ is 3-antipodally symmetric.
\end{corollary}

A straightforward consequence of the above corollary, we obtain that the Hopf link is 3-antipodally symmetric, see Figure \ref{fig16}.

\begin{corollary}\label{cor;wantipol3} Let $n\ge 1$ be an odd integer and let $L(W_n,S_E)$ be a map such that $S_E(v)=-S_E(-v)$ for each pair of antipodal vertices of $med({W_n})$. Then, $L(W_n,S_E)$ is 3-antipodally symmetric. 
\end{corollary}

\begin{proof} By \cite[Proposition 1]{MRAR1}, $W_n$ is an antipodally self-dual map and thus, by Lemma \ref{lem:key}, $med(W_n)$ is a 2-antipodally symmetric map, realized by, say $\alpha_2$. By hypothesis, each pair of antipodal vertices of $med({W_n})$ have opposite signs and thus $\alpha_2$ is sign-reversing. Moreover,  by  Remark \ref{rem;color}, $\alpha_2$ is color-reversing with respect to $med(W_n)$. The result follows by Theorem \ref{thm:antipodalS3} (b). 
\end{proof}

\begin{corollary}\label{cor;Torus} Let $n\ge 2$ be an odd integer. Then, the {\em torus knot} $T(2,n)$ is 3-antipodally symmetric. 
\end{corollary}

\begin{proof} Let $S'_E$ be the edge-signature of $W_n$ where the edges of the exterior cycle have $+$ sign and edges incident with the center have $-$ sign. It can be checked that the diagram $D(W_n,S'_E)$ is the same as the one for $T(2,n)$. 
The result then follows by Corollary \ref{cor;wantipol3}.
\end{proof}

\begin{corollary}\label{cor;sumS3} Let $K$ be a knot. Then, $K\# K$ is 3-antipodally symmetric.
\end{corollary}

\begin{proof} Let $K=D(G,S_E)$ for some edge-signed map $(G,S_E)$. 
We have that $D(G,S_E)=D(G^*,-S_E)$.  
We notice that $K\# K=D(G \diamond G, S_E\cup -S_E)$. Now, by \cite[Theorem 3]{MRAR1}, $G \diamond G^*$ is antipodally self-dual thus, by Lemma \ref{lem:key}, $med(G \diamond G^*)$ is a 2-antipodally symmetric map realized by, say $\alpha_2$. Each pair of antipodal vertices of $med(G \diamond G^*)$ have the opposite signs,  that is, $\alpha_2$ is sign-reverserving. Moreover, by  Remark \ref{rem;color},  $\alpha_2$ is color-reversing with respect to $med(G \diamond G^*)$. The result follows by Theorem \ref{thm:antipodalS3} (b). 
\end{proof}

%%%%%%%%%%%%%%%%%%%%%%%%%%%%%%%%%%%%%%%%%%%%%
\section{$\Gamma$-curve}\label{sec;gamma}

Let $G=(V,E,F)$ be a map. Let $\Gamma$ be a non self-intersecting curve in the plane starting at one vertex in $V(G)$ verifying the following conditions :
\smallskip

a) $\Gamma$ may passes throughout vertices,

b) $\Gamma$ may contain edges (if it contains an edge it also contains its extreme vertices). $\Gamma$ does not intersect edges otherwise.
\smallskip

We have that $\Gamma$ is composed of paths and intervals of consecutive vertices, called {\em vertex-interval}. If $\Gamma$ has not vertex-interval then $\Gamma$ is just a cycle of $G$, called {\em $\Gamma$-cycle}. It will be called {\em $\Gamma$-noncycle} otherwise.

\begin{figure}[H]
\centering
\includegraphics[width=0.2\linewidth]{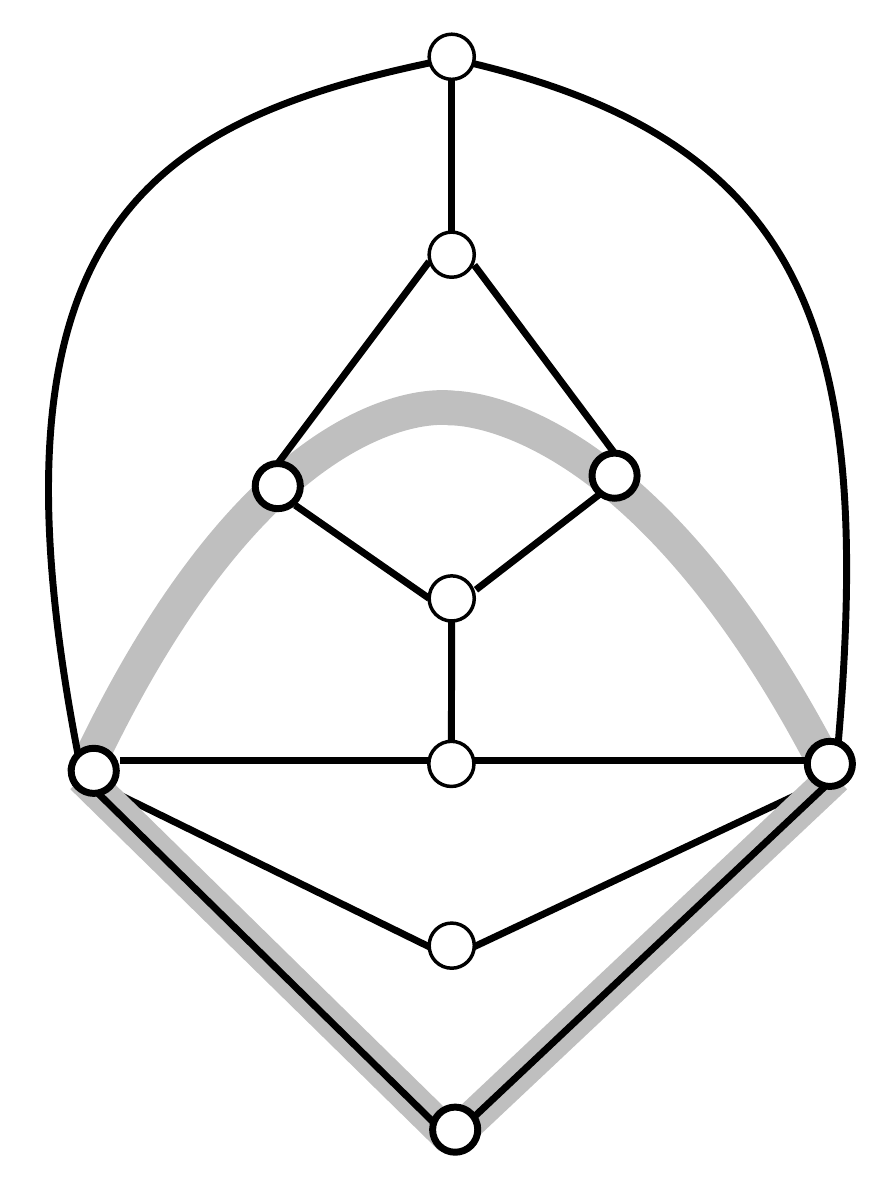}
\caption{A  $\Gamma$-noncycle (in thick gray) of $G$ consisting of two edges and a vertex-interval formed by four consecutive vertices.}
\label{fig25a}
\end{figure}

\begin{figure}[H]
\centering
\includegraphics[width=0.45\linewidth]{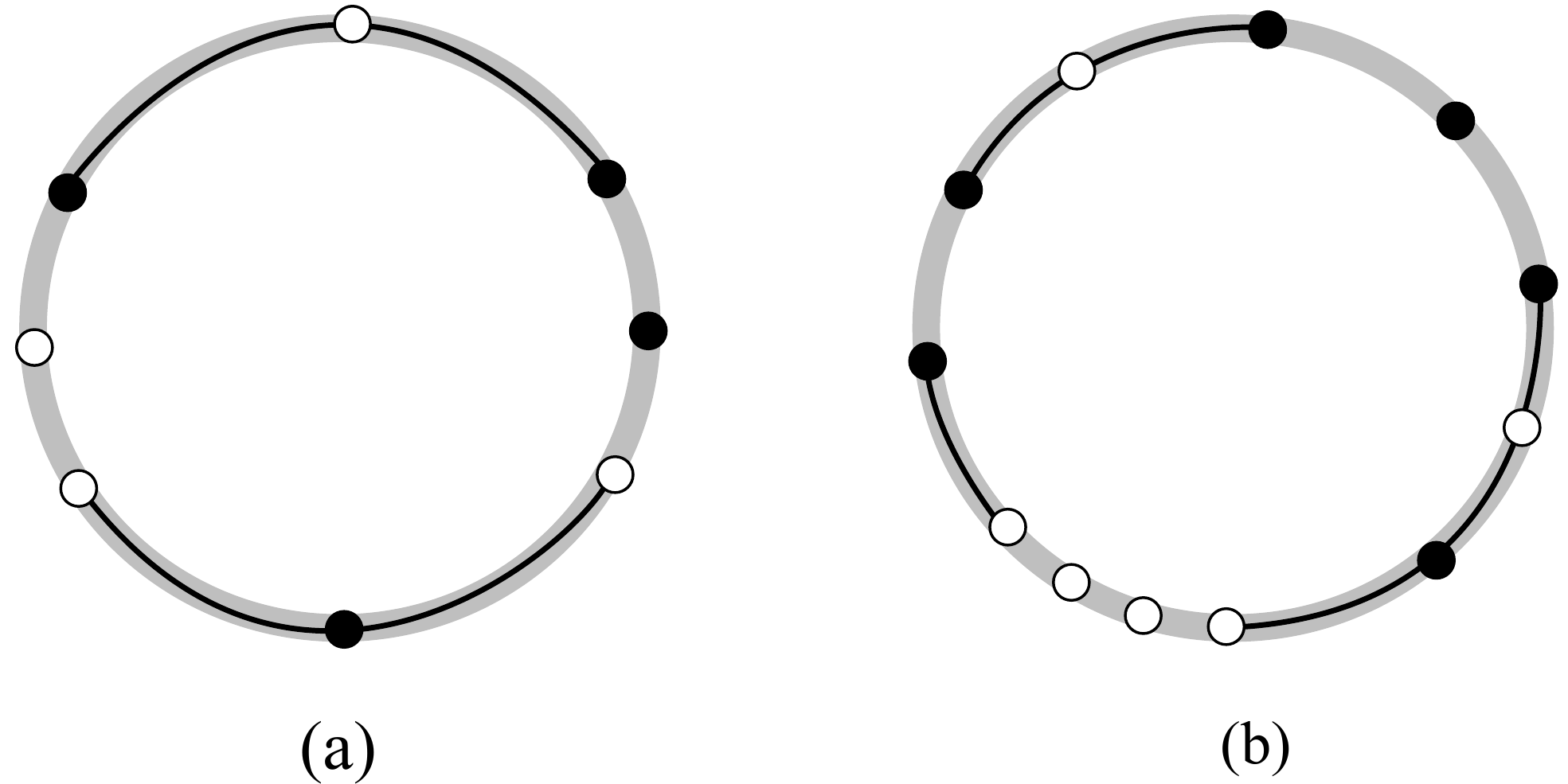}
\caption{(a) An even nonadequate $\Gamma$-noncycle consisting of two paths (each of length 2) and of two vertex-interval formed by three consecutive vertices (none of them monochromatic) (b) An even adequate  $\Gamma$-noncycle consisting of three paths (two of length 2 and one of length 1) and three vertex-interval (all monochromatic), formed by two consecutive black vertices, three consecutive black vertices and four consecutive white vertices).}
\label{fig25b}
\end{figure}

We say that $\Gamma$ is {\em even} if it contains an even number of vertices of $G$. If $V(G)$ is 2-colored, we say that $\Gamma$ is {\em adequate} if the vertices, in each interval of consecutive vertices are monochromatic (they have the same color). We notice that the color of the vertices in each path contained in $\Gamma$ alternates, see Figure \ref{fig25b}.    

For a curve $\Gamma$ in $G$, we denote by $V_{int(\Gamma)}$ (resp. $V_{ext(\Gamma)}$) the set of vertices in the interior (resp. exterior) of $\Gamma$. Let $G_{int(\Gamma)}=[V\setminus V_{ext(\Gamma)}]$ (resp. $G_{ext(\Gamma)}=[V\setminus V_{int(\Gamma)}])$, that is, $G_{int(\Gamma)}$ (resp. $G_{ext(\Gamma)}$) is the the graph induced by vertices $V\cup V_{int(\Gamma)}$ (resp. $V\cup V_{ext(\Gamma)}$).  
\smallskip

%Since $\Gamma_G$ does not intersect edges then $V(\Gamma_G)$ is a {\em cut-set of $G$}, that is, $G\setminus V(\Gamma_G)$ consists of two connected components, say $\Gamma_G^{int}$ (inside $\Gamma_G$) and $\Gamma_G^{ext}$ (outside $\Gamma_G$). 

%\section{Symmetric cycle}\label{sec:symcycle}

%Let $(G,C_V,S_F)$ be a vertex-colored face-signed map. For a cycle $\Gamma$ in $G$, we denote by $V_N$ (resp. $V_S$) be the set of vertices inside (resp. outside ) $\Gamma$. Let $G_N=[V\setminus V_S]$ (resp. $G_S=[V\setminus V_N]$), that is, $G_N$ (resp. $G_S$) is the the graph induced by vertices $V\cup V_N$ (resp. $V\cup V_S$).  

A curve $\Gamma$ of $G$ is said to be {\em symmetric} if there is a bijection
$$\sigma : V(G_{int(\Gamma)}) \longrightarrow V(G_{ext(\Gamma)})$$ 

such that $\{u,v\}\in E(G_{int(\Gamma)})$ if and only if $\{\sigma(u),\sigma(v)\}\in E(G_{ext(\Gamma)})$
with $\sigma(\Gamma)=\Gamma$. 

In other words, $\Gamma$ is symmetric if there is an isomorphism $\sigma$ between $G_{int(\Gamma)}$ and $G_{ext(\Gamma)}$ where $\Gamma$ is fixed. We notice that $\sigma$ can be thought as an automorphism of $G$ with $\sigma |_\Gamma=\Gamma$. In this case, $\sigma$ is color-preserving (resp. color-reversing) if each pair of vertices $v\in V_{int(\Gamma)}$ and $\sigma(v)\in V_{ext(\Gamma)}$ have the same (resp. opposite) color and it is sign-preserving (resp. sign-reversing) if each pair of faces $f\in F(G_{int(\Gamma)})$ and $\sigma(f)\in F(G_{ext(\Gamma)})$ have the same (resp. different) sign. 

%and $\sigma(G_{int(\Gamma)})=G_{ext(\Gamma)}$.
\begin{figure}[H]
    \centering
    \includegraphics[width=0.6\textwidth]{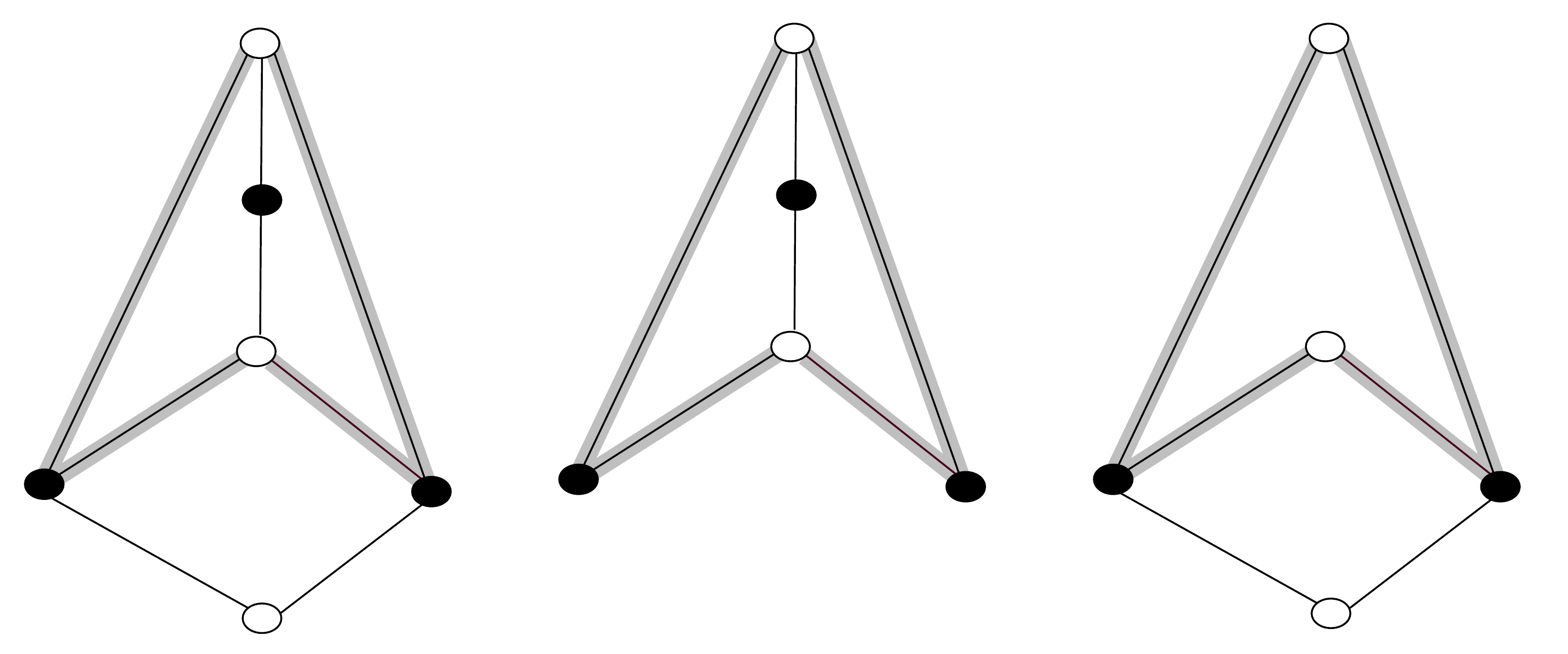}
    \caption{(From left to right) Graph $G$ with a color-reversing symmetric $\Gamma$-cycle (in gray), $G_{int(\Gamma)}$ and $G_{ext(\Gamma)}$.}
    \label{fig26}
\end{figure}

\begin{figure}[H]
    \centering
    \includegraphics[width=.8\textwidth]{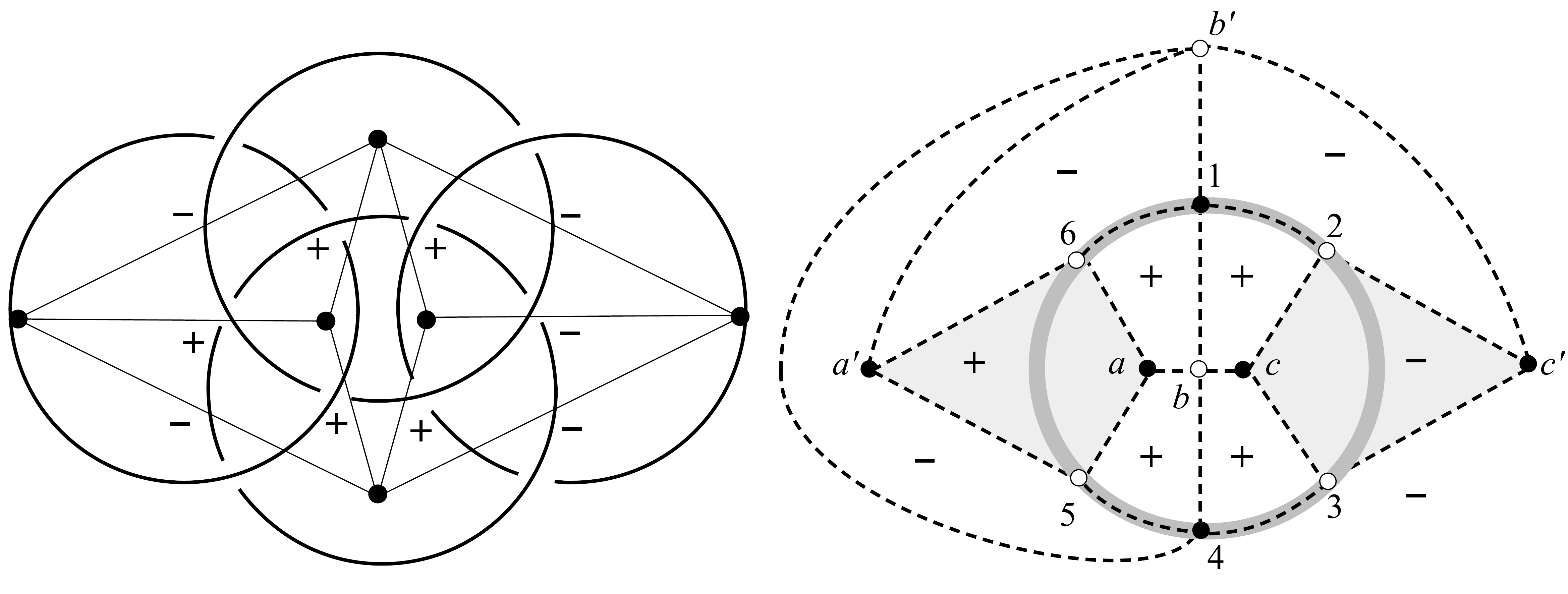}
    \caption{(Left) Edge-signed graph map $(G,S_E)$ together with $D(G,S_E)$ (Right) Face-signed incidence graph $(I(G),C_V,S_F)$ with an even adequate sign-reversing curve $\Gamma$ (gray) with $\sigma(a)=a', \sigma(b)=b', \sigma(c)=c'$ and $\sigma(\Gamma)=\Gamma$ (shared faces are shaded).} 
    \label{fig33a}
\end{figure}

Let us suppose that $\Gamma$ is an even adequate curve of $(I(G),C_V,S_F)$.  Let $I_\Gamma$ be an interval of vertices in $\Gamma$. Since the vertices in $I_\Gamma$ are monochromatic and the faces of $I(G)$ are squares then two consecutive vertices in $I_\Gamma$ belong to a face with two edges in $G_{int(\Gamma)}$ and two edges in $G_{ext(\Gamma)}$.  We call such face a {\em shared face} with respect to $\Gamma$, see Figure \ref{fig33a}.

%We notice that if $\Gamma_\sigma$ is symmetric then the existence of at least one shared face clearly forces $\Gamma_\sigma$ to be color-preserving.
%\smallskip

\begin{figure}[H]
    \centering
    \includegraphics[width=.85\textwidth]{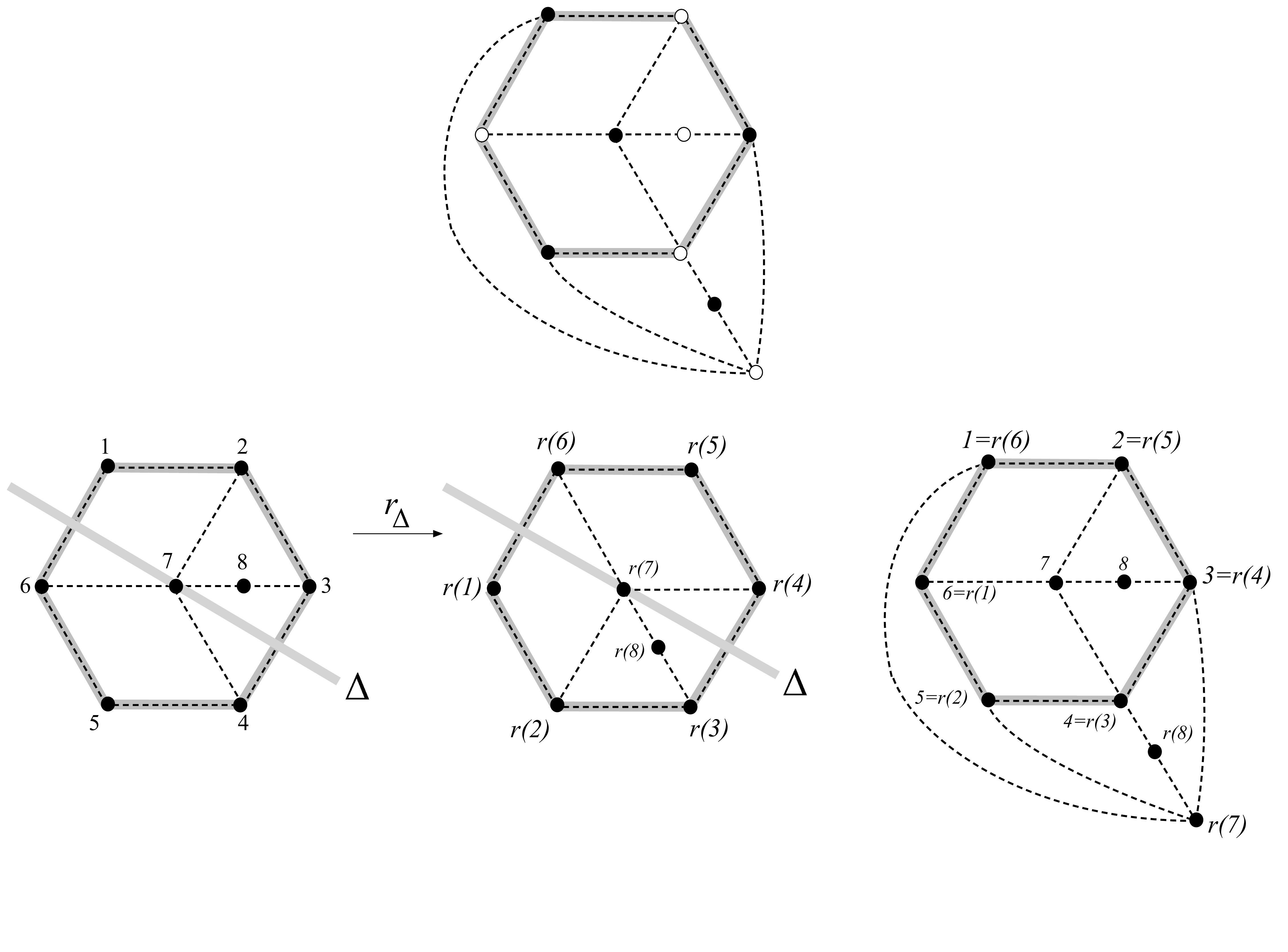}
    \caption{(Top) Bicolored incidence graph $I(G)$ (dashed edges) with a symmetric $\Gamma$-cycle (in gray). (Bottom)  A drawing of $I(G)_{int(\Gamma)}$ with a reflection line $\Delta$, a drawing of $I(G)_{ext(\Gamma)}$ (given by $r_\Delta$) and the graph obtained by identifying $\Gamma$ (in this case we obtain $G$ and thus $r_\Delta$ induces an automorphism $\sigma(I(G))$.} \label{fig28}
\end{figure}

Since $\sigma$ keeps $\Gamma$ invariant then we should have that $\sigma|_\Gamma$ is either a rotation or a reflection of $\Gamma$. Let $\Gamma$ be symmetric curve of a graph $G$. Suppose that we draw $G_{int(\Gamma)}$ in the plane with $\Gamma$ as regular polygon. Let $r_\Delta$ be the reflection of $G_{int(\Gamma)}$ with respect to a line $\Delta$. Since $r_\Delta(G_{int(\Gamma)})=G_{int(\Gamma)}$ and $\Gamma$ is symmetric then $r_\Delta(G_{int(\Gamma)})$ can be thought as $G_{ext(\Gamma)}$.
The latter induce an automorphism $\sigma(G)$ if the graph obtained by identifying the two regular polygones (together with their labels) gives $G$, see Figure \ref{fig28}.

\begin{figure}[H]
    \centering
    \includegraphics[width=.85\textwidth]{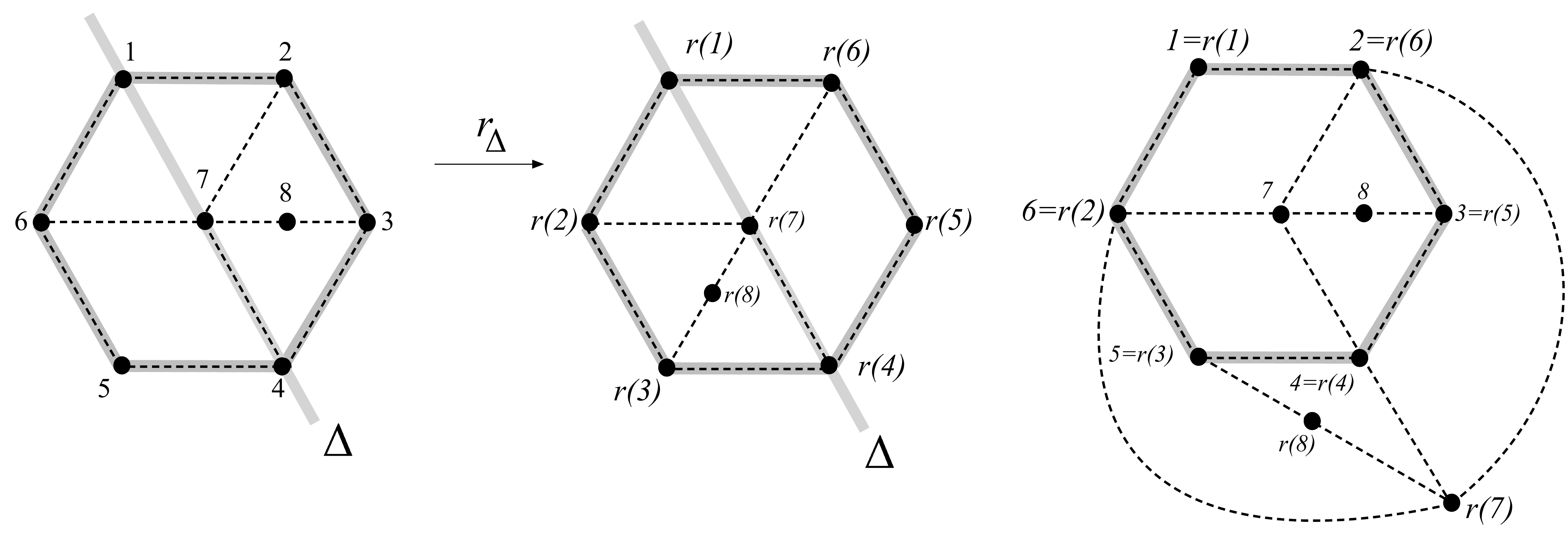}
    \caption{A drawing of $I(G)_{int(\Gamma)}$ (of the incident graph of Figure \ref{fig28}) with a reflection line $\Delta$, a drawing of $I(G)_{ext(\Gamma)}$ (given by $r_\Delta$) and the graph obtained by identifying $\Gamma$. In this case such a graph is not isomorphic to $I(G)$.
%$\sigma_1(r_\Delta(G_{int(\Gamma)}))=G_{ext(\Gamma)}$ where $\sigma_1(1)=6, \sigma_1(2)=5, \sigma_1(3)=4, \sigma_1(4)=3, \sigma_1(5)=2, \sigma_1(6)=1, \sigma_1(7)=10$ and $\sigma_1(8)=9$. \\
} \label{fig288}
\end{figure}

We notice that there might be reflections of a symmetric $\Gamma$ curve of  graph $G$ not inducing an automorphism of $G$, see Figure \ref{fig288}.

\begin{figure}[H]
    \centering
    \includegraphics[width=.5\textwidth]{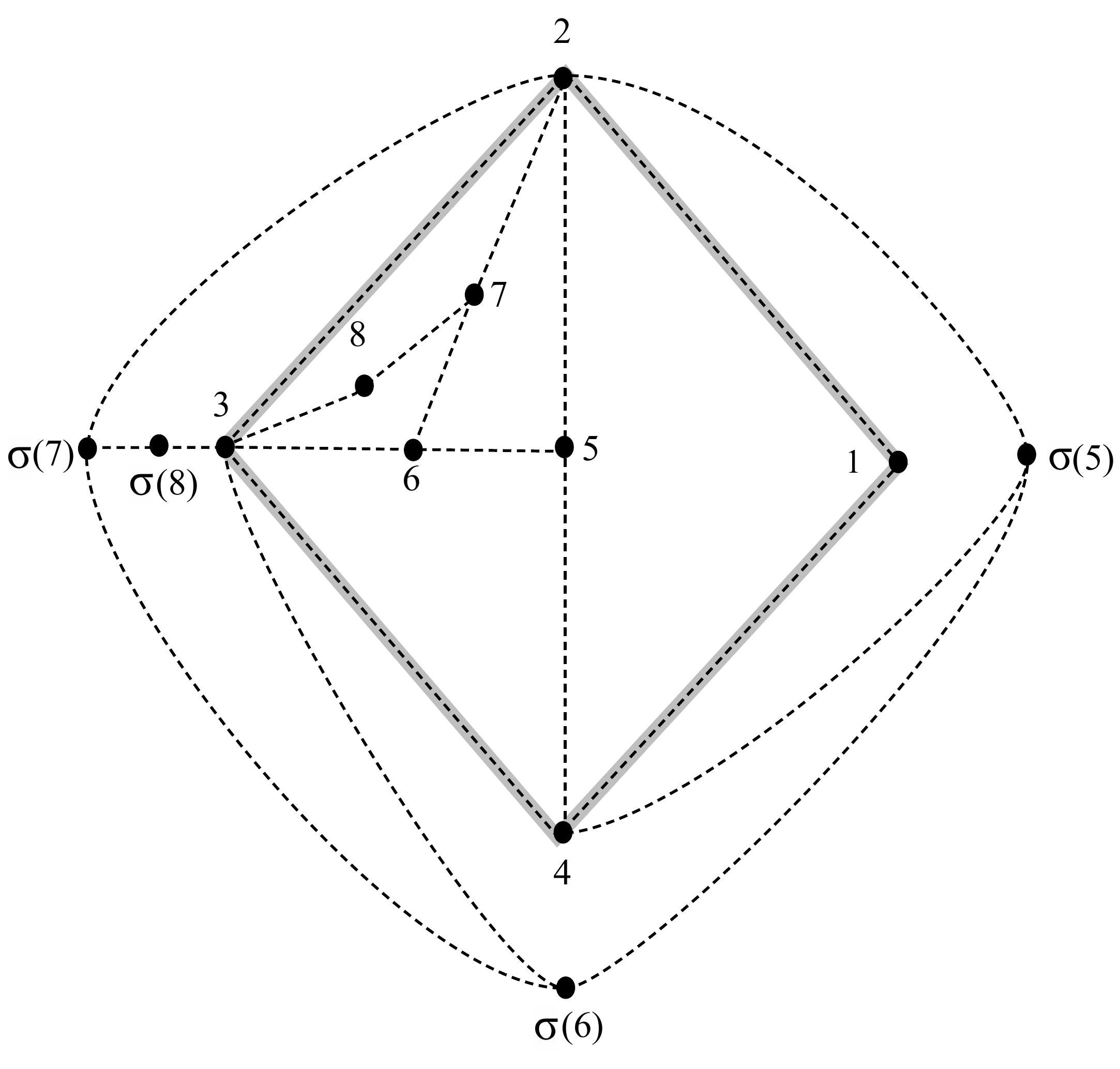}
    \caption{A symmetric nonreflexive $\Gamma$ cycle.} \label{fig28ax}
\end{figure}

We say that $\Gamma$ is {\em reflexive} if $\sigma$ arises on this way. In this case, we have that either  $\Delta$ cuts twice $\Gamma$ between two vertices (either on an edge or between two vertices in a vertex-interval) or $\Delta$ passes either through two vertices or through one vertex and a shared face. We notice that there are graphs admitting a symmetric $\Gamma$ curve which is not reflexive, see Figure \ref{fig28ax}.

%we can pass from each other with a reflexion with respect to a line $\Delta$. 
%We say that $\Gamma_\sigma$ is {\em color-preserving} (resp. {\em color-reversing}) if  each pair of vertices $v\in V_{int(\Gamma)}$ and $\sigma(v)\in V_{ext(\Gamma)}$ have the same (resp. opposite) color and it is {\em sign-preserving} (resp. {\em sign-reversing}) if each pair of faces $f\in F(G_{int(\Gamma)})$ and $\sigma(f)\in F(G_{ext(\Gamma)})$ have the same (resp. different) sign.  Notice that there is not requirement on the colors for the vertices in $\Gamma$, for instance, automorphisms $\sigma_1$ and $\sigma_2$ given in Figure \ref{fig28} are color-preserving of the 2-coloring given in Figure \ref{fig28a} however they do not preserve color on the vertices of $\Gamma$. 
 
%\begin{figure}[H]
  %  \centering
 %   \includegraphics[width=.27\textwidth]{images/ExIncGraph2}
 %   \caption{Bicolored incidence graph $I(G)$.} \label{fig28a}
%\end{figure}

\begin{theorem}\label{theo:amphi} Let $(G,S_E)$ be an edge-signed map. If $(I(G),C_V,S_F)$ admits either

(a) a sign-preserving color-reversing reflexive curve or 

(b) a sign-reversing color-preserving reflexive curve

then $L(G,S_E)$ is amphichiral. 

Moreover, $L(G,S_E)$ and its mirror can be embedded such that one can be obtained from the other by performing a rotation $\pi$ degrees.
\end{theorem}

\begin{proof} (a) Let $\Gamma$ be a sign-preserving color-reversing reflexive curve of $(I(G),C_V,S_F)$ realized by $\sigma$. We shall construct an embedding of $(I(G),C_V,S_F)$ in $\stw$ from which we can pass to  $(I(G),C_V^o,S_F)$ by performing a rotation of $\pi$ degrees. Since $(I(G), C_V, S_F)$ determines $L(G,S_E)$ then, by Remark \ref{rem:I(G)},  $(I(G), C_V^o, S_F)$ determines $L(G,S_E)^*$.  Therefore, this rotation of $\pi$ degrees would give the desired orientation-preserving homeomorphism from $L(G,S_E)$ to $L(G,S_E)^*$. 
\smallskip

The construction of the embedding is in three steps, we will be illustrating these steps by applying them to the map given in Figure \ref{fig28}.
\smallskip

{\bf Step 1)} We draw $G_{int(\Gamma)}$ in the Northern hemisphere of $\stw$ by taking the drawing of $G_{int(\Gamma)}$ together with the reflection line in $\mathbb{B}^1$ with $\Gamma$ forming the vertices of a regular polygon lying in the equator of $\stw$.  We then project this to the Northern hemisphere of $\stw$, see Figure \ref{fig29a}.

\begin{center}
\begin{figure}[H]
    \includegraphics[width=.8\textwidth]{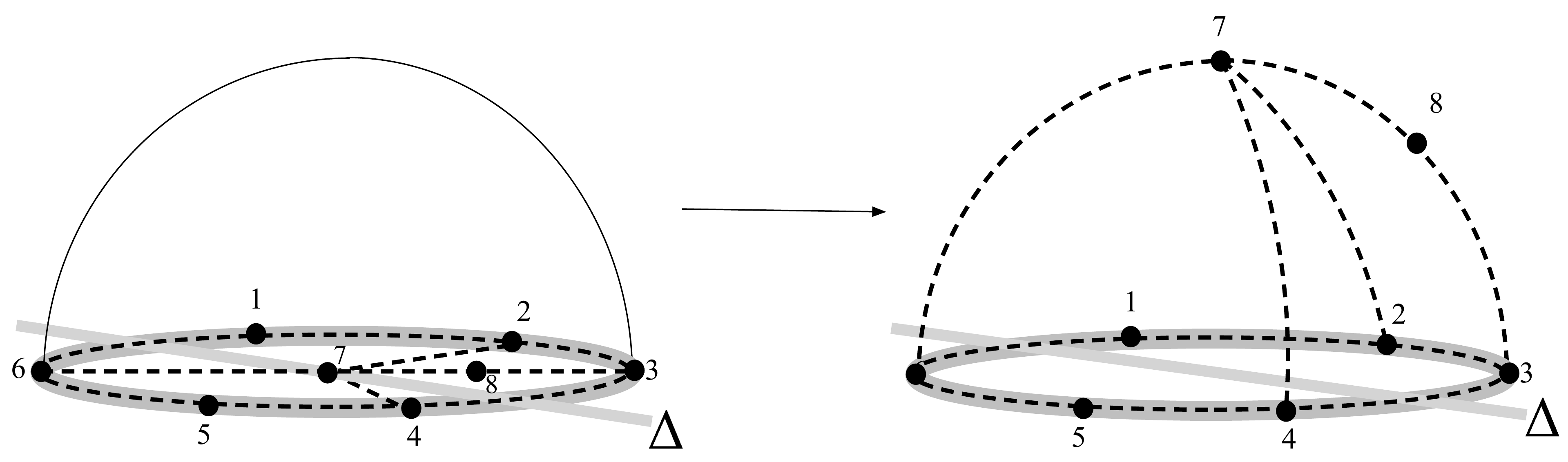}
    \caption{Step 1.}
    \label{fig29a}
\end{figure}
\end{center}

{\bf Step 2)} We draw $G_{ext(\Gamma)}$ in the Southern hemisphere by taking the drawing of $r_\Delta(G_{int(\Gamma)})$ with $r_\Delta(\Gamma)$ forming the vertices of a regular polygon lying in the equator of $\stw$.  We then project this to the Southern hemisphere of $\stw$. We notice that the same drawing can be obtained by  performing a rotation of $\pi$ degrees around $\Delta$ of the drawing of $G_{int(\Gamma)}$ in the Northern hemisphere, see Figure \ref{fig29b}.

\begin{center}
\begin{figure}[H]
    \includegraphics[width=.8\textwidth]{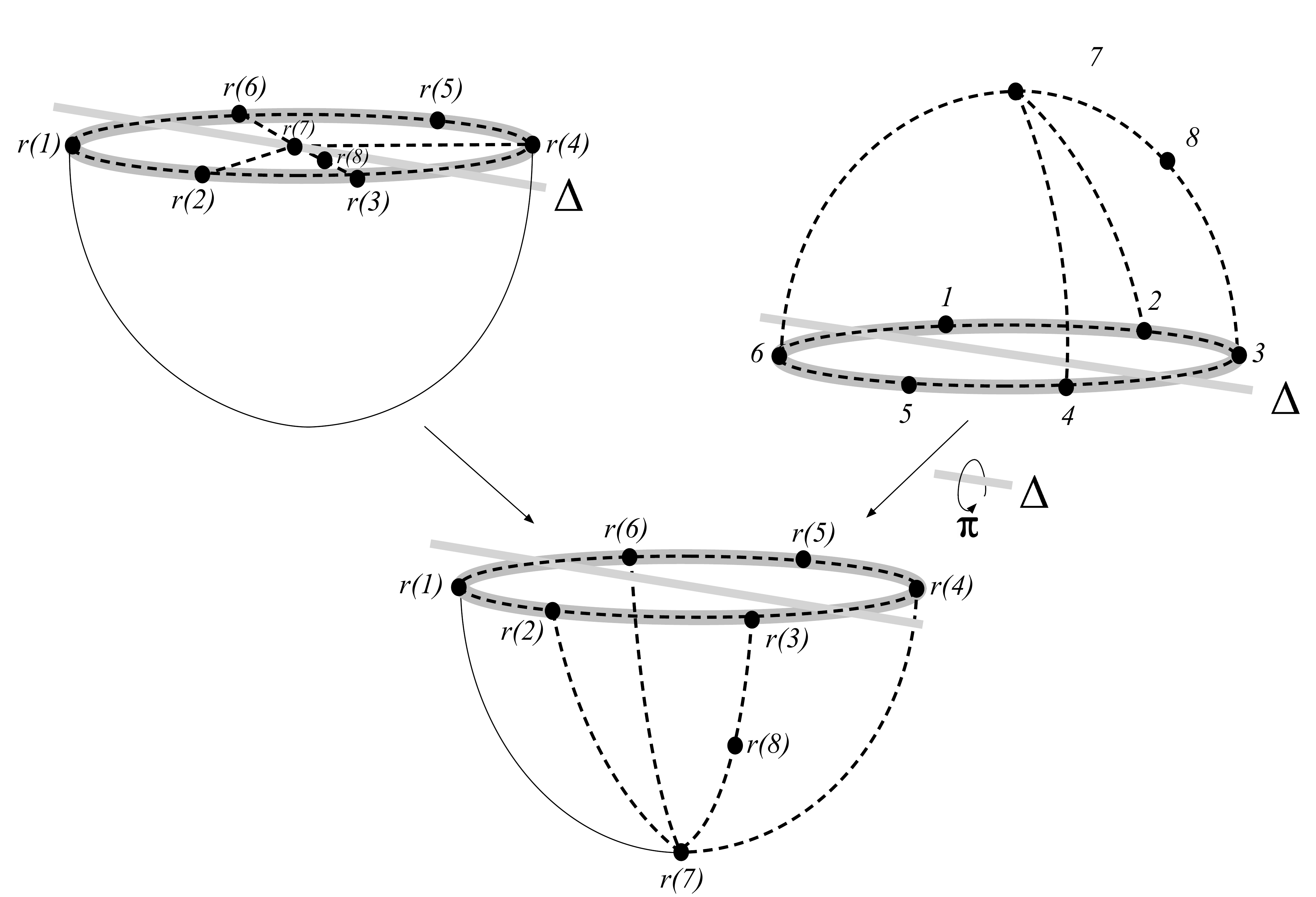}
    \caption{Step 2.}
    \label{fig29b}
\end{figure}
\end{center}

{\bf Step 3)} We color the vertices of both embeddings as given by $C_V$ and update labels according with $r_\Delta$. We finally glue together the two hemispheres. Notice that vertices in $\Gamma$ with the same labels are aligned one above the other, see Figure \ref{fig29c}.

\begin{center}
\begin{figure}[H]
    \includegraphics[width=.45\textwidth]{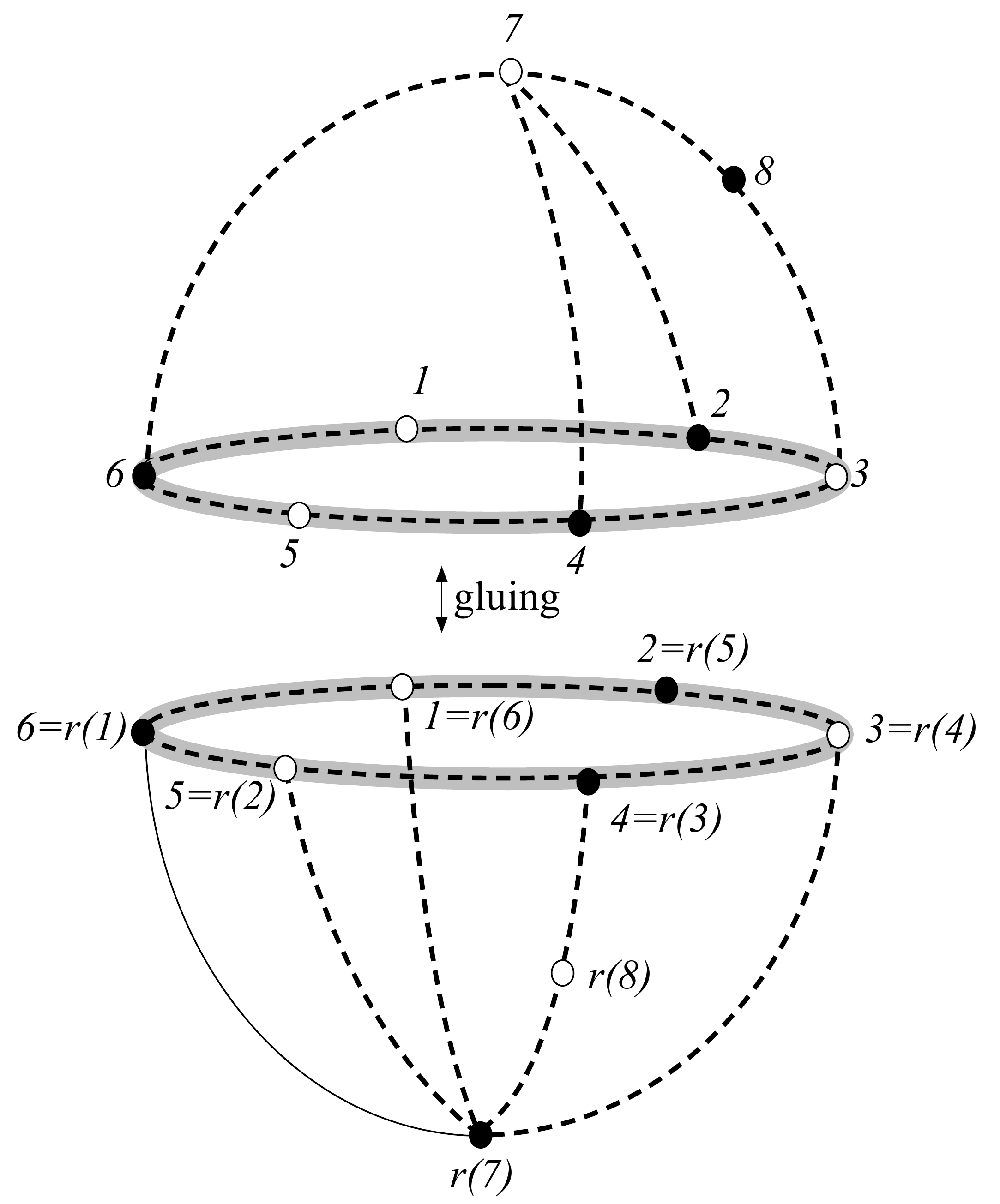}
    \caption{Step 3.}
    \label{fig29c}
\end{figure}
\end{center}

We observe that a rotation of $\pi$ degrees around $\Delta$ of the drawing of $(I(G),C_V,S_F)$ exchanges faces $f\in G_{int(\Gamma)}$ and $\sigma(f)\in G_{ext(\Gamma)}$ both having the same signs (since 
$\Gamma$ is sign-preserving) but swapping the coloring (since $\Gamma$ is color-reversing). We thus have that we can pass from $(I(G),C_V,S_F)$ to $(I(G),C_V^o,S_F)$ by applying a rotation of $\pi$ degrees, see Figure \ref{fig30}. 

\begin{figure}[H]
    \centering
    \includegraphics[width=.8\textwidth]{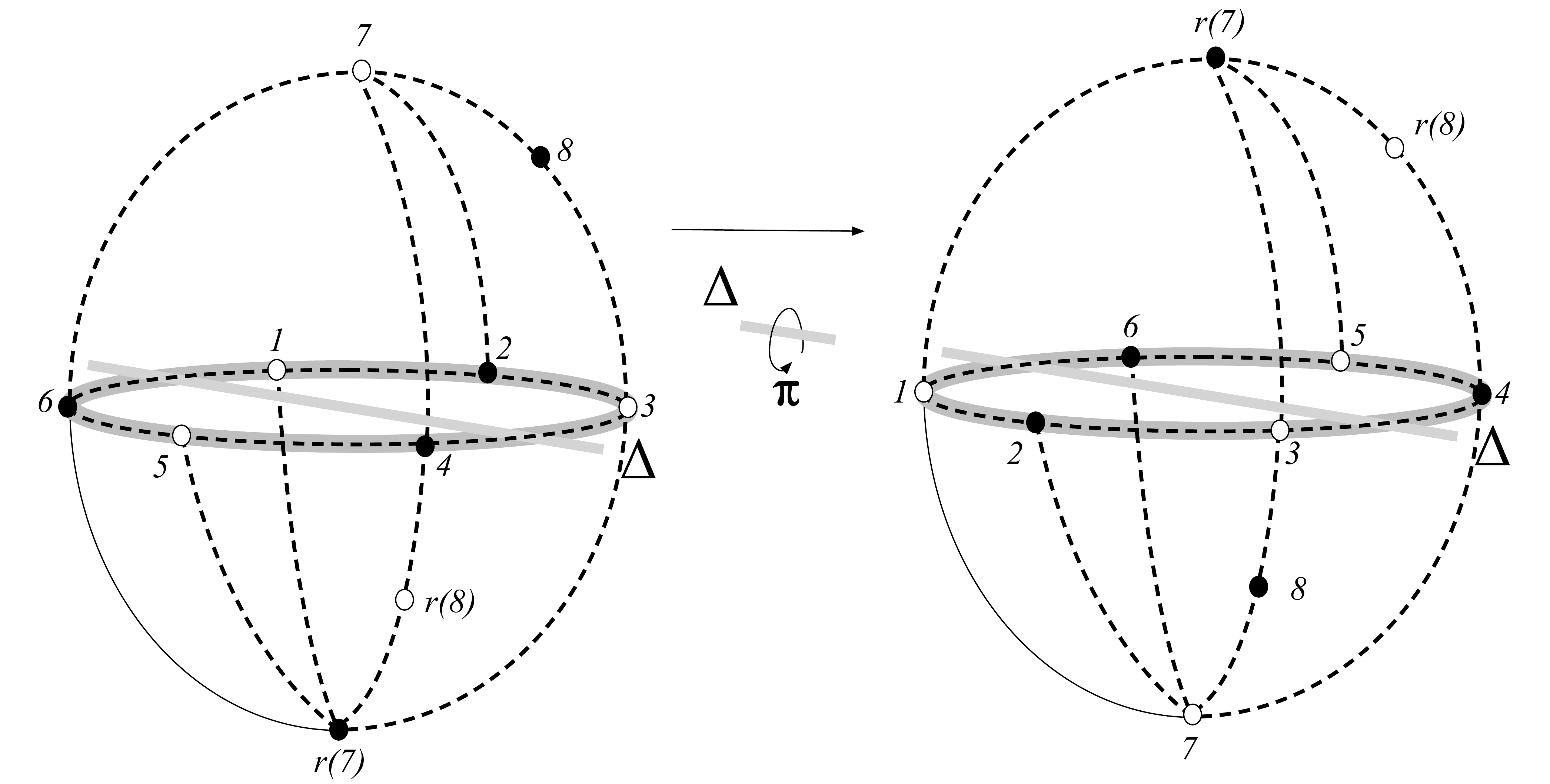}
    \caption{Rotation around $\Delta$ from $(I(G),C_V,S_F)$ to $(I(G),C_V,S_F)$.}
    \label{fig30} 
\end{figure}

(b) Let $\Gamma$ be a sign-reversing color-preserving reflexive curve of $(I(G),C_V,S_F)$. We construct an embedding of $(I(G),C_V,S_F)$ in $\stw$ from which we can pass to  $(I(G),C_V,-S_F)$ by performing a rotation of $\pi$ degrees. Since $(I(G), C_V, S_F)$ determines $L(G,S_E)$ then, by Remark \ref{rem:I(G)}, $(I(G), C_V, -S_F)$ determines $L(G,S_E)^*$. therefore, this rotation of $\pi$ degrees would give the desired orientation-preserving homeomorphism from $L(G,S_E)$ to $L(G,S_E)^*$. 
same  We construct such a map in the way as in the case (a). 
% In this case, the reflexion $\Delta$ must be color-preserving. 
This is illustrated in Figure \ref{fig33b} (applied to the incidence graph given in Figure \ref{fig33a}). 

\begin{figure}[H]
    \centering
    \includegraphics[width=.8\textwidth]{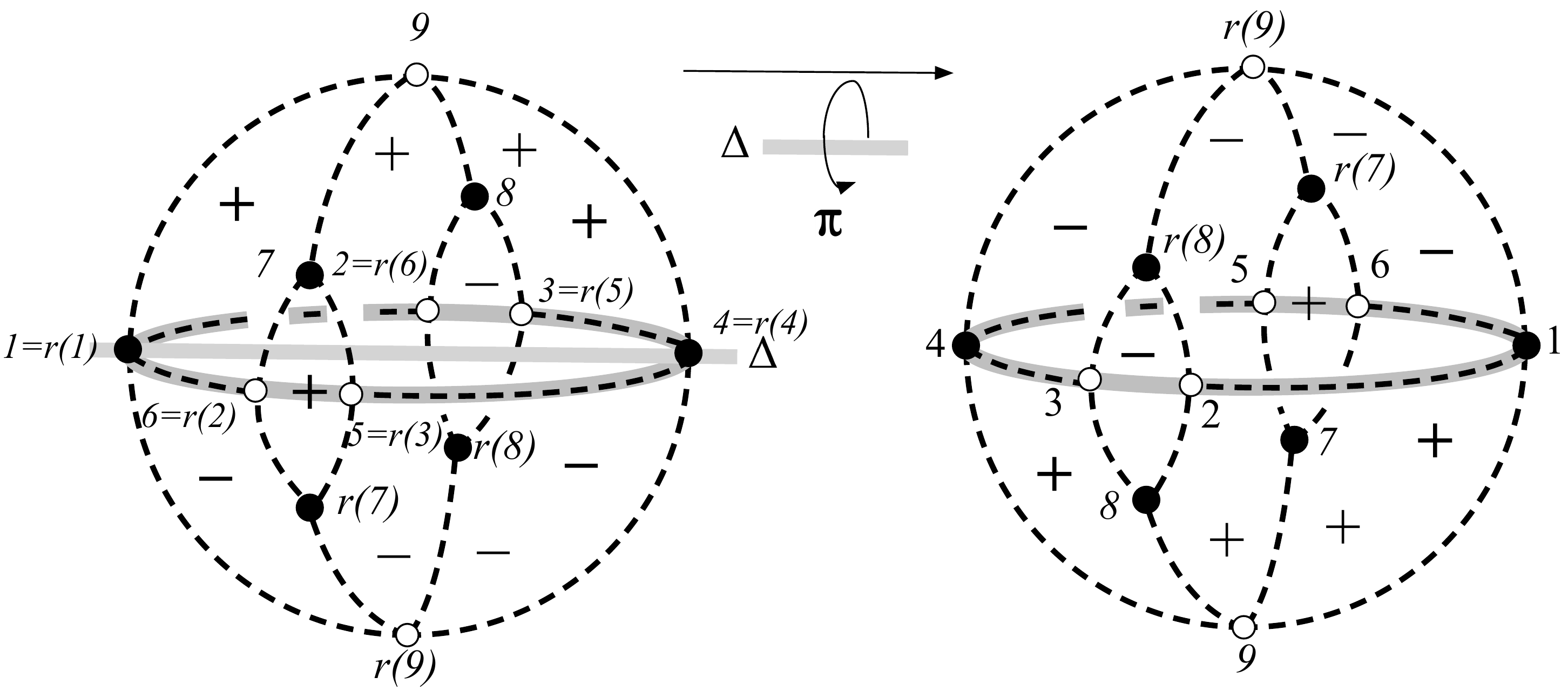}
    \caption{Rotation $\Delta$ from $(I(G),C_V,S_F)$ to $(I(G),C_V,-S_F)$.}
    \label{fig33b}
\end{figure}

We observe that a rotation of $\pi$ degrees around $\Delta$ of the drawing of $(I(G),C_V,S_F)$ exchanges faces $f\in G_{int(\Gamma)}$ and $\sigma(f)\in G_{ext(\Gamma)}$ having different signs (since 
$\Gamma$ is sign-reversing) and  keeping the coloring (since $\Gamma$ is color-preserving). 
We thus have that we can pass from $(I(G),C_V,S_F)$ to $(I(G),C_V,-S_F)$ by applying the rotation of $\pi$ degrees.
\end{proof}

\begin{corollary}\label{cor;w2} Let $n\ge 2$ be an even integer and let $(W_n,S_E^+)$ be an edge-signed map. Then, $L(W_n,S_E^+)$ is amphichiral. 
\end{corollary}

\begin{proof} The result follows by Theorem \ref{theo:amphi} (a) by noticing that $(I(W_n),S_F^+)$ with $n\ge 2$ even always admits a sign-preserving color-reversing reflexive cycle, see Figure \ref{fig31-a}.

\begin{figure}[H]
    \centering
    \includegraphics[width=.3\textwidth]{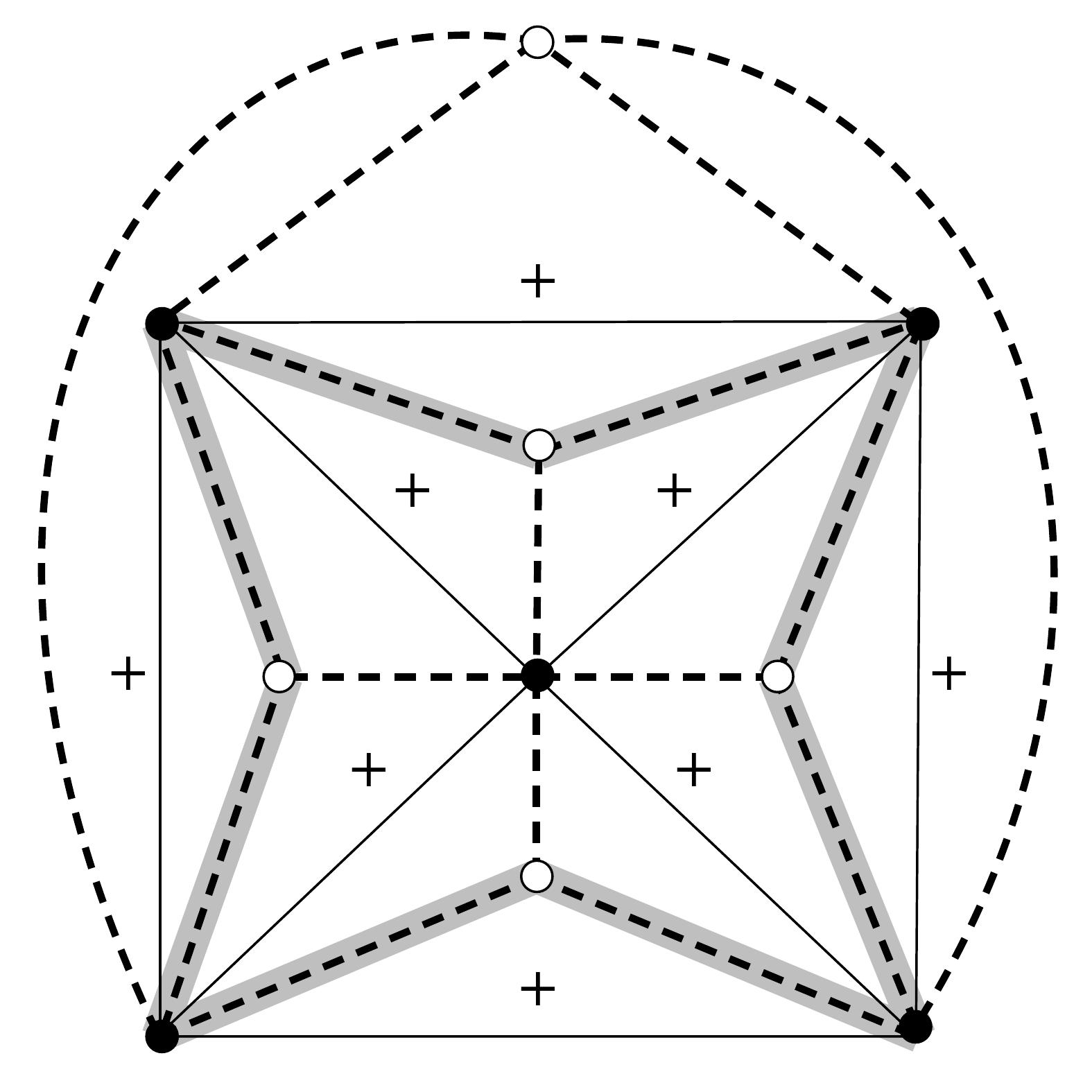}
    \caption{$(W_4,S_E^+)$ (light edges) and $(I(W_4),S_F^+)$ (dashed edges) together with a sign-preserving color-reversing reflexive $\Gamma$-cycle (gray).}
    \label{fig31-a}
\end{figure}
\end{proof}

\begin{figure}
%[H]
    \centering
    \includegraphics[width=1.1\textwidth]{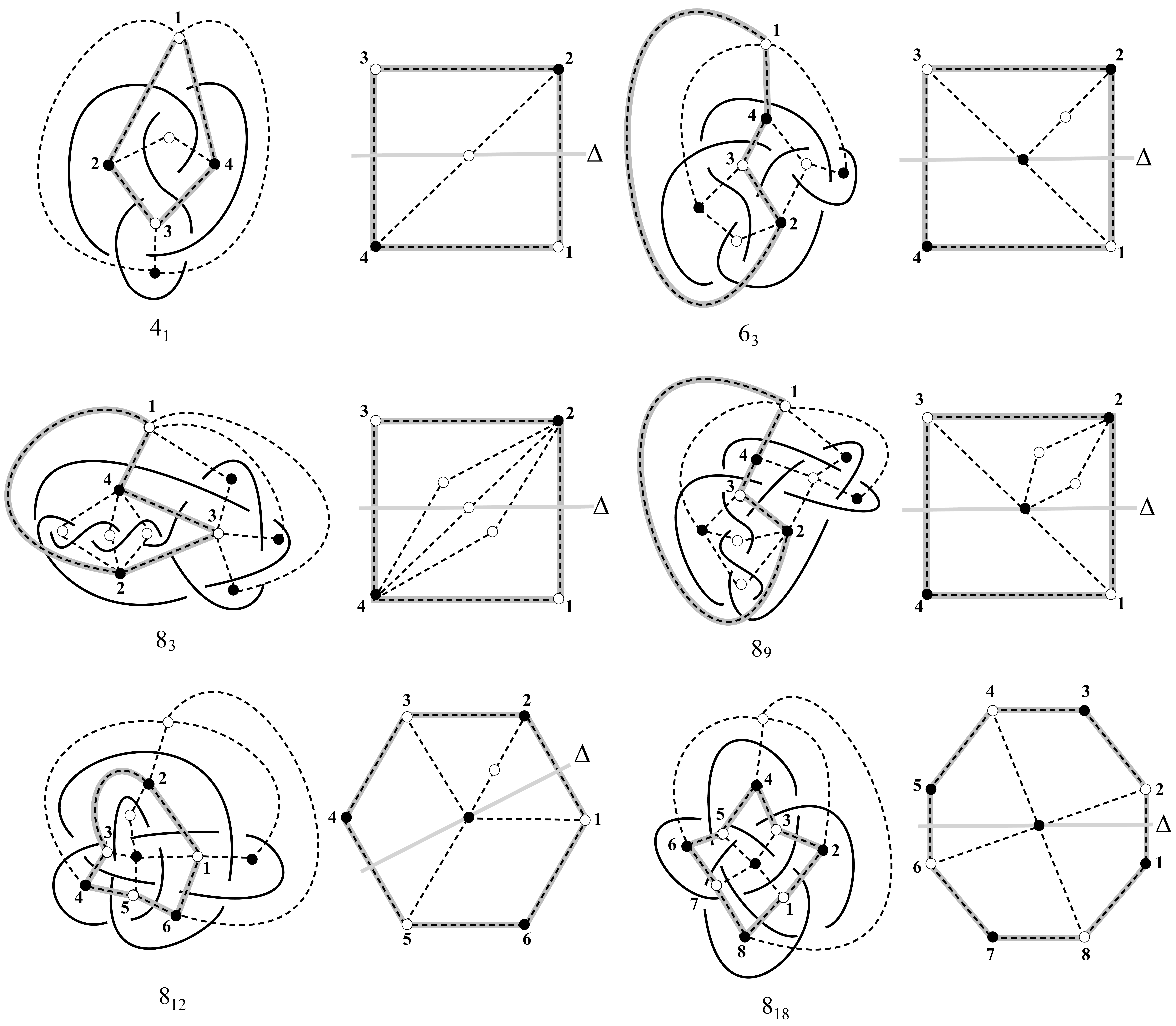}
    \caption{Knots (bold) with the incidence graph (dashed edges) admitting a color-preserving reflexive $\Gamma$-cycle (gray) together with a drawing of $I(G)_{int(\Gamma)}$ with a reflection line.}
    \label{fig31}
\end{figure}

This corollary implies, in particular, that the Eight-figure knot is amphichiral, see first diagram in Figure \ref{fig31}.

%\begin{figure}[H]
   % \centering
%    \includegraphics[width=.5\textwidth]{images/Curve-Link}
 %   \caption{Face-signed map $(I(G),S_F)$ (see Figure \ref{fig33aaa}) together with its associated link $D(G,S_E)$ (bold).} 
   % \label{fig33c}
%\end{figure}

\subsection{Constructing many amphichiral links} 
\begin{figure}[H]
    \centering
    \includegraphics[width=.5\textwidth]{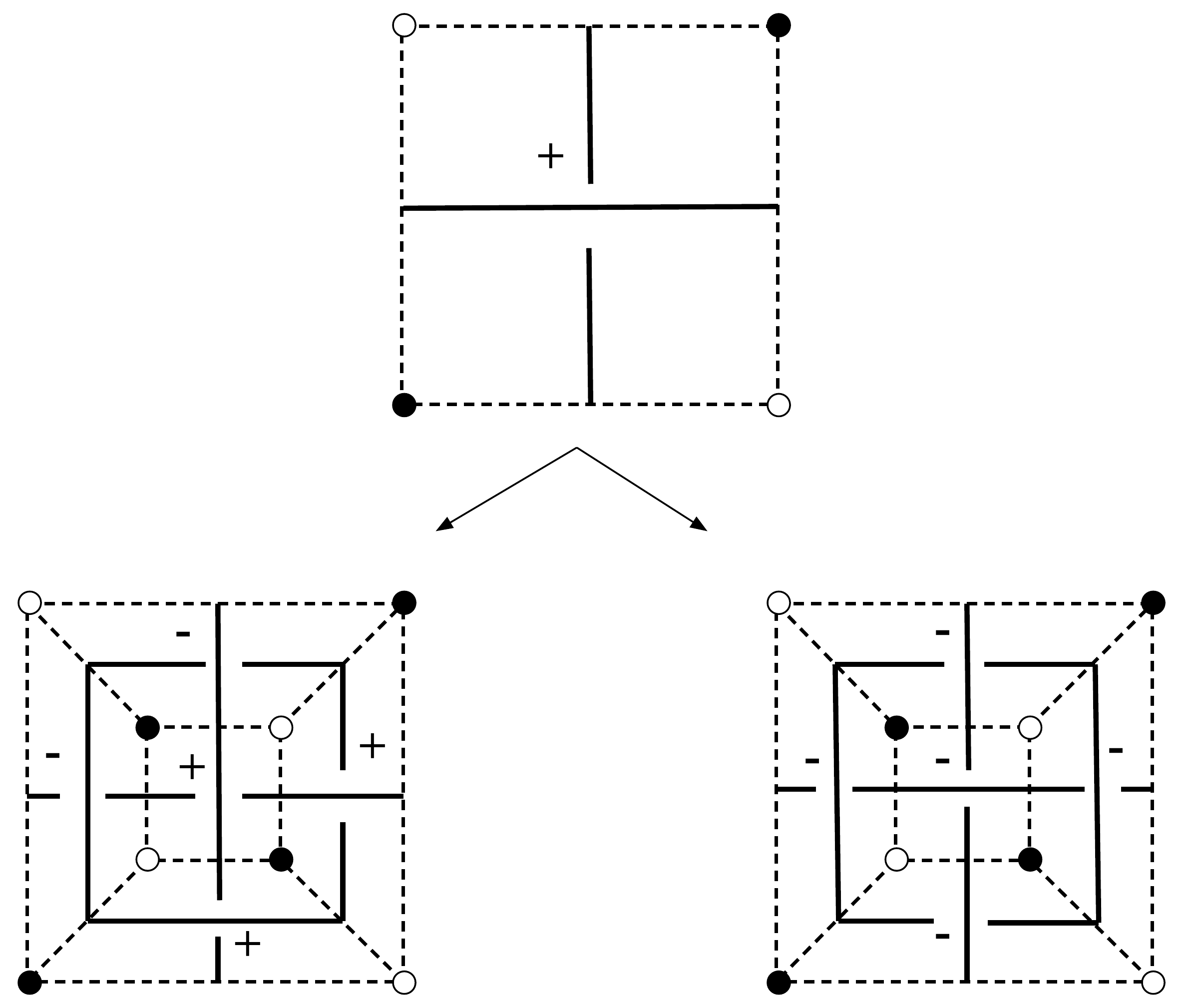}
    \caption{Decomposition of a square into 5 squares. Two different signatures are illustrated inducing a different local modification of the link diagram (in bold lines).}
    \label{fig32}
\end{figure}

Let $(I(G),C_V,S_F)$ be a bicolored face-signed map  where $I(G)$ is the incident graph of $G$. If $I(G)$ admits a  sign-preserving color-reversing reflexive cycle $\Gamma_\sigma$ then we are able to construct another bicolored face-signed map $(I'(G),C_V,S'_F)$ also admitting a sign-preserving color-reversing reflexive cycle $\Gamma_{\sigma'}$ where the isomorphism $\sigma'$ would be a natural extension of $\sigma$. We proceed as follows. 
Let $f$ be a face in the interior of $\Gamma_\sigma$. We partition $f$ into 5 square faces by inserting a new square in the interior and joining corresponding corners (with the appropriate vertex coloring) and assign a sign to each of the new faces, see Figure \ref{fig32}.

Symmetrically, we do the same to decomposition to $\varphi(f)$ where the new squares are properly signed, that is, with the same sign as the corresponding new squares in $f$. Let $\Gamma_{\varphi'}$ be the same cycle as $\Gamma_{\sigma}$ where $\sigma(v)'=\sigma(v)$ for any $v\in V(I(G))$ and $\sigma'(v)=v'$ for each new vertex $v\in f$ and $v'\in \sigma(f)$. We clearly have that $\Gamma_{\sigma'}$ is a sign-preserving color-reversing reflexive cycle in $(I'(G),C_V,S'_F)$. Therefore, by Theorem \ref{theo:amphi},  $(I'(G),C_V,S'_F)$ determines another amphichiral link $L(I'(G),S'_F)$.
 \smallskip
 
\begin{figure}[H]
    \centering
    \includegraphics[width=.34\textwidth]{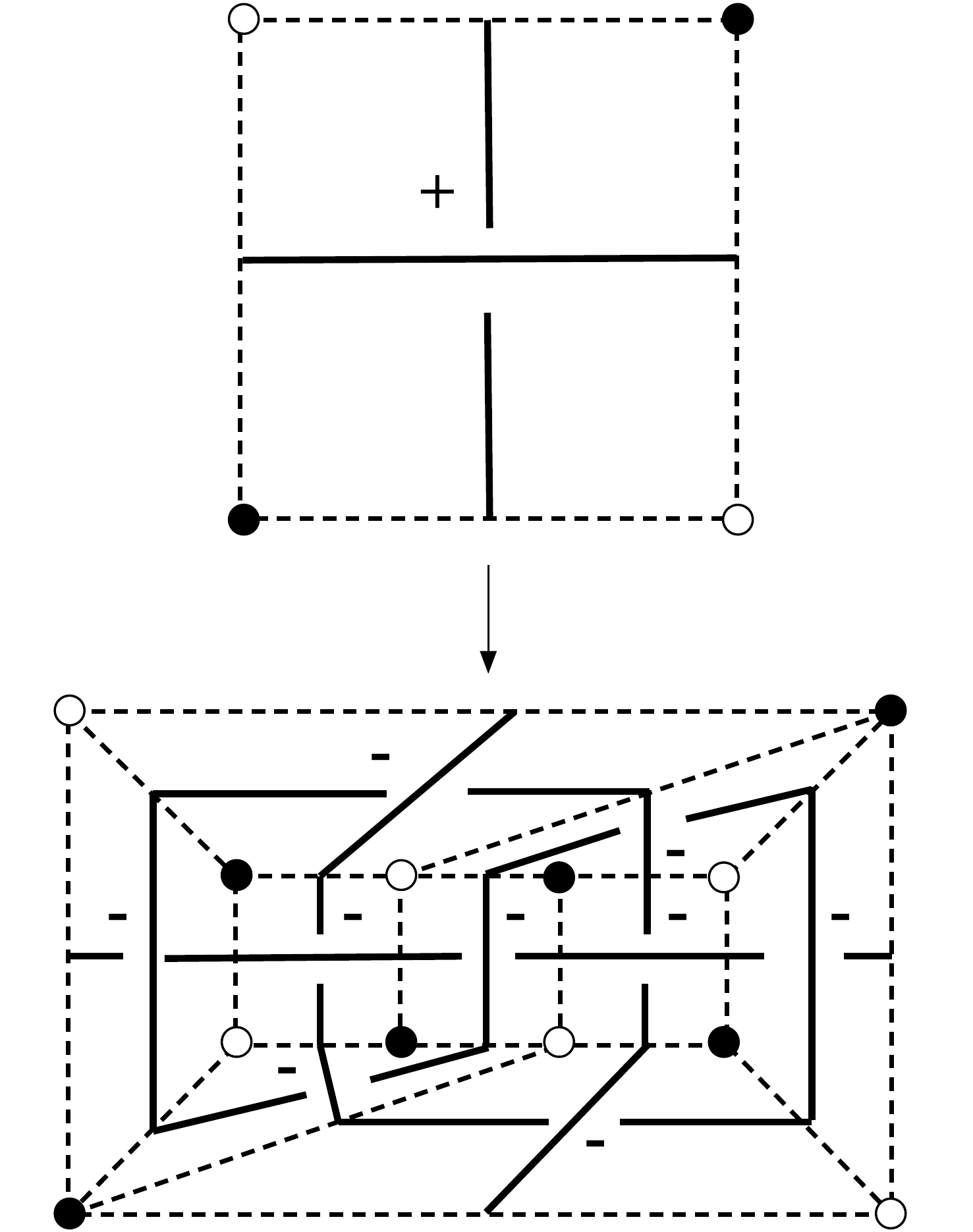}
    \caption{A decomposition of a square into 7 squares together with the link diagram changes (bold lines).}
    \label{fig33}
\end{figure}

 The above construction can be mimic with a different decomposition of a face. For instance, a second decomposition is illustrated in Figure \ref{fig33}.  

\subsection{Amphichiral number} Let $D(L)$ be a diagram of link $L$. 
We define the {\em amphichiral number} of $L$, denoted by $\amph(L)$, as 
\smallskip

$\amph(L):=$ the minimal number of crossing switches in $D(L)$ to become $L$ amphichiral. 
\smallskip

We clearly have that $\amph(L)=0$ if $L$ is already amphichiral and $\amph(K)\le U(K)$ where $U(K)$ is the {\em unknotting number} of $K$, that is, the minimum number of crossing switch to untie $K$.

\begin{proposition}\label{pro:amphinumber} Let $(G,S_E)$ be an edge-signed map with
$G$ antipodally self-dual map. Then,
$$\amph(D(G,S_E))\le \frac{n}{2}$$
where $n$ is the number of crossings of $D(G,S_E)$.
\end{proposition}

\begin{proof} In \cite[Theorem 1]{MRAR1} was proved that if $G$ is an antipodally self-dual map then $I(G)$ admits
% a symmetric cycle of length $2m$ with $m\ge1$ odd. In other words, there is 
an isomorphism $\sigma(I(G))$ such that $\sigma(\Gamma)=\Gamma$ and $\sigma(I(G)_{int(\Gamma)})=I(G)_{ext(\Gamma)}$. Moreover, by Remark \ref{rem;color}, $\Gamma_\sigma$ is color-reversing.  
\smallskip

We may now force $\Gamma_\sigma$ to become sign-preversing by switching the signs of the faces in $I(G)_{ext(\Gamma)}$ properly, that is, such that the face $\sigma(f) \in I(G)_{ext(\Gamma)}$ has the same sign as face $f \in I(G)_{int(\Gamma)}$. We might need to make at most $\frac{n}{2}$ such switches in the worse case (notice that $n$ is even since $G$ is antipodally self-dual). 
Notice that such switchings correspond to crossing switchings of the diagram.The result follows by Theorem \ref{theo:amphi} (a).
\end{proof}

\subsection{Invertible knots} 
\begin{figure}[H]
    \centering
    \includegraphics[width=.7\textwidth]{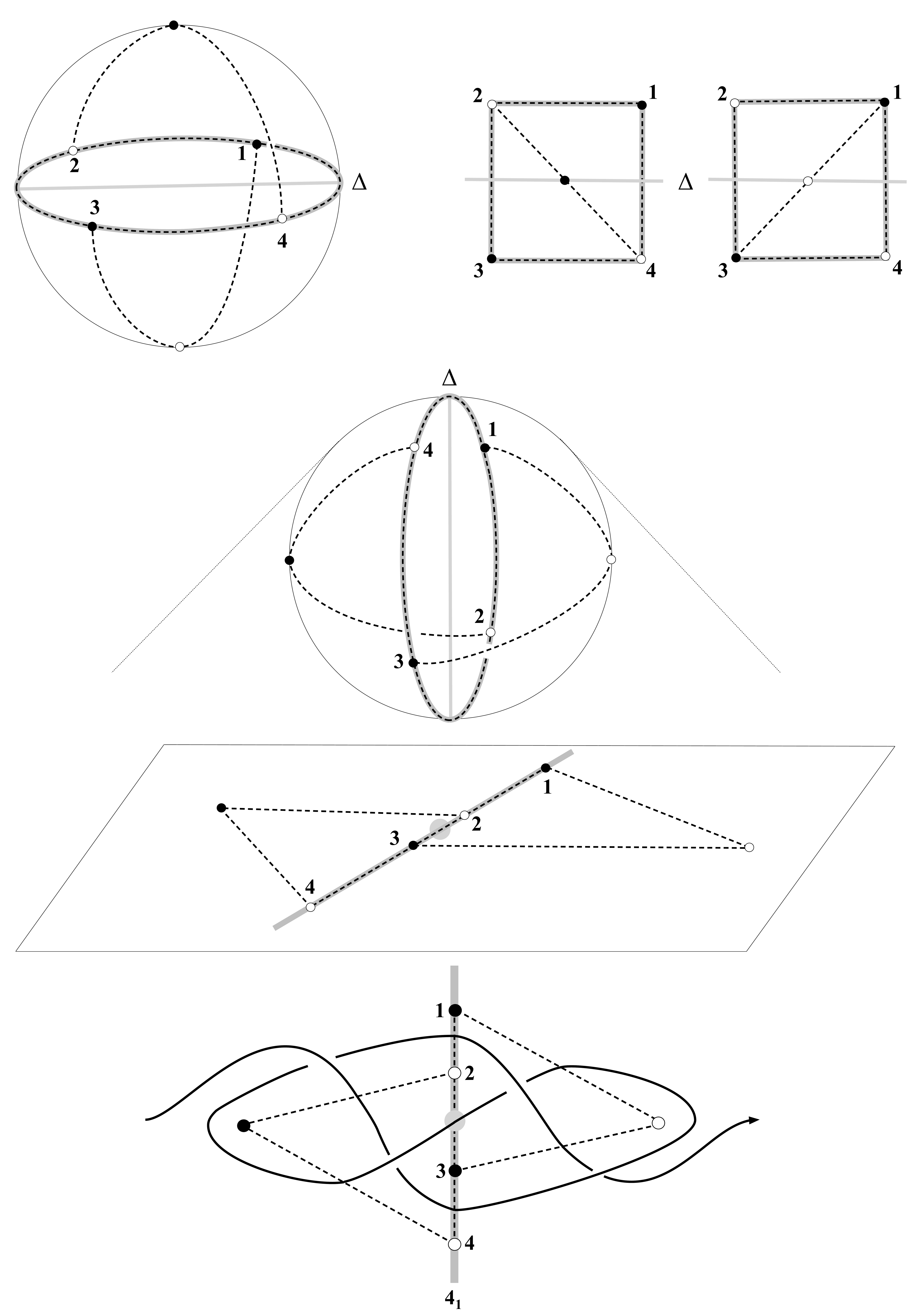}
    \caption{(Top) drawing of $I(G)$ together with the projections of both the Northern and the Southern hemispheres. (Middle) Rotation of drawing with the reflection line $\Delta$ going through the North pole together with the stereographic projection. (Bottom) Diagram of the induced knot $4_1$.}
    \label{fig33a1}
\end{figure}

It turns out that closely related drawings to those considered above for reflexive curves are useful for some other issues. For instance, they are useful to investigate whether a link is rigidly amphichiral in $\sth$ but not necessarily rigidly amphichiral in $\rth$ (the other direction is always true). Indeed, we may construct appropriate drawings  to show the latter as follows. Take a drawing of $G_{int(\Gamma)}$ together with the reflection line $\Delta$ in $\mathbb{B}^1$ with $\Gamma_\sigma$ forming the vertices of a regular polygon lying in the equator of $\stw$.  We then project this to both the Northern and the Southern hemispheres of $\stw$ and make a rotation of $\stw$ such that line $\Delta$ going through the North pole. We finally take a stereographic projection in which $\Gamma$ becomes an infinity line (the ends of the string are closed up at infinity in $\sth$), see Figure \ref{fig33a1}.

\begin{figure}
%[H]
    \centering
    \includegraphics[width=1\textwidth]{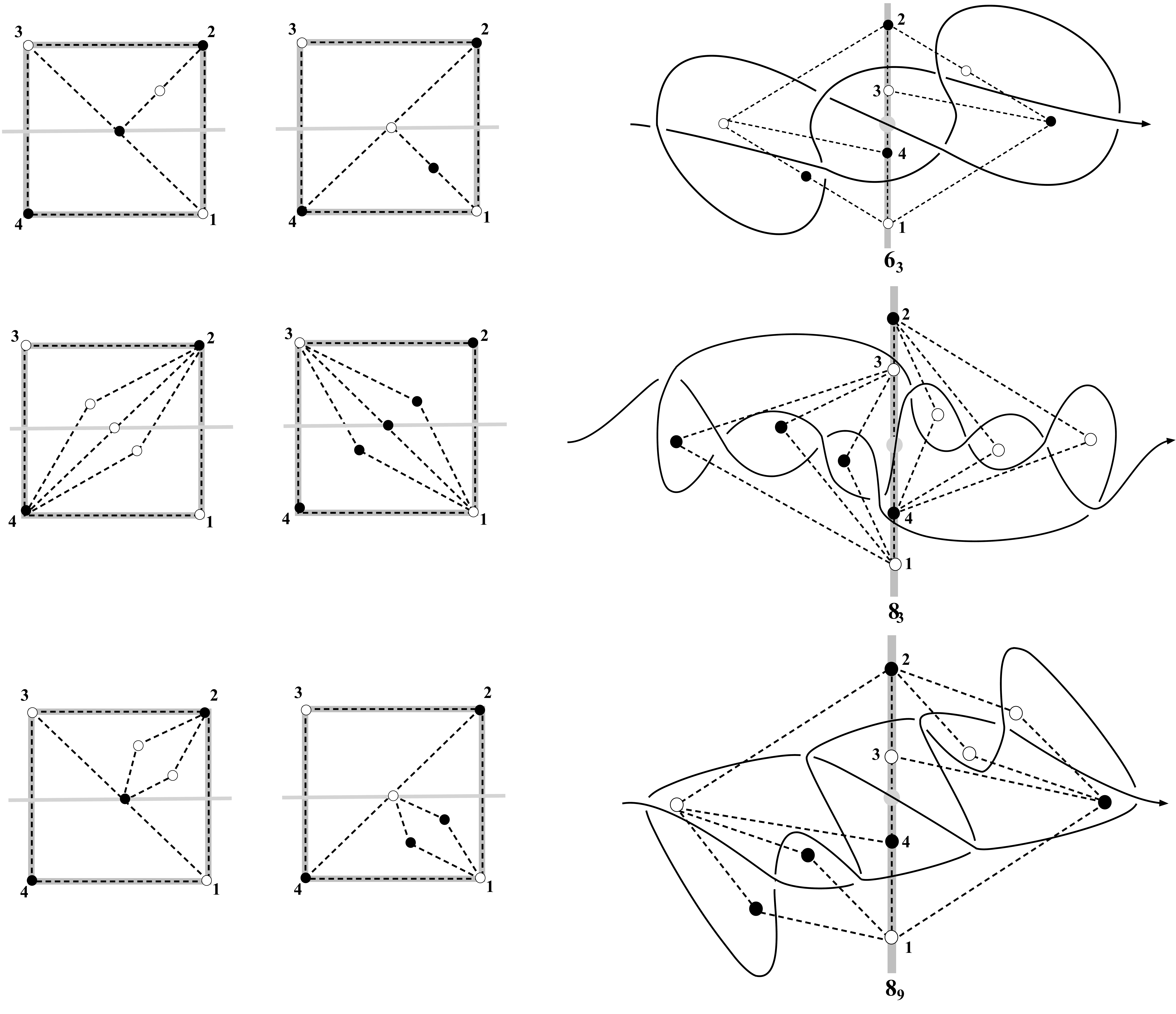}
    \caption{Symmetry representation for knots $6_3, 8_3$ and $8_9$ in $\sth$ that can be rotated by $\pi$ degrees to obtain their mirror images.}
    \label{fig33b1}
\end{figure}

\begin{figure}
%[H]
    \centering
    \includegraphics[width=1\textwidth]{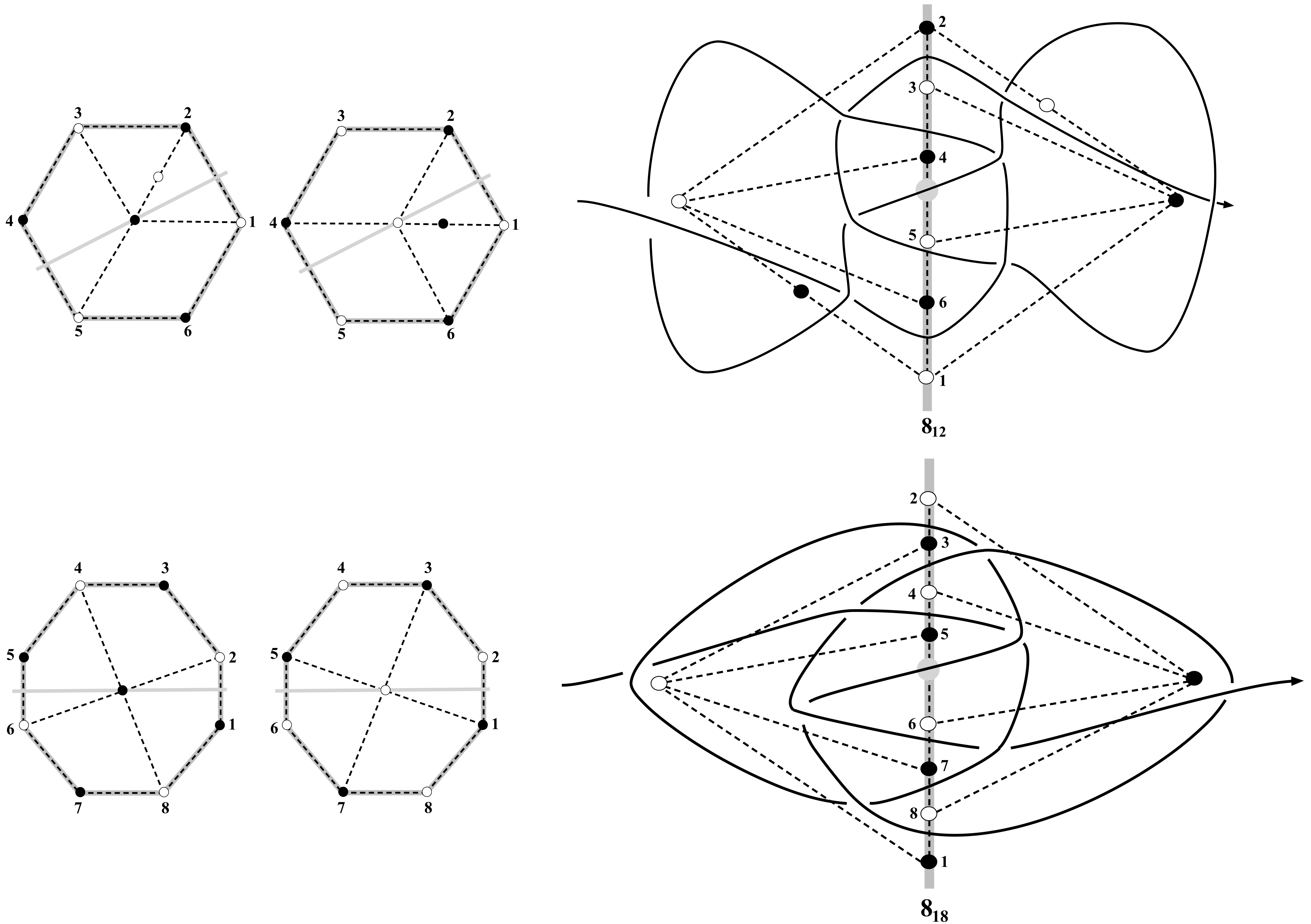}
    \caption{Symmetry representation for knots $8_{12}$ and $8_{18}$ in $\sth$ that can be rotated by $\pi$ degrees to obtain their mirror images.}
    \label{fig33b2}
\end{figure}

Observe that these embeddings are invariant under the orientation-preserving homeomorphism of $\sth$ given by a rotation by $\pi$ degrees about an axis going through the center point of the knot and the point at infinity, composed with a reflection through the plane of the paper. We notice that if we want to restrict a homeomorphism of $\sth$ to a homeomorphism of $\rth$, we must select a point that is fixed by the homeomorphism to be our point at infinity in $\sth$. Both of the fixed points of the homeomorphism in our construction are on the knot, so there is no way to restrict this homeomorphism of $\sth$ to a homeomorphism of the knot in $\rth$, see \cite[Chapter 4]{Flap} for a nice discussion on this topic.
\smallskip

We also observe that a rotation of $\pi$ degree about an axe going through the center point (of the plane) of any of these oriented diagrams give the diagrams of their mirrors with a reversed orientation, that is, the corresponding knots are negative amphichiral.

%Let us suppose that $\Gamma_G$ is placed in an imaginary circle $C$ with $\Gamma_G^{int}$ (resp.$\Gamma_G^{ext}$) inside (resp. outside) $C$.  We say that $\Gamma_G$ is {\em rotatable} if there is a nontrivial rotation of $C$ (dragging $\Gamma_G, \ \Gamma_G^{int}$ and $\Gamma_G^{ext}$ together with the corresponding signs) leaving $I(G)$ invariant and in which the sign of each shared faces is swapped.  

%\begin{theorem}\label{theo:amphi-b} Let $(G,S_E)$ be an edge-signed map. Suppose that $(I(G),S_F,C_V)$ admits an even adequate sign-reversing curve $\Gamma_G$. If $\Gamma_G$ is rotatable then $D(G,S_E)$ is amphichiral. 

%\jorge{Moreover, $D(G,S_E)$ and its mirror can be embedded such that one can be obtained from the other by one rotation.}\end{theorem}

%see Figure \ref{fig30c}. 

%\begin{figure}[H]
  %  \centering
   % \includegraphics[width=1\textwidth]{Steps123inverse}
 %   \caption{Isometrie from $Map(I(G))$ with signature $S_F$ to $Map(I(G))$ with signature $-S_F$.}
 %   \label{fig30c}
%\end{figure}

%%%%%%%%%%%%%%%%%%%%%%%%%%%%%%%%%%%%%%%%%%%%%%%%
\section{Self-dual pairing}\label{sec;isomap}

Let us quickly recall some notions on the classification of the self-dual maps.
\smallskip

Let $Aut(G)$ be the group formed by the set of all {\em automorphism} of $G$ (i.e., the set  of isomorphisms of $G$ into itself). Let $Dual(G)$ be the set of all {\em duality isomorphisms} of $G$ into $G^*$. We notice that $Dual(G)$ is not a group since the composition of any two of them is an automorphism. 
\smallskip

Let us suppose that $G=(V,E,F)$ is a self-dual map so that there is a bijection 
 $\phi:(V,E,F)\rightarrow (F^*,E^*,V^*)$. Following $\phi$ with the correspondence $*$ gives a permutation on $V\cup E\cup F$ which preserve incidences but reverses dimension of the elements. The collection of all such permutations or {\em self-dualities} generate a group $Cor(G)=Aut(G)\cup Dual(G)$ in which the automorphisms $Aut(G)$ are contained as a subgroup of index 2.   
 \smallskip
 
It is known \cite[Lemma 1]{SS} that for a given map $G$ there is an homeomorphism $\rho$ of $\mathbb{S}^2$ to itself such that for every $\sigma\in Aut(G)$ we have that $\rho\sigma\in Isom(\mathbb{S}^2)$ where $Isom(\stw)$ is the group of isometries of the 2-sphere. In other words, any planar graph $G$ can be drawn on the 2-sphere such that any automorphism of $G$ act as an isometry of the sphere.  This was extended in \cite{SS} by showing that given any self-dual graph $G$ there are maps $G$ and $G^*$ so that $Cor(G)$ is realized as a group of spherical isometries. 

From now on, we will denote by $\widehat{G}=\rho(G)$ and $\widehat{\sigma}=\rho \sigma$ for a certain homeomorphism $\rho$ satisfying the above property.
\smallskip

The couple $Cor(G) \rhd Aut(G)$ is called the {\em self-dual pairing} of the map $G$. In \cite{SS} were enumerated and classified all self-dual maps. In the notation of \cite{Cox} the possible 24 pairings are :
 \smallskip
 
 - among the infinite classes $[2,q]\rhd [q], [2,q]^+\rhd [q]^+, [2^+,2q]\rhd [2q], [2,q^+]\rhd [q]^+$ and $[2^+,2q^+]\rhd [2q]^+$ or 
 \smallskip
 
 - among the special pairings $[2]\rhd [1],  [2]\rhd [2]^+, [4]\rhd [2], [2]^+\rhd [1]^+, [4]^+\rhd [2]^+, [2,2]\rhd [2,2,]^+, [2,4]\rhd [2^+,4], [2,2]\rhd [2,2^+], [2,4]\rhd [2,2], [2,4]^+\rhd [2,2]^+, [2^+,4]\rhd [2,2]^+, [2^+,4]\rhd [2^+,4^+], [2,4^+]\rhd [2^+,4^+], [2,2^+]\rhd [2^+,2^+], [2,4^+]\rhd [2,2^+], [2,2^+]\rhd [1], [3,4]\rhd [3,3], [3,4]^+\rhd [3,3]^+$ and $[3^+,4]\rhd [3,3]^+$.
\smallskip

We may translate the notion of symmetric $\Gamma$ curve, of the precede section, into a dual-pairing terms. To this end, we first notice that $Isom(\mathbb{S}^2)$ consists of three types of elements: $(i)$ rotations, $(ii)$ reflexions and $(iii)$ rotary reflexion (a rotation and a reflexion). We also observe that there is always a plane $H$ stable under these isometries: for a rotation $H$ is the plane perpendicular to an axe, for a reflexion $H$ is either the plane used for the reflexion or any plane perpendicular to it and for a rotary reflexion $H$ is either the plane used for the reflexion or any plane perpendicular to it if $\theta=0$ where $\theta$ is the rotation angle or any plane containing the origin if $\theta=\pi$.
\smallskip

Let $G$ be self-dual map and let $\widetilde{med}(G)$ be a redraw of $med(G)$ such that $Aut(med(G)) < Isom(\mathbb{S}^2)$. For $\sigma\in Aut(\widetilde{med}(G))$, we define 
$T_{\sigma,H}=H\cap \mathbb{S}^2$ where $H$ is a stable plane.  We have that if $T_{\sigma,H}$ does not go through any vertex of $med(G)$ then $T_{\sigma,H}$ induces a symmetric cycle $\Gamma$ in $I(G)$ and if $T_{\sigma,H}$ goes through vertices of $med(G)$ then $T_{\sigma,H}$ induce a curve $\Gamma$ curve with interval vertices of $I(G)$.

\subsection{Equivalent maps}

We say that two maps $G_1$ and $G_2$ of a same graph are {\em equivalent} if there exists an orientation-preserving homeomorphism $\varphi: \stw\rightarrow \stw$ with $\varphi(G_1)=G_2$. Similarly as done for links, we shall write $[G]$ (called {\em map-type}) the class of maps equivalent to $G$.

\begin{remark} We have that $[G]$ is an equivalent relation on the set of all embeddings of $G$ and it depends entirely on the embedding and not on the abstract graph itself. 
\end{remark}

\begin{figure}[H]
    \centering
    \includegraphics[width=0.73\textwidth]{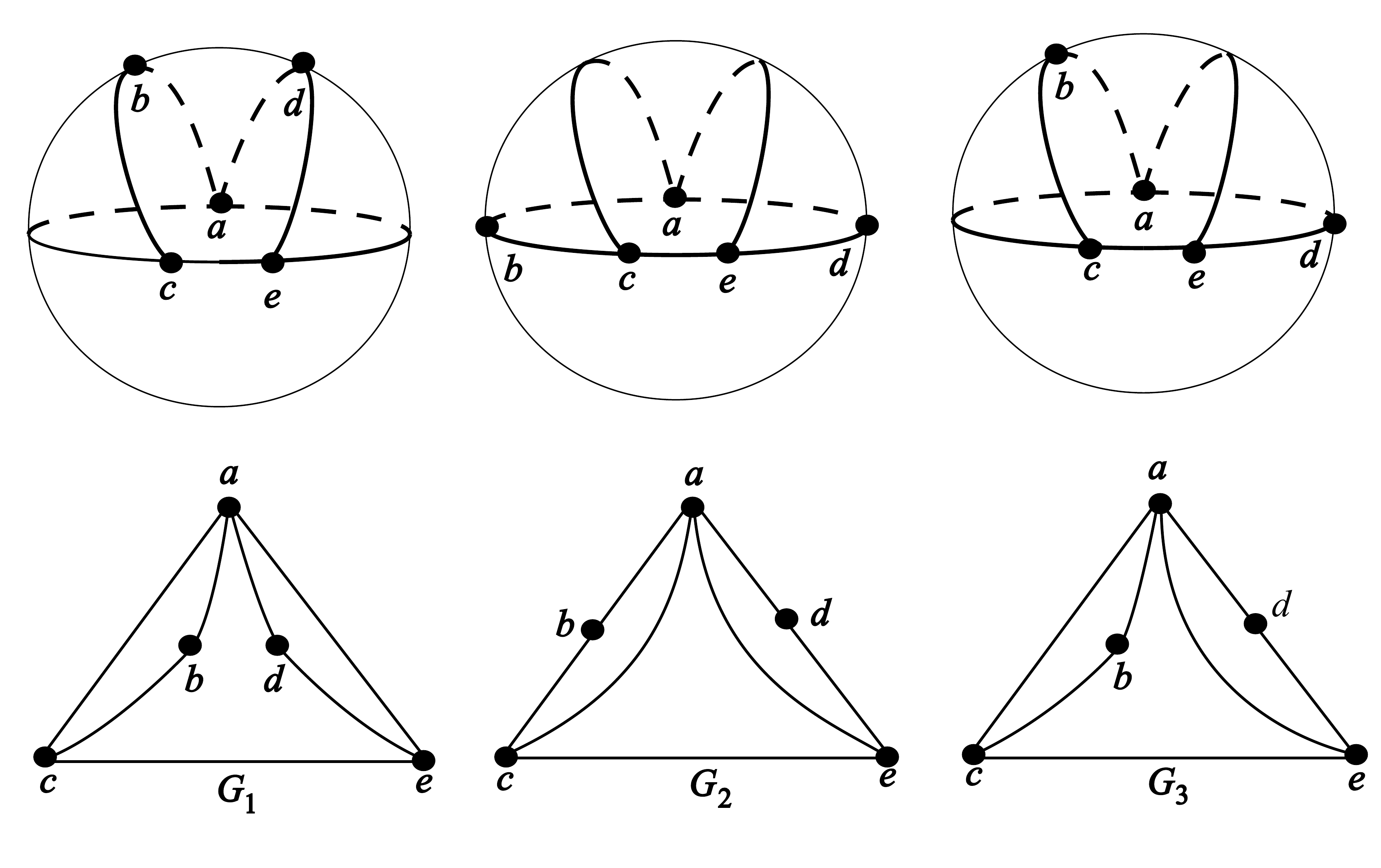}
    \caption{Graphs $G_1, G_2$ and $G_3$ are isomorphic to each other (bijections are given by the same labeled vertices). However, $G_1\in [G_2]$ and $G_1\not\in [G_3]$.}
    \label{equiv}
\end{figure}

\begin{lemma}\label{dual} Let $G_1$ and $G_2$ be two maps of a same graph. Then,
$G_1\in [G_2]$ if and only if $G_1^*\in [G_2^*]$.
\end{lemma}

\begin{proof} Let $F_1$ and $F_2$ be the set of faces of map $G_1$ and map $G_2$ respectively.
For every homeomorphism $h$ of $\stw$ to itself sending $G_1$ to $G_2$ we have that $h|_f\in F_2$ for every face $f\in F_1$. Furthermore, if two faces $f_1$ and $f_2$ of $G_1$ share an edge $e$ in their boundaries then $f_1^*$ and $f_2^*$ share the edge $h|_e\in E(G_2)$. It follows that the image of every geometric dual of $G_1$ is also a geometric dual of $G_2$.
\end{proof}

We have the following easy consequence of the above lemma.

\begin{corollary}\label{corodual} Let $G_1\in [G_2]$. Then, $med(G_1)\in [med(G_2)]$ and 
$I(G_1)\in [I(G_2)]$.
\end{corollary}

%By Lemma \ref{dual},  we have that for each equivalent class $\overline{G}$ there is a unique class $\overline{G^*}$ and therefore $\overline{(G)}^*:=\overline{(G^*)}$ is well defined.  Moreover,
%by Corollary \ref{corodual} we can have 

%$$\hbox{$\overline{(G)}^\square:=\overline{(G^\square)},\ med{\overline{(G)}}:=\overline{(med(G))}$ and $I(\overline{G)}:=\overline{(I(G))}$.}$$ 
%Without further precisions, in the following we will confuse a plane graph $G$ with its equivalence class $\overline{G}$. 
  
 \begin{lemma}\label{2to3} Let $(G_1,S_E)$ and $(G_2,S_E)$ be two maps such that $[G_1]=[G_2]$.
Then, $[L(G_1,S_E)]=[L(G_2,S_E)]$.
 \end{lemma}
 
 \begin{proof}  Let $D(G_1,S_E)$ be the link diagram induced by $(med(G_1), S_V, C_F)$ and let $\theta$ be the map from $D(G_1,S_E)$ to $\embed(G_1,S_E)$.  By Corollary \ref{corodual}, there is an orientation-preserving homeomorphism $h: \stw\rightarrow\stw$ such that $h(med(G_1))=med(G_2)$. We clearly have that $h$ preserves edge-signs and face-colorings.
\smallskip
  
The composition $\theta\circ h\circ\theta^{-1}$ induce an homeomorphism from $\rth$ to itself such that $$\theta\circ h\circ\theta^{-1}(\embed(G_1,S_E))=L(G_2,S_E).$$
\end{proof}

 \begin{theorem}\label{thm:self-dual-amphi} Let $(G, S_E)$ be an edge-signed map such that $G\in [G^*]$. Then, $L(G,S_E)$ is amphichiral for every signature $S_E$.
 \end{theorem}
 
 \begin{proof} 
$$\begin{array}{llll}
L(G,S_E)^*&\in &[L(G,-S_E)]& \mbox{(by definition of mirror)}\\
            &=&[L(G^*,S_E)]& \mbox{(by the sign rules)}\\
            &=&[L(G,S_E)]&\mbox{(by Lemma \ref{2to3} since $G\in [G^*]$)}.              
 \end{array}$$
\end{proof}

\subsection{Self-dual maps}\label{subsec;permutation2}
We shall give an amphichirality result by using the classification of the self-dual maps. 

 \begin{lemma}\label{lem;dualperserving} Let $G$ be a self-dual map. If either
 
(a) there exists $\tau \in Dual(G)$ such that the isometry $\widehat \tau$ is oriented-preserving or
 \smallskip
 
(b) there exists $\sigma\in Aut(G)$ such that the isometry $\widehat \sigma$ is not oriented-preserving 
 \smallskip
 
 then $L(G,S)$ is amphichiral for every signature $S$. 
 \end{lemma}
 
 \begin{proof} (a) The fact that the isometry $\widehat \tau$, arising from the duality isomorphism $\tau(G)=G^*$, is oriented-preserving gives an isotopy between $\widehat G$ and $\widehat G^*$. The result follows by Theorem \ref{thm:self-dual-amphi}. 
\smallskip

(b) Let $\eta\in Dual(G)$ be a duality isomorphism. If the isometry $\widehat \eta$ is oriented preserving then the result follows by case a). Otherwise, we consider 
 $$\widehat \sigma\circ \widehat \eta(\widehat G)=\widehat \sigma(\widehat G^*)\stackrel{(G \text{ self-dual})}{=}\widehat \sigma(\widehat G)\stackrel{(\sigma\in Aut(G))}{=}\widehat G\stackrel{(G \text{ self-dual})}{=}\widehat G^*$$
 
which is clearly oriented-preserving isometry (composition of two non oriented-preserving isometries) implying that $\widehat G\in [\widehat G^*]$. The result follows by Theorem \ref{thm:self-dual-amphi}. 
\end{proof}
 
% \begin{lemma}\label{lem;pocket}  ({\em Pocket Lemma}) Let $G$ be a self-dual map. If there exists $\sigma\in Aut(G)$  such that the isometry $\widehat \sigma$ is not oriented-preserving then $D(G,S)$ is amphichiral for every signature $S$. 
 %\end{lemma}
 
% \begin{proof} Let $\rho\in Dual(G)$ be a duality isomorphism. If the isometry $\widehat \rho$ is oriented preserving then the result follows by Lemma \ref{lem;dualperserving}. Otherwise, we take $\widehat \sigma$ (out of the pocket) and consider 
%$$\widehat \sigma\circ \widehat \rho(G)=\widehat \sigma(G^*)\stackrel{(G \text{ self-dual})}{=}\widehat \sigma(G)\stackrel{(\widehat \sigma\in Aut(G))}{=}G\stackrel{(G \text{ self-dual})}{=}G^*$$
%which is clearly oriented-preserving isometry (composition of two non oriented-preserving isometries) implying that $G\simeq G^*$. The result follows by Theorem \ref{thm:self-dual-amphi}. 
 %\end{proof}

\begin{theorem}\label{theo;amphi1} Let $(G,S_E)$ be an edge-signed self-dual map. If the self-dual pairing of the map $G$ is other than $[2,q^+]\rhd [q]^+,[2^+,2q^+]\rhd [2q]^+, [2]\rhd [2]^+, [2,2]\rhd [2,2]^+,
[2^+,4]\rhd [2,2]^+, [3^+,4]\rhd [3,3]^+$ then $L(G,S_E)$ is amphichiral for every signature $S_E$. 
\end{theorem}

\begin{proof}  It can checked that in any self-dual pairing other than the 6 given in the hypothesis we always have that there is either $\rho\in Dual(G)$ with $\widehat \rho$ oriented-preserving or $\sigma\in Aut(G)$ with $\widehat \sigma$ non oriented-preserving. Table \ref{tab1} points out the appropriate isometry in each case. The result then follows by Lemma \ref{lem;dualperserving}.  

\begin{table}[H]
%[htbp]
\begin{center}
\begin{tabular}{|l|c|}
\hline
\text{ Self-dual pairing} & Type of isometry\\
\hline
\hline
$[2,q]\rhd [q]$ &  (2) \\
\hline
$[2,q]^+\rhd [q]^+$ & (1)\\
\hline 
$[2^+,2q]\rhd [2q]$ & (2)\\
\hline
 $[2]\rhd [1]$ & (1)\\
\hline 
 $[4]\rhd [2]$ & (1)\\
\hline
 $[2]^+\rhd [1]^+$ & (1)\\
 \hline
 $[4]^+\rhd [2]^+$ & (1)\\
 \hline
 $[2,4]\rhd [2^+,4]$ & (2)\\
\hline
 $[2,2]\rhd [2,2^+]$ & (2)\\
 \hline
 $[2,4]\rhd [2,2]$ & (2)\\
 \hline
 $[2,4]^+\rhd [2,2]^+$ & (1)\\
\hline
 $[2^+,4]\rhd [2^+,4^+]$ & (2)\\
 \hline
 $[2,4^+]\rhd [2^+,4^+]$ & (2)\\
 \hline
 $[2,2^+]\rhd [2^+,2^+]$ & (2)\\
 \hline
 $[2,4^+]\rhd [2,2^+]$ & (2)\\
 \hline
 $[2,2^+]\rhd [1]$ & (2)\\
 \hline
 $[3,4]\rhd [3,3]$ & (2)\\
 \hline
 $[3,4]^+\rhd [3,3]^+$ & (1)\\
 \hline
\end{tabular}
\end{center}
\caption{Self-dual pairings having an isometry of type either (1) if 
there is $\rho\in Dual(G)$ with $\widehat \rho$ oriented-preserving or (2) if there is $\sigma\in Aut(G)$ with $\widehat \sigma$ not oriented-preserving.}\label{tab1}
\end{table}
\end{proof}

We notice that in each of the 6 special self-dual pairings in Theorem \ref{theo;amphi1}, we have that both all $\rho \in Dual(G)$ arise non oriented-preserving isometries and all $\sigma\in Aut(G)$ arise  
oriented-preserving isometries. 
%\smallskip

%\jorge{
%\subsection{Symmetry group} In \cite{GS}, Gr\"unbaum and Sherphard studied the {\em symmetry groups} of a knot $K$ defined as the groups of isometries of $\mathbb{R}^3$ that map the knot onto itself. They determined the only fixed-point groups of isometries $\mathbb{R}^3$ that can occur as symmetry groups of knots. They mentioned that it would be interesting to determine all possible groups for all links with prescribed number of loops. On this direction, we present the following result.

%\begin{theorem}\label{thm:symgroup} Let $(G,S_+)$ be a signed map with $G$ antipodally self-dual. Then, $D(G,S_+)$ admits  "dar la lista de los grupos que tienen la funcion antipoda" as symmetry groups.
%Moreover, $D(G,S)$ belongs to one of the following self-dual pairing : "dar la lista de los self-dual pairings que contienen la funcion antipoda"
%\end{theorem} 

%\begin{proof} By Corollary \ref{cor:antipodal}, $D(G,S_+)$ is 3-centrally symmetric. We thus have that any isometry of $med(G)$ admitting the antipodal function $\alpha_2$ would leave $D(G,S)$ invariant. We determine the list of groups containing the antipodal mapping. 
%Moreover, by using the above, it can be verified that $a\in Dual(G)$ and $a\not\in Aut(G)$ for each self-dual pairing in the desired list.
%\end{proof}
%} 

%%%%%%%%%%%%%%%%%%%%%%%%%%%%%%%%%%%%%%%%%%%%%%%
%%%%%%%%%%%%%%%%%%%%%%%%%%%%%%%%%%%%%%%%%%%%%%%

\section{Concluding Remarks}

\begin{figure}[H]
\centering
\includegraphics[width=.7\linewidth]{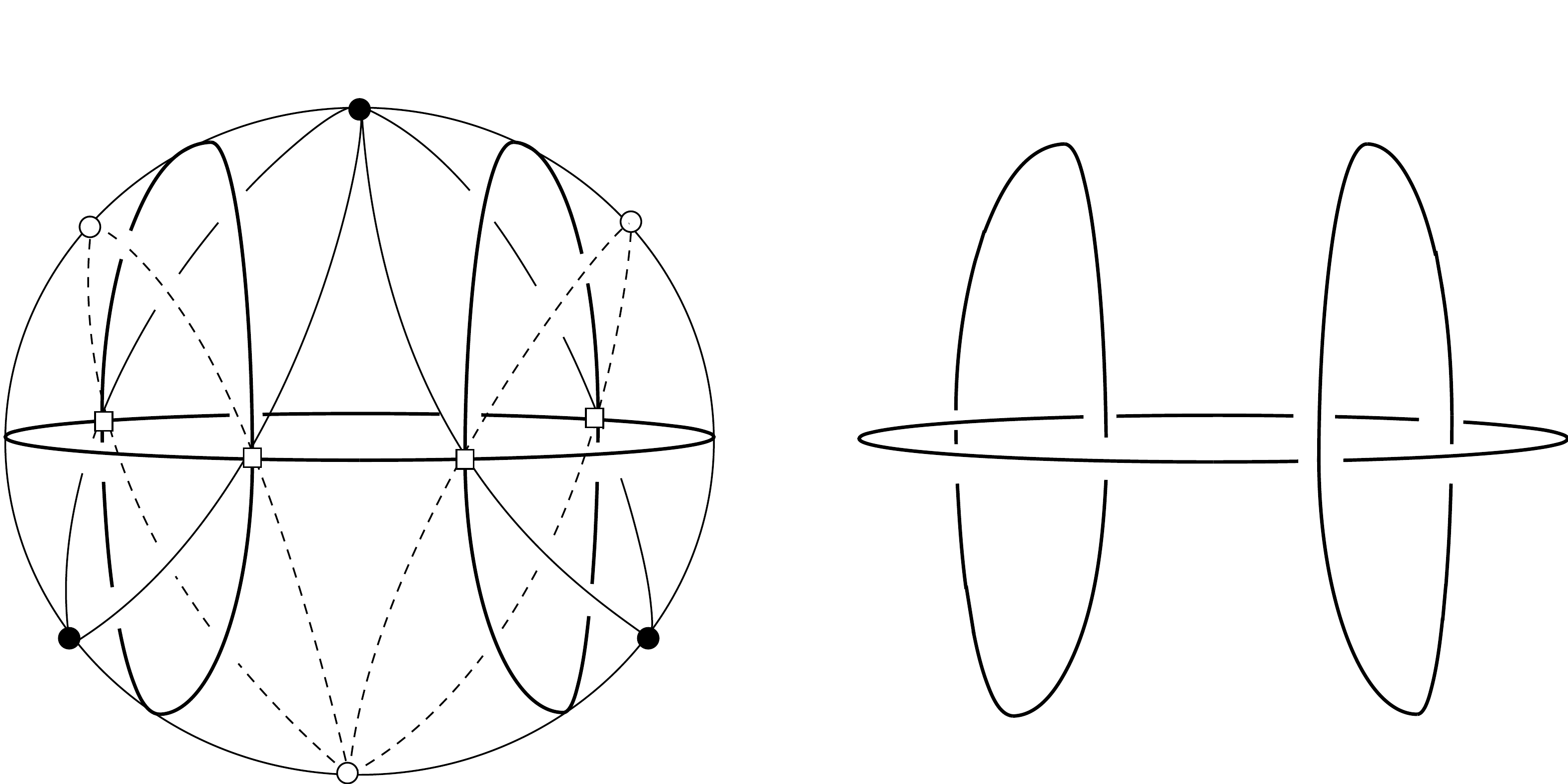}
\caption{(Right) An antipodally self-dual map $G$ (black vertices and straight edges) and $G^*$ (white vertices and dotted edges) and a map of $med(G)$ 3-antipodally symmetric  (squared vertices and bold edges). (Left) Induced link 2-antipodally symmetric admitting two noncentrally symmetric components but antipodally symmetric between them.}
\label{fig16a}
\end{figure}

We say that a link $L(G,S_E)$ is {\em antipodal} if the map $G$ is antipodally self-dual. In the process of our investigations a great deal of structure of antipodal links has been revealed and some problems and questions arose. We notice that there are antipodal links 3-centrally symmetric with a component not necessarily 2-centrally symmetric to itself, see Figure \ref{fig16a}.

\begin{question} Is true that conditions (a) and (b) in Theorem \ref{thm:antipodalR3} are also necessary if a knot is 3-centrally symmetric ?
\end{question} 

\begin{question} Let $(G,S_E^+)$ be an antipodally self-dual edge-signed map. Is it true that alternating link $L(G,S_E^+)$ has always an odd number of components ?
\end{question}

\begin{problem} Characterize all the antipodal links $L(G,S_E)$, in particular, when $G$ is 1-, 2- or 3-connected.
\end{problem}

In view of our results, the following natural question arise.

\begin{question} Let $L$ be an amphichiral link. Does there exists a diagram $D$ of $L$ such that $G_D$ is self-dual ?
\end{question}

This question is on the same flavor as a conjecture due to Kauffman \cite{Kauff} (see also \cite{MR}) claiming that for any alternating amphichiral knot there is a reduced alternating diagram $D$ of the knot such that $G_D$ is self-dual. Dasbach and Hougardy \cite{DH} provided a 14-crossing alternating knot that is a counterexample. Nevertheless, the latter does not seem to be also a counterexample for the above question since it might admits a nonalternating or a no reduced diagram whose Tait's graph is self-dual. 
\smallskip

The conditions given in Theorem \ref{theo:amphi} may lead to an algorithm to detect whether a link is amphichiral. However, its computational complexity depends on how efficiently the desired symmetric cycle is obtained.  

\begin{problem} Let $G$ be a bipartite, planar graph having only square faces. Can it be found either 
a consistent color-reversing symmetric cycle or an opposite color-preserving symmetric cycle in polynomial time ?
\end{problem}

We finally point out that both antipodally self-dual and antipodally symmetric maps play an important r\^ole in {\em real projective} links, that is, links embedded in $\mathbb{R}\mathbb{P}^3$ (work in progress \cite{MRAR3}).   

%%%%%%%%%%%%%%%%%%%%%%%%%%%%%%%%%%%%%%%%%%%%%%%


\begin{thebibliography}{99}

\bibitem{Adam} C.A. Adams, Knot Book : An Elementary Introduction to the Mathematical Theory of Knots, {\em Amer. Math. Soc., Providence , Rhod e Island} (2000).

\bibitem{Cox} H.S.M. Coxeter and W.O.J. Moser, Generators and relations for discrete groups, Ergebnisse des Mathematik und ihrer Grensgebiete, Bd. 14, Springer-Verlag (1972).

%\bibitem{GS} B. Gr\"unbaum and G.C. Shephard, Symmetry Groups of Knots, {\em Mathematics Magazine} {\bf 58}(3) (1985), 161-165.

\bibitem{DH} O. T. Dasbach and S. Hougardy, A conjecture of Kauffman on amphichiral alternating knots, {\em J. Knot Theory Ramifications} {\bf 5}(5) (1996), 629–635.

\bibitem{Flap} E. Flapan, {\em When Topology Meets Chemistry} Cambridge University Press (2000).

\bibitem{HTW} J. Hoste, M. Thistlethwaite and J. Weeks, The first 1, 701, 936 knots, {\em The Mathematical Intelligencer} {\bf 20}(4) (1998) 33--48.

\bibitem{Kauff} L.H. Kauffman, Problems in Knot Theory, in Open problems in topology, 487–522, North-Holland, Amsterdam, 1990.

\bibitem{Liv} C. Livingston, Knot theory, {\em The Carus Math. Monographs} {\bf 24}, Math. Assoc. of Amer. (1993).

\bibitem{MR} J. van Mill, G.M. Reed, Open problems in topology, {\em Topology Appl.} {\bf 42}(3)  (1991), 301–307. 

\bibitem{MRAR1} L. Montejano, J.L. Ram\'irez Alfons\'in and I. Rasskin, Self-dual maps I: antipodality,
arXiv:2008.12853.

\bibitem{MRAR3} L. Montejano, J.L. Ram\'irez Alfons\'in and I. Rasskin, Self-dual maps III: projective links, in preparation.

\bibitem{SS} B. Servatius and H. Servatius, The 24 symmetry pairs of self-dual maps on the sphere, {\em Disc. Math.} {\bf 140} (1995), 167-183.

\end{thebibliography}
\end{document}